\documentclass[a4paper,twoside,11pt,english,intlimits]{article}

\usepackage[utf8]{inputenc}
\usepackage[T2A,T1]{fontenc}
\DeclareSymbolFont{cyrillic}{T2A}{cmr}{m}{n}
\DeclareMathSymbol{\Sha}{\mathalpha}{cyrillic}{216}

\usepackage[english]{babel}

\usepackage{amsmath,amsthm,amssymb,stmaryrd}
\usepackage{mathtools}
\usepackage[ttscale=.875]{libertine}
\usepackage[libertine]{newtxmath}

\usepackage{setspace}

\usepackage{enumerate}
\usepackage{fancyhdr,titlesec,url}
\usepackage[all,knot,poly]{xy}
\usepackage{array}
\usepackage{hyperref}
\usepackage{graphicx}
\usepackage[numbers,sort&compress]{natbib}
\usepackage{multirow}
\usepackage{ifthen}
\usepackage{time}
\usepackage{algorithm2e}

\usepackage{etoolbox}

\usepackage{cancel}
\usepackage{tikz}

\newcommand{\periodafter}[1]{\ifstrempty{#1}{}{#1.}}
\titleformat{\section}[block]{\scshape\filcenter\LARGE\boldmath}{\thesection.}{.5em}{}
\titleformat{\subsection}[block]{\bfseries\filcenter\large\boldmath}{\thesubsection.}{.5em}{\medskip}
\titleformat{\subsubsection}[runin]{\bfseries\boldmath}{\thesubsubsection.}{.5em}{\periodafter}%{}[.]
\titlespacing{\subsubsection}{0pt}{\topsep}{.5em}

\newtheoremstyle{ntheorem}%
	{\topsep}{\topsep}{\itshape}{0pt}{\bfseries}{.}{.5em}%
	{\thmnumber{#2.\hspace{.5em}}\thmname{#1}\thmnote{ (#3)}}
	
\newtheoremstyle{ndefinition}%
	{\topsep}{\topsep}{\normalfont}{0pt}{\bfseries}{.}{.5em}%
	{\thmnumber{#2.\hspace{.5em}}\thmname{#1}\thmnote{ (#3)}}
	
\newtheoremstyle{nremark}%
	{\topsep}{\topsep}{\normalfont}{0pt}{\itshape}{.}{.5em}%
	{\thmnumber{}\thmname{#1}\thmnote{ (#3)}}

\theoremstyle{ntheorem}
  	\newtheorem{theorem}[subsubsection]{Theorem}
  	\newtheorem{proposition}[subsubsection]{Proposition}
	\newtheorem{lemma}[subsubsection]{Lemma}
  	\newtheorem{corollary}[subsubsection]{Corollary}

\theoremstyle{ndefinition}

\makeatletter
\def\@equationname{equation}
\newenvironment{eqn}[1]{%
    \def\mymathenvironmenttouse{#1}%
    \ifx\mymathenvironmenttouse\@equationname%
        \refstepcounter{subsubsection}%
    \else
        \patchcmd{\@arrayparboxrestore}{equation}{subsubsection}{}{}%          doesn't change output?
        \patchcmd{\print@eqnum}{equation}{subsubsection}{}{}%
        \patchcmd{\incr@eqnum}{equation}{subsubsection}{}{}%
    \fi
    \csname\mymathenvironmenttouse\endcsname%
}{%
    \ifx\mymathenvironmenttouse\@equationname%
        \tag{\thesubsubsection}%
    \fi
    \csname end\mymathenvironmenttouse\endcsname%
}
\makeatother

\fancypagestyle{plain}{
\fancyhf{}\fancyfoot[LC,RC]{}
\fancyfoot[LE,RO]{$\thepage$}

}

\UseTips
\SelectTips{eu}{11}

\newdir{ >}{{}*!/-10pt/@{>}}
\newdir{ -}{{}*!/-10pt/@{}}
\newdir{> }{{}*!/+10pt/@{>}}

\makeatletter

\xyletcsnamecsname@{dir4{}}{dir{}}
\xydefcsname@{dir4{-}}{\line@ \quadruple@\xydashh@}
\xydefcsname@{dir4{.}}{\point@ \quadruple@\xydashh@}
\xydefcsname@{dir4{~}}{\squiggle@ \quadruple@\xybsqlh@}
\xydefcsname@{dir4{>}}{\Tttip@}
\xydefcsname@{dir4{<}}{\reverseDirection@\Tttip@}

\xydef@\quadruple@#1{%
	\edef\Drop@@{%
		\dimen@=#1\relax
		\dimen@=.5\dimen@
		\A@=-\sinDirection\dimen@
		\B@=\cosDirection\dimen@
		\setboxz@h{%
			\setbox2=\hbox{\kern3\A@\raise3\B@\copy\z@}%
			\dp2=\z@ \ht2=\z@ \wd2=\z@ \box2
			\setbox2=\hbox{\kern\A@\raise\B@\copy\z@}%
			\dp2=\z@ \ht2=\z@ \wd2=\z@ \box2
			\setbox2=\hbox{\kern-\A@\raise-\B@\copy\z@}%
			\dp2=\z@ \ht2=\z@ \wd2=\z@ \box2
			\setbox2=\hbox{\kern-3\A@\raise-3\B@ \noexpand\boxz@}%
			\dp2=\z@ \ht2=\z@ \wd2=\z@ \box2
		}%
		\ht\z@=\z@ \dp\z@=\z@ \wd\z@=\z@ \noexpand\styledboxz@
	}%
}

\xydef@\Tttip@{\kern2pt \vrule height2pt depth2pt width\z@
	\Tttip@@ \kern2pt \egroup
	\U@c=0pt \D@c=0pt \L@c=0pt \R@c=0pt \Edge@c={\circleEdge}%
	\def\Leftness@{.5}\def\Upness@{.5}%
	\def\Drop@@{\styledboxz@}\def\Connect@@{\straight@{\dottedSpread@\jot}}}
	
\xydef@\Tttip@@{%
	\dimen@=.25\dimen@
 	\B@=\cosDirection\dimen@
	\setboxz@h\bgroup\reverseDirection@\line@ \wdz@=\z@ \ht\z@=\z@ \dp\z@=\z@
	{\vDirection@(1,-1)\xydashl@ \xyatipfont\char\DirectionChar}%
	{\vDirection@(1,+1)\xydashl@ \xybtipfont\char\DirectionChar}%
}

\xydef@\ar@form{
	\ifx \space@\next \expandafter\DN@\space{\xyFN@\ar@form}%
	\else\ifx ^\next \DN@ ^{\xyFN@\ar@style}\edef\arvariant@@{\string^}%
	\else\ifx _\next \DN@ _{\xyFN@\ar@style}\edef\arvariant@@{\string_}%
	\else\ifx 0\next \DN@ 0{\xyFN@\ar@style}\def\arvariant@@{0}%
	\else\ifx 1\next \DN@ 1{\xyFN@\ar@style}\def\arvariant@@{1}%
	\else\ifx 2\next \DN@ 2{\xyFN@\ar@style}\def\arvariant@@{2}%
	\else\ifx 3\next \DN@ 3{\xyFN@\ar@style}\def\arvariant@@{3}%
	\else\ifx 4\next \DN@ 4{\xyFN@\ar@style}\def\arvariant@@{4}%
	\else\ifx \bgroup\next \let\next@=\ar@style
	\else\ifx [\next \DN@[##1]{\ar@modifiers{[##1]}}%]
	\else\ifx *\next \DN@ *{\ar@modifiers}%
	\else\addLT@\ifx\next \let\next@=\ar@slide
	\else\ifx /\next \let\next@=\ar@curveslash
	\else\ifx (\next \let\next@=\ar@curveinout %)
	\else\addRQ@\ifx\next \addRQ@\DN@{\ar@curve@}%
	\else\addLQ@\ifx\next \addLQ@\DN@{\xyFN@\ar@curve}%
	\else\addDASH@\ifx\next \addDASH@\DN@{\defarstem@-\xyFN@\ar@}%
	\else\addEQ@\ifx\next \addEQ@\DN@{\def\arvariant@@{2}\defarstem@-\xyFN@\ar@}%
	\else\addDOT@\ifx\next \addDOT@\DN@{\defarstem@.\xyFN@\ar@}%
	\else\ifx :\next \DN@:{\def\arvariant@@{2}\defarstem@.\xyFN@\ar@}%
	\else\ifx ~\next \DN@~{\defarstem@~\xyFN@\ar@}%
	\else\ifx !\next \DN@!{\dasharstem@\xyFN@\ar@}%
	\else\ifx ?\next \DN@?{\ar@upsidedown\xyFN@\ar@}%
	\else \let\next@=\ar@error
	\fi\fi\fi\fi\fi\fi\fi\fi\fi\fi\fi\fi\fi\fi\fi\fi\fi\fi\fi\fi\fi\fi\fi \next@}

\makeatother

\newcommand{\fl}{\rightarrow}
\newcommand{\fll}{\longrightarrow}
\newcommand{\ofl}[1]{\overset{\displaystyle #1}{\fll}}

\newcommand{\dfl}{\Rightarrow}

\newcommand{\qfl}{\xymatrix@1@C=10pt{\ar@4 [r] &}}

\newcommand{\ens}[1]{\left\{ #1 \right\}}

\newcommand{\op}[1]{{#1}^{o}}

\newcommand{\cl}[1]{\overline{#1}}
\newcommand{\rep}[1]{\widehat{#1}}

\DeclareMathOperator{\id}{Id}

\DeclareMathOperator{\Hom}{\mathrm{Hom}}

\DeclareMathOperator{\Tor}{\mathrm{Tor}}

\renewcommand{\phi}{\varphi}
\renewcommand{\epsilon}{\varepsilon}

\newcommand{\Nb}{\mathbb{N}}

\newcommand{\dr}{\partial}

\newcommand{\Br}{\mathcal{B}}
\newcommand{\Cr}{\mathcal{C}}

\newcommand{\Fr}{\mathcal{F}}
\newcommand{\Gr}{\mathcal{G}}

\newcommand{\Lr}{\mathcal{L}}

\newcommand{\Ur}{\mathcal{U}}
\newcommand{\Vr}{\mathcal{V}}
\renewcommand{\Wr}{\mathcal{W}}

\newcommand{\K}{\mathbb{K}}

\def\catego#1{\mathsf{#1}}

\newcommand{\Cat}{\catego{Cat}}
\newcommand{\Pol}{\catego{Pol}}

\newcommand{\Set}{\catego{Set}}
\newcommand{\Ord}{\mathbf{Ord}}

\newcommand{\Vect}{\catego{Vect}}
\newcommand{\grVect}{\catego{grVect}}

\newcommand{\ho}{\mathrm{H}}

\newcommand{\ifthen}[2]{\ifthenelse{#1}{#2}{}}

\renewcommand{\id}{id}

\def\Cat{\catego{Cat}}
\def\Gpd{\catego{Gpd}}
\def\glob{\mathbb{O}}
\def\Rglob{\mathbb{RO}}
\def\Glob{\catego{Glob}}
\def\RGlob{\catego{RGlob}}
\newcommand\Bimod[1]{\catego{Bimod}(#1)}

\newcommand{\Item}[1]{{\upshape\textbf{(#1)}}}
\DeclareMathOperator{\Red}{Red}
\DeclareMathOperator{\red}{Red_m}
\DeclareMathOperator{\Sph}{Sph}

\newcommand{\undermath}[2]{\underset{\mathclap{#1}}{#2}}
\def\part{\odot}
\def\SS{\mathfrak{S}}

\newcommand{\triNormal}[1]{\begin{tikzpicture}[scale=#1]
\draw[line cap=round]
  (-1ex,.5ex) -- (0,0)
  (-1ex,.5ex) -- (0,1ex)
  (-1ex,.5ex) -- (0,.33ex)
  (-1ex,.5ex) -- (0,.67ex);
\end{tikzpicture}}
\def\tri{\mathbin{\mathchoice
{\triNormal{1.1}}
{\triNormal{1.1}}
{\triNormal{.8}\,}
{\triNormal{.6}}}}
\def\Unit{\mathbb{I}}

\newcommand{\range}[1]{\{1,\ldots,#1\}}

\usepackage{forest}
\forestset{
	default preamble={
		begin draw/.code={\begin{tikzpicture}[baseline={([yshift=-.5ex]current bounding box.center)}]}, 
		for tree={grow'=north, math content}
	},
	skeleton/.style={
		for tree={l=0, l sep=1, s=-2, s sep=-2, if n children=0
			{child anchor=south, node options={font=\scriptsize}}
			{parent anchor=south, child anchor=south}}
	},
	input tree/.style={
		for tree={edge=-, l=0, l sep=4, if n children=0
			{node options={font=\scriptsize}, inner ysep=1}
			{node options={inner ysep=1}}}
	},
	no input tree/.style={
		for tree={edge=-, l=0, l sep=4,
		node options={inner ysep=1}}
	},
	big input tree/.style={
		for tree={edge=-, if n children=0
			{node options={font=\scriptsize}, inner ysep=1}
			{node options={inner ysep=1}}}
	}
}
\newcommand{\skele}[1]{
	\begin{forest}
		skeleton,
		#1
	\end{forest}}
\newcommand{\inputtree}[1]{
	\begin{forest}
		input tree,
		#1
	\end{forest}}
\newcommand{\tree}[1]{
	\begin{forest}
		no input tree,
		#1
	\end{forest}}

\newcommand{\segment}[2]
{%
\tree{[#1[#2, xshift=2.5ex]]}%
}

\newcommand{\pullbackcorner}[1][dr]{\save*!/#1-1.6pc/#1:(-1,1)@^{|-}\restore}

\def\Ov{\mathcal{O}v}

\DeclareMathOperator{\supp}{Supp}
\DeclareMathOperator{\Std}{Std}

\DeclareMathOperator{\lm}{lm}
\DeclareMathOperator{\lc}{lc}
\DeclareMathOperator{\lt}{lt}

\renewcommand{\leq}{\leqslant}
\renewcommand{\geq}{\geqslant}

\newcount\hh
\newcount\mm
\mm=\time
\hh=\time
\divide\hh by 60
\divide\mm by 60
\multiply\mm by 60
\mm=-\mm
\advance\mm by \time
\def\hhmm{\number\hh:\ifnum\mm<10{}0\fi\number\mm}

\def\Abel{\catego{Ab}}
\def\K{\mathbf{k}}

\def\SymO{\SS\catego{Op}}
\def\ShuO{\Sha\catego{Op}}
\def\grShuO{\catego{gr}\Sha\catego{Op}}
\def\ShuT{\Sha\catego{Tree}}
\def\Coll{\catego{Coll}}
\def\Ind{\catego{Ind}}
\def\Fin{\catego{Fin}}
\def\Ord{\catego{Ord}}
\def\SSColl{\SS\Coll} 
\def\ShuPol{\Sha\Pol}

\def\setShuffleComp{\times_\Sha}
\newcommand{\fopShu}[1]{#1^{\Sha}}
 
\newcommand{\ftree}[1]{#1^{\ast_\Sha}}

\newcommand{\auteur}[3]{
\noindent
\begin{minipage}[t]{.45\textwidth}
\begin{flushright}
\textsc{#1} \\
{\footnotesize\textsf{#2}}
\end{flushright} 
\end{minipage}
\qquad
\begin{minipage}[t]{.45\textwidth}
#3
\end{minipage}
}

\begin{document}
\thispagestyle{empty}

\begin{center}

\begin{doublespace}
\begin{huge}
{\scshape Shuffle polygraphic resolutions}
\end{huge}

\vskip+2pt

\begin{huge}
{\scshape for operads}
\end{huge}

\bigskip
\hrule height 1.5pt 
\bigskip

\vskip+5pt

\begin{Large}
{\scshape Philippe Malbos - Isaac Ren}
\end{Large}
\end{doublespace}

\vfill

\vskip+20pt

\begin{small}\begin{minipage}{14cm}
\noindent\textbf{Abstract --}
Shuffle operads were introduced to forget the symmetric group actions on symmetric operads while preserving all possible operadic compositions.
Rewriting methods were then applied to symmetric operads via shuffle operads: in particular, a notion of Gröbner basis was introduced for shuffle operads with respect to a total order on tree monomials.
In this article, we introduce the structure of shuffle polygraphs as a categorical model for rewriting in shuffle operads, which generalizes the Gröbner bases approach by removing the constraint of a monomial order for the orientation of the rewriting rules.
We define $\omega$-operads as internal $\omega$-categories in the category of shuffle operads. 
We show how to extend a convergent shuffle polygraph into a shuffle polygraphic resolution generated by the overlapping branchings of the original polygraph.
Finally, we prove that a shuffle operad presented by a quadratic convergent shuffle polygraph is Koszul.
\medskip


\smallskip\noindent\textbf{M.S.C. 2020 --} 18M70, 68Q42, 18N30.
\end{minipage}
\end{small}
\end{center}

\vskip+0pt

\begin{center}
\begin{small}\begin{minipage}{12cm}
\renewcommand{\contentsname}{}
\setcounter{tocdepth}{2}
\tableofcontents
\end{minipage}
\end{small}
\end{center}

\clearpage

\section{Introduction}

Algebraic rewriting theory provides methods to compute cofibrant replacements of algebraic structures from presentations that take into account computational properties of these structures. This rewriting approach gives algebraic algorithmic methods to solve decidability and computational problems, such as the ideal membership problem, and the computation of linear bases and of (co)homological properties. Abelian resolutions for monoids~\cite{Squier87,Kobayashi90, Brown92}, groups \cite{DehornoyLafont03}, small categories~\cite{GuiraudMalbos12advances}, associative algebras~\cite{Anick86, GuiraudHoffbeckMalbos19}, and linear operads \cite{DotsenkoKhoroshkin10, DotsenkoKhoroshkin13} have been constructed using rewriting methods.
The machinery at the heart of these constructions consists in presenting an algebraic structure by a system of \emph{generators} and \emph{rewriting rules}, and producing a cofibrant replacement that involves the overlappings occurring in the applications of the rewriting rules. 
Rewriting approaches for linear structures were developed in many algebraic algorithmic contexts, notably by Janet and Buchberger for commutative algebras~\cite{Janet20,Buchberger65}, Shirshov, Bokut, and Bergman for associative algebras~\cite{Shirshov62,Bokut76,Bergman78}, Dotsenko-Khoroshkin for linear operads~\cite{DotsenkoKhoroshkin10,KhoroshkinPiontkovski2015}. In all of these works, the rewriting systems are formulated in terms of \emph{Gröbner bases}, and thus are defined with respect to a given monomial order.
Rewriting approaches have also been used in the categorical context to present higher categories by higher-dimensional rewriting systems, called \emph{polygraphs} (or \emph{computads}) \cite{Street76,Burroni93}. In this context, the cofibrant replacements of a higher category are generated by \emph{polygraphic resolutions} introduced in~\cite{Metayer03,LafontMetayerWorytkiewicz10,GuiraudMalbos12advances}.

An important issue when studying algebras or operads is the automatic construction of small abelian resolutions.
There exist some inductive constructions that start from presentations with certain computational properties.
In particular, Anick introduced a general machinery that computes a resolution for an associative algebra, whose $n$-dimensional generators correspond to overlappings of applications of $n$ defining relations. This resolution is a compromise between the bar resolution, which is easy to compute but very large, and the minimal resolution, which is difficult to make explicit in general. Moreover, it is difficult to determine conditions for which Anick's resolution is minimal~\cite{Tamaroff21}. 
Some conditions have been shown to be sufficient: as an immediate consequence of its construction, Anick's resolution is minimal for monomial algebras, and a quadratic algebra with a convergent presentation has a minimal resolution given by its Koszul dual. For $N$-homogeneous algebras, minimality is harder, and Berger introduces an extra condition in~\cite{Berger01}.
Dotsenko and Khoroshkin constructed resolutions for shuffle monomial operads by the inclusion-exclusion principle and for operads presented by a Gröbner basis by deformation of the monomial case~\cite{DotsenkoKhoroshkin13}. In this operadic case, the question of minimality of the resolution is even more difficult due to the combinatorial complexity of the underlying tree structure of shuffle operad terms. In particular, unlike for algebras the constructed resolution for a monomial operad is not necessarily minimal.

In this work, we combine the polygraphic and the Gröbner bases approaches in order to compute higher-dimensional presentations of shuffle operads using the polygraphic machinery.
We define \emph{shuffle $\omega$-operads} as internal $\omega$-categories in the category of shuffle operads. We introduce the structure of \emph{shuffle $\omega$-polygraphs} as systems of generators and relations for shuffle $\omega$-operads. Unlike the Gröbner bases approach, the orientation of the relations in a shuffle polygraph does not depend on a given monomial order. 
The main construction of this article extends a confluent and terminating shuffle polygraph presenting a shuffle operad into a shuffle polygraphic resolution generated by the overlapping branchings of the original polygraph.
In order to address the question of minimal resolutions, we make explicit these overlappings in all dimensions of the polygraphic resolution.
We then give an inductive method to compute a bimodule resolution that allows us to state a minimality result for shuffle operads, as well as a condition for Koszulness.

Now we present the organization and the main results of this article.

\subsection*{Higher operads}

The notion of a \emph{symmetric operad} appears in many situations to describe operations in several arguments, with symmetric group actions, acting on topological or algebraic objects~\cite{May72, Loday96}. \emph{Shuffle operads} were introduced by Dotsenko and Khoroshkin in \cite{DotsenkoKhoroshkin10} to forget the symmetric group actions on the arguments while preserving all possible operadic compositions. The shuffle version allows us to define \emph{monomials} and \emph{oriented relations} in order to present symmetric operads by rewriting systems. 
Symmetric and shuffle operads are defined as internal monoids in the monoidal presheaf categories of \emph{collections} and \emph{symmetric collections} respectively, as recalled in Section~\ref{SS:ShuffleOperads}. Explicitly, a \emph{collection} is a presheaf on the category $\Ord$ of nonempty finite ordered sets and order-preserving bijections, with values in the category $\Vect$ of vector spaces. The monoidal product on collections is the \emph{shuffle composition} recalled in \textsection~\ref{SSS:MonoidalProducts}. A \emph{symmetric collection} is a presheaf on the category $\Fin$ of nonempty finite sets and bijections, with values in the category $\Vect$. The functor $-^u:\Ord\rightarrow\Fin$ that forgets the order induces a functor $-^u:\SSColl\rightarrow\Coll$ from the category of symmetric collections to the category of collections. Since we restrict to nonempty sets, none of the operads considered in this work have operations of arity $0$.

In Section~\ref{SS:HigherDimensionalShuffleOperads}, we introduce the notion of a (strict) higher shuffle operad. We define a \emph{shuffle $\omega$-operad} as an internal $\omega$-category in the category $\ShuO$ of shuffle operads. Shuffle $\omega$-operads, with internal $\omega$-functors, form a category denoted by $\ShuO_\omega$.
In Section~\ref{S:HigherDiomensionalOperadsAsBimobules}, we study the interaction between the higher-categorical structure of $\omega$-operads and its underlying linear structure. The object of $n$-cells of a shuffle $\omega$-operad has a shuffle operad structure, and the $n$-cells can be \emph{$\star_k$-composed} along $k$-dimensional cells for $0\leq k<n$. Due to the linear structure, the $\star_k$-composition of two $n$-cells $a$ and~$b$ in a shuffle $\omega$-operad can be written as the following linear combination:
\[
a\star_k b=a-t_k(a)+b,
\]
where $t_k(a)$ denotes the $k$-dimensional target of $a$, which coincides with the $k$-dimensional source of~$b$. In particular, every $n$-cell in a shuffle $\omega$-operad is invertible. Moreover, for $n\geq 1$, the compatibility between the shuffle composition and the $\star_0$-composition implies that the \emph{elementary composition} $a\circ_{i,\tau}b$ of $n$-cells $a$ and $b$, as defined in \textsection~\ref{SSS:ShuffleCompositions}, can be seen either one of two \emph{orthogonal reduction paths} from $s_0(a)\circ_{i,\tau}s_0(b)$ to $t_0(a)\circ_{i,\tau}t_0(b)$, pictured as follows:

\vskip-10pt

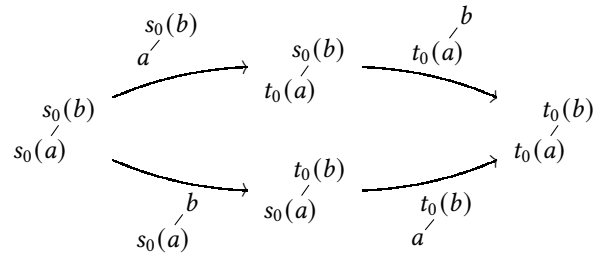
\[
\begin{small}
\xymatrix@C=5em @R=-1em{
&
\segment{t_0(a)}{s_0(b)}
  \ar@/^2ex/ [dr] ^-{\segment{t_0(a)}{b}}
&
\\
\segment{s_0(a)}{s_0(b)}
  \ar@/^2ex/ [ur] ^-{\segment{a}{s_0(b)}}
  \ar@/_2ex/ [dr] _-{\segment{s_0(a)}{b}}
&& 
\segment{t_0(a)}{t_0(b)}
\\
&
\segment{s_0(a)}{t_0(b)}
  \ar@/_2ex/ [ur] _-{\segment{a}{t_0(b)}}
&
}
\end{small}
\]
The \emph{linear exchange relation} introduced in \textsection~\ref{SSS:LinearExchangeRelation} states that these two reductions paths are equal.
With these remarkable relations, the axioms of shuffle $\omega$-operads can be simplified. We deduce a characterization of the structure of $\omega$-operad in terms of bimodules over shuffle operads. Our first result, Theorem~\ref{T:MainTheoremA}, proves that the category $\ShuO_\omega$ is isomorphic to the full subcategory of $\RGlob(\Bimod{\ShuO})$, whose objects are pairs $(P,A)$ where $P$ is a shuffle operad and $A=(A_n)_{n\geq 0}$ is a reflexive globular $P$-bimodule such that $A_0=P$, and $A_n$ satisfies the linear exchange relation for all $n\geq 1$.
 
\subsection*{Shuffle operadic polygraphs and rewriting}

The notion of a \emph{polygraph} was introduced in the set-theoretical context by Street and Burroni as systems of generators and relations for presentations of higher (strict) categories \cite{Street76, Burroni91}. A linear version of polygraphs was introduced in \cite{GuiraudHoffbeckMalbos19} for the presentation of associative $\omega$-algebras. In Section~\ref{SS:ShufflePolygraphs}, we define an analogous notion for shuffle $\omega$-operads, which we call \emph{shuffle polygraphs}. Explicitly, for $n\geq 0$, a \emph{shuffle $n$-polygraph} is a data $X=(X_0,\ldots X_n)$ defined by induction, where~$X_k$, the set of \emph{$k$-generators}, forms a globular extension of the free shuffle $(k-1)$-operad generated by the shuffle $(k-1)$-polygraph $(X_0,\ldots X_{k-1})$. Such a data can be pictured as a diagram
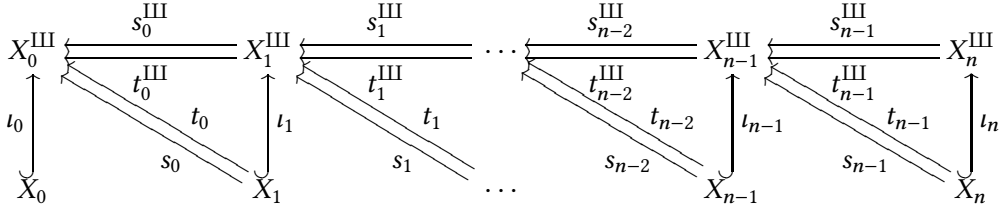
\[
\xymatrix @C=2.25cm @R=1.25cm{
\fopShu{X_0}
&
\fopShu{X_1}
  \ar@<-0.5ex>[l] _-{s_0^{\Sha}}
  \ar@<+0.5ex>[l] ^-{t_0^{\Sha}}
&
{\vphantom{\fopShu{X_2}}\cdots}
  \ar@<-0.5ex>[l] _-{s_1^{\Sha}}
  \ar@<+0.5ex>[l] ^-{t_1^{\Sha}}
&
\fopShu{X_{n-1}}
  \ar@<-0.5ex>[l] _-{s_{n-2}^{\Sha}}
  \ar@<+0.5ex>[l] ^-{t_{n-2}^{\Sha}}
&
\fopShu{X_n}
  \ar@<-0.5ex>[l] _-{s_{n-1}^{\Sha}}
  \ar@<+0.5ex>[l] ^-{t_{n-1}^{\Sha}}
\\
X_0
  \ar@{^(->}[u] ^-{\iota_0}%)
&
X_1
  \ar@<+0.5ex> [ul] ^(0.35){s_0}
  \ar@<-0.5ex> [ul] _(0.35){t_0}
  \ar@{^(->}[u] _-{\iota_1}
&
{\vphantom{X_2}\cdots}
  \ar@<+0.5ex> [ul] ^(0.35){s_1}
  \ar@<-0.5ex> [ul] _(0.35){t_1}
&
X_{n-1}
  \ar@<+0.5ex> [ul] ^(0.35){s_{n-2}}
  \ar@<-0.5ex> [ul] _(0.35){t_{n-2}}
  \ar@{^(->}[u] _-{\iota_{n-1}}
&
X_n
  \ar@<+0.5ex> [ul] ^(0.35){s_{n-1}}
  \ar@<-0.5ex> [ul] _(0.35){t_{n-1}}
  \ar@{^(->}[u] _-{\iota_n}
}
\]
where $s_i$ and $t_i$ denote the source and target maps of the globular extensions, and the horizontal diagram corresponds to the underlying globular operad of the free $n$-operad generated by the $n$-polygraph~$X$, denoted by $\fopShu{X}_n$. As for set-theoretical polygraphs, in Section~\ref{SS:ShufflePolygraphs} we define the category $\ShuPol_n$ of \emph{shuffle $n$-polygraphs} and the free $n$-operad functor $\fopShu{(-)} : \ShuPol_n \fl \ShuO_n$ by induction on the dimension $n$, and the category of shuffle $\omega$-polygraphs as the limit of the forgetful functors $\ShuPol_n \fl  \ShuPol_{n-1}$ for $n\geq 1$.

The shuffle polygraphic approach lets us present shuffle operads by oriented presentations, called \emph{rewriting systems}: the shuffle operad $\cl{X}$ presented by a shuffle $1$-polygraph $X$ is defined as the coequalizer of the source and target morphisms $s_0^\Sha,t_0^\Sha : \fopShu{X}_1 \rightrightarrows \fopShu{X_{0}}$ in the category $\ShuO$.
Note that, in this work, we consider operads with only one color. The $0$-generators correspond to the generators of the shuffle operads, and the $1$-generators correspond to the oriented relations. For presentations of multicolored shuffle operads, we need to consider shuffle two-dimensional polygraphs, whose $0$-generators correspond to colors, $1$-generators to generators, and $2$-generators to oriented relations.

We use rewriting theory on shuffle $1$-polygraphs to deduce global rewriting properties, such as \emph{confluence} and \emph{termination} from local properties of the $1$-generators, also called \emph{rewriting rules}. Without loss of generality, in Section~\ref{S:ShuffleOperadicRewriting} we will consider $1$-polygraphs with \emph{left-monomial} rules, reducing a single monomial into a linear combination of monomials. A \emph{rewriting step} of a left-monomial $1$-polygraph $X$ is a $1$-cell $f$ of the free shuffle $1$-operad $\fopShu{X}_1$ of size $1$, and of the form $f = \lambda g + 1_c$, where $\lambda$ is a nonzero scalar, $g$ is a $1$-monomial of $\fopShu{X}_1$, and $c$ is a $0$-cell of the free shuffle operad $\fopShu{X}_0$ such that the $0$-monomial $s_0(u)\notin \supp(c)$. A $1$-cell of the free $1$-operad $\fopShu{X}_1$ is \emph{positive} if it is the $\star_0$-composition of rewriting steps. A polygraph is \emph{terminating} if there is no infinite sequence of $\star_0$-composition of rewriting steps.

This shuffle polygraphic approach generalizes that of Gröbner bases introduced by Dotsenko and Khoroshkin in \cite{DotsenkoKhoroshkin10}. Indeed, the orientation of the polygraphic rules does not depend on a given monomial order. However, termination is not ensured by a monomial order, so it must be proven by considering the rewriting rules themselves. Beyond the property of termination, the \emph{confluence} property of a $1$-polygraph $X$ states that for every \emph{branching} of two positive $1$-cells $f$, $g$ of $\fopShu{X}_1$ with the same source~$a$, there exist two positive $1$-cells $h$ and $k$ of $\fopShu{X}_1$ as in the following \emph{confluence diagram}:
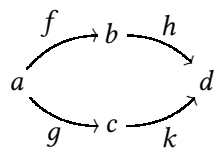
\[
\vcenter{\xymatrix @C=2.2em @R=0.35em{
& b
	\ar@/^/ [dr] ^-{h}
\\
a
	\ar@/^/ [ur] ^-{f}
	\ar@/_/ [dr] _-{g}
&&
d
\\
&
c
	\ar@/_/ [ur] _-{k}
}}
\]

When the system is terminating, confluence can be deduced from \emph{local confluence}, that is, when all the branchings of rewriting steps are confluent~\cite{Newman42, Huet80}. Local confluence can be proven by the confluence of all branchings involving minimal overlappings of the rules, called \emph{critical branchings}. This is the \emph{critical branching theorem} proved in many algebraic contexts \cite{Nivat73,KnuthBendix70, ChenavierDupontMalbos20}. Coherent versions of this result where introduced in~\cite{GaussentGuiraudMalbos15,GuiraudHoffbeckMalbos19}. 
In \textsection~\ref{SSS:ClassificationLocalBranchings} we introduce \emph{essential branchings} that refine the notion of a critical branching, and generate all critical branchings by transitivity.
Theorem~\ref{T:CriticalBranchings} proves a coherent essential branching theorem for shuffle polygraphs. It states that we can extend a terminating, convergent, left-monomial shuffle $1$-polygraph $X$ into an \emph{acyclic shuffle $2$-polygraph} by considering confluences of essential branchings. 

In Section~\ref{SS:MonomialOrder}, we give several algebraic interpretations of the confluence property of a terminating, left-monomial shuffle $1$-polygraph $X$. Proposition~\ref{P:ConfluenceAlgebraicCharacterisation} proves that the confluence of $X$ is equivalent to having a decomposition of the free shuffle operad $\fopShu{X}_0$ into a direct sum of the ideal generated by the $1$-generators of $X$ and the collection of normal forms with respect to this $1$-generators. Proposition~\ref{P:ConvergentGröbner} proves that the notion of a convergent shuffle $1$-polygraph, where the rules are oriented with respect to a given monomial order, is equivalent to the notion of Gröbner bases introduced in \cite{DotsenkoKhoroshkin10}. Proposition~\ref{P:ConvergentPBW} gives a polygraphic interpretation of the Poincaré-Birkhoff-Witt (PBW) criterion introduced by Hoffbeck in  \cite{Hoffbeck10} as a generalization of Priddy's PBW criterion for associative algebras~\cite{Priddy70}.

\subsection*{The overlapping polygraphic resolution and Koszulness}

An $\omega$-polygraph $X$ is \emph{acyclic} if, for every $n\geq 1$, the quotient of the free shuffle $n$-operad~$\fopShu{X}_n$ by the ideal generated by the cellular extension $X_{n+1}$ is aspherical, that is all parallel $n$-cells are equal. We say that an acyclic $\omega$-polygraph $X$ is a \emph{shuffle polygraphic resolution} of the shuffle operad it presents.
Section~\ref{S:ShufflePolygraphicResolutionConvergence} presents the main result of this article, Theorem~\ref{T:MainTheoremD}, which extends a reduced convergent left-monomial shuffle $1$-polygraph $X$ into a polygraphic resolution, denoted by $\Ov(X)$, of the presented shuffle operad. The generators of this polygraphic resolution correspond to higher-dimensional overlappings induced by the rewriting rules of $X$, defined in \textsection~\ref{SSS:HigherDimensionalOverlappings}: an $n$-generator of $\Ov(X)$, called an \emph{$n$-overlapping}, is a sequence of monomials, written
\[
u_0\tri \vec v_1 \tri \vec v_2 \tri \cdots \tri \vec v_n,
\]
where, when seen as planar trees, each sequence of monomials $\vec v_i = (v_{i,1}, \ldots, v_{i,k})$ is attached to the leaves of the $(i-1)$-overlapping $u_0\tri\vec v_1\tri\cdots\tri\vec v_{i-1}$ in a manner that adds exactly enough to apply a new rewriting rule. Explicitly, for low dimensions, the $0$-overlappings correspond to the $0$-generators~$X_0$, the $1$-overlappings correspond to the sources of the rewriting rules of $X_1$, and the $2$-overlappings correspond to critical branchings of $X$.

The acyclicity of the shuffle $\omega$-polygraph $\Ov(X)$ is proven by the construction of a \emph{homotopical contraction}.
Thus the shuffle polygraphic resolution $\Ov(X)$ provides an alternative construction to the differential-graded shuffle operads constructed by Dotsenko and Khoroshkin \cite{DotsenkoKhoroshkin13}. Moreover, since associative algebras are particular cases of shuffle operads, the construction of $\Ov(X)$ is another way to obtain the polygraphic resolutions for associative algebras introduced in \cite{GuiraudHoffbeckMalbos19}.

The Quillen homology of a symmetric operad can be computed from its associated shuffle operad. Indeed, the reduced bar complex of a symmetric operad $P$ is isomorphic, as a shuffle differential-graded cooperad, to the reduced bar complex of the associated shuffle operad $P^u$, that is, $\overline B(P)^u \simeq \overline B(P^u)$, {\cite[Prop. 1.4]{DotsenkoKhoroshkin13}}. 
In Section~\ref{SS:BimodulesResolution}, we associate to an acyclic shuffle $\omega$-polygraph~$X$ that presents a shuffle operad~$P$ a $P$-bimodule resolution $P\langle X\rangle$ of the \emph{trivial $P$-bimodule}~$\Omega_P$.
Thus shuffle polygraphic resolutions provide a constructive way to compute the homology of symmetric operads.
Proposition~\ref{P:concentratedDiagonal} proves that when there exists an increasing function $w:\Nb\fl\Nb\setminus\{0\}$ such that the $n$-generators of $X$ are concentrated in weight~$w(n)$, then the resolution $P\langle X\rangle$ is minimal.
Finally, we define in Section~\ref{SS:ConfluenceKoszul} a criterion of Koszulness in terms of quadratic convergence: Theorem~\ref{T:MainTheoremE} states that shuffle operads presented by quadratic convergent $1$-polygraphs are Koszul. This result generalizes those obtained by Dotsenko and Khoroshkin in \cite{DotsenkoKhoroshkin10} for shuffle operads with quadratic Gröbner bases defined with respect to a given monomial order. This new rewriting-based sufficient condition for Koszulness does not depend on a monomial order, which is required to define Gröbner bases. 

\subsection*{Conventions and notations}

Throughout this article $\K$ denotes a field of characteristic zero. All vector spaces are over this field~$\K$, and we denote by $\Vect$ the category of vector spaces and linear maps as morphisms.
We denote by~$\Nb$ the set of nonnegative integers. We denote by~$\Ord$ the category of nonempty finite ordered sets, whose morphisms are order-preserving bijections. We denote by $\Fin$ the category of nonempty finite sets, whose morphisms are bijections.

\section{Higher shuffle operads}
\label{S:HigherDimensionalShuffleOperads}

In this section we introduce the notion of a higher shuffle operad. We first recall the structure of shuffle operads from \cite{DotsenkoKhoroshkin10} and we decompose the shuffle composition into elementary compositions. We then define the category of shuffle $\omega$-operads and characterize it as a certain subcategory of globular bimodules over shuffle operads.

\subsection{Shuffle operads}
\label{SS:ShuffleOperads}

In this preliminary subsection we recall from \cite{DotsenkoKhoroshkin10} the definitions on shuffle operads used in this article. We refer the reader to \cite{LodayVallette12} or \cite{BremnerDotsenko16} for a complete account on symmetric and shuffle operads.

\subsubsection{Presheaves on finite sets}
\label{SSS:Modules}
In general, a presheaf $X$ on $\Ord$ or $\Fin$ with values in a category~$\Cr$ are determined by the family of objects $(X(k))_{k\geq 1}$, where $X(k):=X(\range{k})$. We will adopt this notation for the following definitions.

An \emph{set indexed by the category $\Ord$}, or \emph{indexed set} for short, is a presheaf on $\Ord$ with values in the category $\Set$. We denote by $\Ind$ the category of indexed sets with natural transformations as morphisms.
A \emph{collection} is a presheaf on $\Ord$ with values in the category $\Vect$. We denote by $\Coll$ the category of collections and their natural transformations. 
A \emph{basis} of $V$ is a indexed set $B=(B(k))_{k\geq 1}$ such that, for each $k$, $B(k)$ is a basis of the space $V(k)$.

A \emph{symmetric collection} is a presheaf on $\Fin$ with values in the category $\Vect$. We denote by $\SSColl$ the category of symmetric collections and their natural transformations. The functor $-^u:\Ord\rightarrow\Fin$ that forgets the order induces a functor
\begin{eqn}{equation}
\label{E:forgetfulFunctorPetitU}
-^u:\SSColl\rightarrow\Coll.
\end{eqn}
In addition, we denote by $\K\langle-\rangle:\Ind\rightarrow\Coll$ the left adjoint functor of the forgetful functor $\Coll \fl \Ind$. 

\subsubsection{Operads (\cite{DotsenkoKhoroshkin10})}
\label{SSS:MonoidalProducts}
The categories $\Coll$ and $\SSColl$ are monoidal with the following products:
\begin{enumerate}
\item the \emph{shuffle composition} on $\Coll$ denoted by $\circ_\Sha$, and defined for $V,W\in\Coll$ by
\[
(V \circ_\Sha W)(I) := \bigoplus_{k = 1}^\infty V(k)\otimes \left(\bigoplus_{\substack{f : I \twoheadrightarrow \range{k} \\ (*)}} W(f^{-1}\{1\}) \otimes \cdots \otimes W(f^{-1}\{k\})\right),
\]
where $I\in\op{\mathbf{Ord}}$ and the sum $(*)$ is taken on \emph{shuffle surjections}, \emph{i.e.}, surjections $f:I\twoheadrightarrow \range{k}$ such that
\[
\min f^{-1}\{1\}<\ldots<\min f^{-1}\{k\}.
\]
\item the \emph{symmetric composition} on $\SSColl$, denoted by $\circ_\SS$, is defined for $V,W\in\SSColl$ by
\[
(V \circ_\SS W)(I) := \bigoplus_{k = 1}^\infty V(k) \otimes_{\K[\SS_k]}\left(\bigoplus_{f : I \twoheadrightarrow \range{k}} W(f^{-1}\{1\}) \otimes \cdots \otimes W(f^{-1}\{k\})\right),
\]
where $I\in\op{\mathbf{Fin}}$, and the sum is taken on all surjections. 
\end{enumerate}
In both case, the unit is the collection $\Unit$ concentrated in arity $1$ with $\Unit(1)=\K$.
	
A \emph{shuffle} (resp.\ \emph{symmetric}) \emph{operad} is an internal monoid $(P,\mu_P,\eta_P)$ in $(\Coll,\circ_\Sha,\Unit)$ (resp.\ $(\SSColl,\circ_\SS,\Unit)$), where $\mu_P$ is the multiplication morphism and $\eta_P$ is the unit morphism. We denote respectively by $\SymO$ and $\ShuO$ the category of symmetric operads and shuffle operads and their morphisms.

The free operad functors $\ftree{-}:\Coll\rightarrow\ShuO$ and $-^{\ast_\SS}:\SSColl\rightarrow\SymO$ are defined using the free monoid functor on left distributive categories as detailed in \cite[Appendix B]{BauesJibladzeTonks97}. 
For an indexed set $X$, we denote by $\fopShu{X}$ the free shuffle operad on $X$ given by the composite of free functor
\begin{eqn}{equation}
\label{E:FreeFunctorOperadFromSet}
\xymatrix{
\Ind
\ar[r] ^-{\K}
&
\Coll 
\ar[r] ^-{\ftree{-}}
&
\ShuO.
}
\end{eqn}

Recall from \cite{DotsenkoKhoroshkin10}, see also \cite{BremnerDotsenko16}, that the forgetful functor $-^u$ is monoidal in the sense that for all symmetric collections $V,W$, we have $(V\circ_\SS W)^u=V^u\circ_\Sha W^u$ in $\Coll$, and in particular that it commutes with free operad functors $-^{*_\SS}$ and $\ftree{-}$, in the sense that for every symmetric collection $V$, we have the isomorphism
\begin{eqn}{equation}
\label{E:CommutationFunctorUwithFreeFunctors}
(V^{*_\SS})^u=\ftree{(V^u)}.
\end{eqn}

\subsubsection{Shuffle composition on indexed sets}
\label{SSS:ShuffleCompositionOmegaSets}
We define a monoidal \emph{shuffle composition} on $\Ind$, also denoted by $\setShuffleComp$, by setting, for indexed sets $X,Y$
\[
(X \setShuffleComp Y)(I) := \coprod_{k = 1}^\infty X(k)\times \left(\coprod_{\substack{f : I \twoheadrightarrow \range{k} \\ (*)}} Y(f^{-1}\{1\}) \times \cdots \times Y(f^{-1}\{k\})\right)
\]
where the coproduct $(*)$ is taken on shuffle surjections. The composition $\setShuffleComp$ has for unit the indexed set concentrated in arity $1$, denoted by $\mathbf{1}$, and such that $\mathbf 1(1)$ is a singleton, whose only element is denoted by~$\varepsilon$.
The functor $\K$ is compatible with product and coproduct, hence the following diagram commutes:
\begin{eqn}{equation}
\label{E:ShuffleCompositionbCommute}
\raisebox{0.6cm}{
\xymatrix@R=1.6em @C=4em
{
\Ind\times\Ind
  \ar[d]_-{\setShuffleComp}
  \ar[r]^-{\K\times\K}
&
\Coll\times\Coll
  \ar[d]^-{\circ_\Sha}
\\
\Ind
  \ar[r]_-{\K}
&
\Coll
}}
\end{eqn}

\noindent
Note that the adjunction between the monoidal categories $(\Set,\times,\{*\})$ and $(\Vect,\otimes,\K)$ is compatible with the canonical isomorphisms of units, associativity, and distributivity, so the induced functors between $(\Ind,\setShuffleComp,\Unit)$ and $(\Coll,\circ_\Sha,\Unit)$ make a lax monoidal adjunction.

\subsubsection{Tree monomials}
The shuffle composition is monoidal, and we denote by $\ShuT$ the category of internal monoids in $\Ind$ with respect to this composition.
The functor $\K$ preserves colimits as a left adjoint, and sends $\setShuffleComp$ to $\circ_\Sha$ as a consequence of commutativity of (\ref{E:ShuffleCompositionbCommute}). 
Free internal monoids in $\Ind$ and $\Coll$ are constructed by colimits and shuffle composition, thus the linearization functor $\K$ induces a linearization functor $\ShuT \fl \ShuO$ such that the following square commutes:
\[
\xymatrix@R=1.6em @C=4em{
\Ind
  \ar[d]_-{\ftree{-}}
  \ar[r]^-{\K}
&
\Coll
  \ar[d]^-{\ftree{-}}
\\
\ShuT
  \ar[r]_-{\K}
&
\ShuO
}
\]

For an indexed set $X$, the elements of the free internal monoid $\ftree{X}$ in $\ShuT$ are called \emph{tree monomials on $X$}. We have $\ftree{X}=(\ftree{X}_k)_{k\geq 0}$, and elements of $\ftree{X}_k$ are said to be of \emph{arity $k$}. In particular the unit in $\ftree{X}$ corresponding to the indexed set $\mathbf{1}$ is called the \emph{trivial tree monomial}.
Elements of the free operad $\fopShu{X}$ are linear combination of tree monomials having a same given arity $k$, and called \emph{terms on~$X$ of arity~$k$}.

\subsubsection{Graphical representation of tree monomials}
The elements of an indexed set can be represented graphically as (planar rooted) trees. For an indexed set $X$ and an ordered set $I=\{i_1<i_2<\cdots<i_k\}$, an element $x\in X(I)$ is depicted by a \emph{corolla}, that is, a tree with only one vertex:
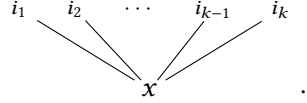
\[
\raisebox{0.5cm}
{\begin{forest} big input tree,
	[x
	  [i_1]
	  [i_2]
	  [\cdots, no edge]
	  [i_{k-1}]
	  [i_k]
	]
\end{forest}}
.\]
For indexed sets $X,Y$, elements of $(X\setShuffleComp Y)(I)$ have the form
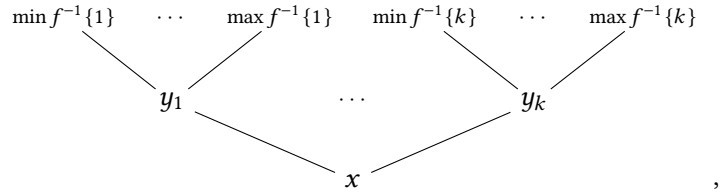
\[
\raisebox{1.1cm}
{\begin{forest} big input tree,
	[x
	  [y_1
	    [\min f^{-1}\{1\}]
	    [\cdots, no edge]
	    [\max f^{-1}\{1\}]
	  ]
	  [\cdots, no edge, before computing xy={s/.average={s}{siblings}}]
	  [y_k
	    [\min f^{-1}\{k\}]
	    [\cdots, no edge]
	    [\max f^{-1}\{k\}]
	  ]
	]
\end{forest}
}
,\]
where $k\geq 1$, $f:I\twoheadrightarrow \range{k}$ is a shuffle surjection, $x\in X(k)$ and $y_i\in Y(f^{-1}\{i\})$ for all $i\in \range{k}$.
In this way, a tree monomial $u$ on $X$ can be represented by a planar tree $T(u)$, whose vertices are elements of $X$, and its arity is the number of its leaves. More generally, for $V,W$ two collections, $V\circ W$ has a basis of tree monomials.

The \emph{weight} of a tree monomial $u$ on an indexed set $X$ is the number of vertices of $T(u)$. A tree monomial $v$ is a (resp.\ \emph{rooted}) \emph{submonomial} of $u$ if $T(v)$ is a (resp.\ rooted) subtree of $T(u)$. If $v$ is a rooted submonomial of $u$, we write $v \subseteq u$. When listing submonomials of a tree monomial $u$, we distinguish the different occurrences of a subtree of $T(u)$: for instance, the tree monomial
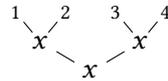
\[
\inputtree{[x[x[1][2]][x[3][4]]]}
\]
contains three distinct occurrences of the submonomial $\inputtree{[x[1][2]]}$, one of which is rooted.

\subsubsection{Inline notation for tree monomials}
For $n_1,\ldots,n_k\geq 1$, a \emph{shuffle surjection of type \linebreak $(n_1,\ldots,n_k)$} is a surjection $f:\range{n_1+\cdots+n_k}\twoheadrightarrow \range{k}$ such that, for all $i\in\range{k}$, $|f^{-1}(\{i\})|=n_i$, and
\[
\min f^{-1}\{1\}<\cdots<\min f^{-1}\{k\}.
\]
Denote by $S(n_1,\ldots,n_k)$ the set of shuffle surjections of type $(n_1,\ldots,n_k)$.

Let $X$ be an indexed set. The \emph{inline notation} for the indexed set of tree monomials $\ftree{X}$ is the term algebra in indexed sets given by the \emph{Backus-Naur form}
\[\ftree{X}::=\mathbf{1}\,\big|\,(X(k)\mid_f \ftree{X}(n_1)\cdots \ftree{X}(n_k)),\]
where $\mathbf{1}$ is the indexed set defined in \ref{SSS:ShuffleCompositionOmegaSets}, and $f$ is a shuffle surjection of type $(n_1,\ldots,n_k)$. When possible, we omit subscript $f$, and we write
\[
(u\mid\vec v):=(u\mid v_1\cdots v_k).
\]
We will use also the notation $\vec v$ for the list of tree monomials $v_1,\ldots,v_k$. Finally, note that for two indexed sets $X,Y$ the indexed set $X\setShuffleComp Y$ can be written in an explicit way as, for $n\geq 1$,
\[
(X\setShuffleComp Y)(n) = \left\{\;(x\mid_f y_1\cdots y_k) \;\left|\; \begin{array}{l}
n_1,\ldots,n_k\geq 1,\;n_1+\cdots+n_k=n, \\
x\in X(k), y_1\in Y(n_1), \ldots, y_k\in Y(n_k), \\
f\in S(n_1,\ldots,n_k)
\;\end{array}\right.\right\}.
\]

\subsubsection{Explicit associativity of shuffle composition}
Let $u,\vec v,\vec v'$ be monomials, and $f,f'$ shuffle surjections such that the tree monomial $((u\mid_f\vec v)\mid_{f'}\vec v')$ is well defined. Denoting the minimal and maximal elements of $f^{-1}\{i\}$ by $f_- i$ and $f_+ i$ respectively for all $i \in \range{k}$, this monomial is represented graphically as
\[
\begin{forest}
for tree={edge=-, l=0, l sep=4, if n children=0
	{node options={font=\scriptsize}, inner ysep=1, child anchor=south}
	{node options={inner ysep=1}}},
[u
  [v_1
    [v'_{f_- 1}
      [f'_- f_- 1]
      [\cdots, no edge]
      [f'_+ f_- 1]
    ]
    [\cdots, no edge, before computing xy={s/.average={s}{siblings}}]
    [v'_{f_+ 1}
      [f'_- f_+ 1]
      [\cdots, no edge]
      [f'_+ f_+ 1]
    ]
  ]
  [\cdots, no edge, before computing xy={s/.average={s}{siblings}}]
  [v_k
    [v'_{f_- k}
      [f'_- f_- k]
      [\cdots, no edge]
      [f'_+ f_- k]
    ]
    [\cdots, no edge, before computing xy={s/.average={s}{siblings}}]
    [v'_{f_+ k}
      [f'_- f_+ k]
      [\cdots, no edge]
      [f'_+ f_+ k]
    ]
  ]
]
\end{forest}
\]
For $i\in\range{k}$, the restriction of $f'$ to $f^{\prime -1}\{f^{-1}\{i\}\}$ is a shuffle surjection of domain $f^{-1}\{i\}$. We introduce the notation $(\vec v\mid_{f,f'}\vec v')$ for the list of monomials $\{(v_i\mid_{f'_{|f^{\prime -1}(f^{-1}\{i\})}} v'_{\min f^{-1}\{i\}}\cdots v'_{\max f^{-1}\{i\}})\}_{1\leq i\leq k}$. Then
\[
((u\mid_f\vec v)\mid_{f'}\vec v') = (u\mid_{ff'}(\vec v\mid_{f,f'}\vec v')),
\]
where $ff'$ is the composition of shuffle surjections.

\subsubsection{Bimodules and ideals}
\label{SSS:Bimodules}
Recall that a \emph{$P$-bimodule} over a shuffle operad $P$, called a \emph{linear module over $P$} in {\cite[Def. 2.13]{BauesJibladzeTonks97}}, is a collection $A$ equipped with two families of morphisms of collections
\begin{align*}
\lambda & :P(k)\otimes P(f^{-1}\{1\})\otimes\cdots\otimes A(f^{-1}\{i\})\otimes\cdots P(f^{-1}\{k\})\fl A(I), \\
\rho & : A(k)\otimes P(f^{-1}\{1\})\otimes\cdots\otimes P(f^{-1}\{k\})\fl A(I),
\end{align*}
for all shuffling surjections $f:I\twoheadrightarrow\range{k}$, defining a \emph{left crossed action} and a \emph{right action} respectively, satisfying compatibility axioms with each other, and associativity and unit axioms with the product of $P$. A morphism of $P$-modules is a morphism of collections compatible with the left and right actions.
We denote by $\Bimod{P}$ the category of $P$-bimodules and morphisms of $P$-bimodules, and by $\Bimod{\ShuO}:=\coprod_{P\in\ShuO}\Bimod{P}$ the category of pairs $(P,A)$ composed of an operad $P$ and a $P$-bimodule $A$.
We denote by $\Lr^P: \Coll \fl \Bimod{P}$ the \emph{free bimodule functor} defined in {\cite[Prop. 2.11]{BauesJibladzeTonks97}}, see also \cite{Markl96}, and given for every $V\in\Coll$ and $I\in\Ord$ by
\[
\Lr^P(V)(I):=\bigoplus_{k\geq 1}P(k)\otimes\left(\bigoplus_{\substack{f:I\twoheadrightarrow\range{k} \\ 1\leq i\leq k}}P(f^{-1}\{1\})\otimes\cdots\otimes (A\circ_\Sha P)(f^{-1}\{i\})\otimes\cdots\otimes P(f^{-1}\{k\})\right)
\]
For an indexed set $X$, we denote by
\[
P\langle X\rangle:=\Lr^P(\K X)
\]
the free $P$-bimodule on~$X$.

An \emph{ideal} of an operad $P$ is a $P$-bimodule $\mathcal{I}$ equipped with an inclusion of $P$-bimodules $\mathcal{I} \hookrightarrow P$.

\subsubsection{Graphical description of bimodules}
Let $X$ be an indexed set an $P$ a shuffle operad. The free $P$-bimodule $P\langle X\rangle$ is the collection generated by tree monomials of the form
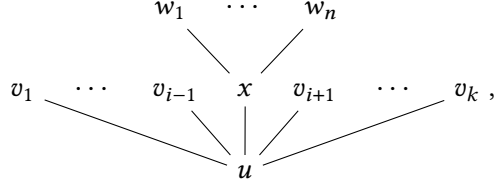
\[\begin{forest}
[u
  [v_1]
  [\cdots, no edge]
  [v_{i-1}]
  [x
    [w_1]
    [\cdots, no edge]
    [w_n]
  ]
  [v_{i+1}]
  [\cdots, no edge]
  [v_k]
]
\end{forest},\]
where $n,k\geq 1$, $x\in X(n)$, $i\in\range{k}$, $u,v_1,\ldots,\check v_i,\ldots,v_k,w_1,\ldots,w_n\in P$, and the inputs are omitted. We use the check notation $\check v_i$ to indicate that we omit $v_i$, as opposed to the usual hat notation $\hat v_i$, in order to avoid confusion with a notation in Section~\ref{S:ShufflePolygraphicResolutionConvergence}.

\subsection{Compositions in shuffle operads}
\label{SS:CompositionShuffleOperads}

In this subsection, we decompose shuffle composition into partial compositions, and we introduce notations for composition and terms in an operad.

\subsubsection{Shuffle partial composition}
\label{SSS:ShuffleCompositions}
Recall from~{\cite[Prop. 2]{DotsenkoKhoroshkin10}}, that for $k,\ell\geq 1$, a \emph{shuffle permutation of type $(k,\ell)$} is a permutation $\tau\in\SS_{k+\ell}$ such that
\[
\tau(1)<\cdots<\tau(k),\quad\text{and}\quad\tau(k+1)<\cdots<\tau(k+\ell).
\]
Denote by $\Sha(k,\ell)$ the set of shuffle permutations of type $(k,\ell)$.
Given indexed sets $X,Y$, and $x\in X(k)$, $y\in Y(\ell)$, for $i\in\range{k}$ and $\tau\in\Sha(\ell-1,k-i)$, we define the \emph{elementary composition} $x\circ_{i,\tau} y$ as the following tree
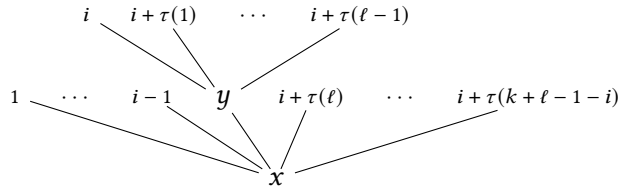
\[
\begin{forest} big input tree,
[x
  [1]
  [\cdots, no edge]
  [i-1]
  [y
    [i]
    [i+\tau(1)]
    [\cdots, no edge]
    [i+\tau(\ell-1)]
  ]
  [i+\tau(\ell)]
  [\cdots, no edge]
  [i+\tau(k+\ell-1-i)]
]
\end{forest}
\]
Elementary compositions are extended to collections by linearity and bidistributivity. 
We denote by $V(k)\circ_{i,\tau}W(\ell)$ the collection composed by elementary compositions of the form $v\circ_{i,\tau}w$, for $v\in V(k)$ and $w\in W(\ell)$. Then the \emph{shuffle partial composition} of collections $V,W$ is defined by
\[
(V\part_\Sha W)(n) :=\bigoplus_{\mathclap{\substack{k,\ell,i\geq 1 \\ k+\ell-1=n \\\tau\in\Sha(\ell-1,k-i)}}} V(k)\circ_{i,\tau} W(\ell).
\]

\subsubsection{Properties of partial compositions}
Note that there are isomorphisms
\begin{align*}
(V\part_\Sha W)(n) & \overset{{\bf (i)}}{\simeq}\bigoplus_{k=1}^\infty V(k)\otimes\left(\bigoplus_{\substack{i,\ell\geq 1 \\ k+\ell-1=n \\\tau\in\Sha(\ell-1,k-i)}}\Unit(\{1\})\otimes\cdots\otimes \Unit(\{i-1\})\otimes W(\{i,i+\tau(1),\ldots,i+\tau(\ell-1)\})\right. \\
& \left.\phantom{:=\bigoplus_{k=1}^\infty V(k)\otimes\left(\bigoplus_{\tau\in\Sha(\ell-1,k-i)}\right.}\otimes\Unit(\{i+\tau(\ell)\})\otimes\cdots\otimes\Unit(\{i+\tau(k+\ell-1-i)\})\right) \\
& \overset{{\bf (ii)}}{\simeq} \bigoplus_{\mathclap{\substack{k,\ell,i\geq 1 \\ k+\ell-1=n \\\tau\in\Sha(\ell-1,k-i)}}} V(k)\otimes W(\ell).
\end{align*}
The isomorphism {\bf (i)} implies that there is an injection of collections
\[
V\part_\Sha W\hookrightarrow V\circ_\Sha (\Unit\oplus W).
\]
The isomorphism {\bf (ii)} implies that partial composition a bidistributive bifunctor
$\part_\Sha:\Coll\times\Coll\fl\Coll$. The partial composition $\part_\Sha$ is not associative. However, if there is no possible confusion, we will use the left bracket rule, that is, $U\part V\part W:=(U\part V)\part W$. 

\subsubsection{Decomposition of shuffle compositions}
\label{SSS:DecompositionShuffleComposition}
Let $A$ be a collection equipped with a morphism $\eta_A:\Unit\fl A$, that is, $\eta_A$ is an object of $\Unit/\Coll$. We can express the shuffle composition $A\circ_\Sha A$ in terms of partial compositions. There exists a natural transformation $\varphi$ from the functor $A\mapsto\bigoplus A\part_\Sha A^{\part_\Sha p}$ to the functor $A\mapsto A\circ_\Sha A$ defined as follows:
\[
\varphi_A:\bigoplus_{p=1}^\infty A\part_\Sha A^{\part_\Sha p}\rightarrow\bigoplus_{p=1}^\infty A\circ_\Sha(\Unit\oplus A)^{\circ_\Sha p}\xrightarrow{\sum1\circ_\Sha(\eta_A+1)^{\circ_\Sha p}} A\circ_\Sha A.
\]
In order to express $A\circ_\Sha A$ in terms of partial compositions, it suffices to define a right inverse to $\varphi$, that is, a natural transformation $\sigma$ from the functor  $A\mapsto A\circ_\Sha A$ to the functor $A\mapsto\bigoplus A\part_\Sha A^{\part_\Sha p}$ such that $\varphi_A\sigma_A=\id_{A\circ_\Sha A}$ for all $A$.
	
Define the morphism
\[
\sigma_A:A\circ_\Sha A\rightarrow\bigoplus_{p=1}^\infty A\part_\Sha A^{\part_\Sha p},
\]
natural in $A$, as follows. An element $a$ of $(A\circ_\Sha A)(n)$ can be written 
\[
a=\begin{forest} big input tree,
[a_0
  [a_1
    [\min f^{-1}\{1\}]
    [\cdots, no edge]
    [\max f^{-1}\{1\}]
  ]
  [\cdots, no edge, before computing xy={s/.average={s}{siblings}}]
  [a_p
    [\min f^{-1}\{p\}]
    [\cdots, no edge]
    [\max f^{-1}\{p\}]
  ]
]
\end{forest}
\]
where $f:\range{n}\twoheadrightarrow\range{p}$ is a shuffle surjection. Set
$
\sigma_A(a):=a_0\circ_{p,\tau_p}a_p\circ_{p-1,\tau_{p-1}}\cdots\circ_{1,\tau_1}a_1,
$
where
\begin{align*}
\tau_p & \in\Sha\left(|f^{-1}\{p\}|-1,0\right)=\{\id_{|f^{-1}\{p\}|-1}\}, \\
\tau_{p-1} & \in\Sha\left(|f^{-1}\{p-1\}|-1,|f^{-1}\{p\}|\right), \\
\vdots \\
\tau_1 & \in\Sha\left(|f^{-1}\{1\}|-1,|f^{-1}\{2,\ldots,p\}|\right),
\end{align*}
are the appropriate shuffle permutations. We check that, for every morphism $f:A\rightarrow B$ of collections, the square
\[
\xymatrix @C=2.7em
{
A\circ_\Sha A
  \ar[r]^-{\sigma_A}
  \ar[d]_-{f\circ_\Sha f}
&
\bigoplus_{p\geq 0} A\part_\Sha A^{\part_\Sha p}
  \ar[d]^-{\sum_{p\geq 0} f\part_\Sha f^{\part_\Sha p}}
\\
B\circ_\Sha B
  \ar[r]_-{\sigma_B}
&
\bigoplus_{p\geq 0} B\part_\Sha B^{\part_\Sha p}
}
\]
commutes. This defines the natural transformation $\sigma$, and we check that it is a right inverse to the natural transformation $\varphi$.

\subsubsection{Example}
\label{SSS:DecompositionExample}
Let $A$ be an object of $\Unit/\Coll$ and
\[
a=\inputtree{[a_0[a_1[1][3]][a_2[2][4]]]}
\]
and element of $A\circ_\Sha A$, where $a_0,a_1,a_2\in A(2)$. Then we have
\[
\sigma_A(a)=a_0\circ_{2,\id_1}a_2\circ_{1,(1\,2)}a_1
\qquad\text{and}\qquad
\varphi_A\sigma_A(a)=\inputtree{[a_0[1[a_1[1][3]]][a_2[2][4]]]}=\inputtree{[a_0[a_1[1][3]][a_2[2][4]]]}.
\]

\subsection{Higher shuffle operads}
\label{SS:HigherDimensionalShuffleOperads}

In this subsection, we introduce the structure of (strict) shuffle $\omega$-operads. 

\subsubsection{Globular objects}
We denote by $\Rglob$ the \emph{reflexive globe category}, whose objects are natural numbers, denoted by $\underline n$, for $n\in \Nb$, and morphisms are generated by 
\[ 
\sigma_n : \underline n \fl \underline{n+1},
\qquad
\tau_n : \underline n \fl \underline{n+1},
\qquad
\iota_{n+1} : \underline{n+1} \fl \underline n,
\]
for all $n$ in $\Nb$, and submitted to the following \emph{globular} and \emph{identities relations}:
\begin{align*}
\sigma_{n+1} \circ \sigma_n = \tau_{n+1} \circ \sigma_n,
&\qquad
\sigma_{n+1}\circ \tau_n = \tau_{n+1} \circ \tau_n,
\\
\iota_n \circ \sigma_n = \id_{\underline n},
&\qquad
\iota_n\circ \tau_n = \id_{\underline n},
\end{align*}
for all $n$ in $\Nb$. 
Omitting the identity maps $\iota_n$ gives the definition of the \emph{globe category} $\glob$. We denote by~$\Rglob_n$ (resp.\ $\glob_n$) the full subcategory of $\Rglob$ (resp.\ $\glob$) whose objects are $\underline 0$, $\underline 1$, $\ldots$, $\underline n$.

A \emph{reflexive globular object} in a category $\Cr$ is a functor $\Rglob^{op} \fl \Cr$, whose restriction to the category~$\Rglob_n^{op}$ is called a \emph{reflexive $n$-globular object}. Explicitly, a reflexive globular object is given by a sequence $A=(A_k)_{k\in \Nb}$ of objects of $\Cr$, equipped with indexed morphisms
\[
s = (s_k : A_{k+1} \fl A_k )_{k\in\Nb},
\qquad
t = (t_k : A_{k+1} \fl A_k )_{k\in\Nb},
\qquad
i = (i_k : A_{k-1} \fl A_{k} )_{k\in\Nb},
\]
of degree $-1$, $-1$ and $1$ respectively, and satisfying the following globular and identities relations
\begin{eqn}{equation}
\label{E:GlobularIdentityRelations}
s^2 = st,
\quad
t^2 = ts,
\quad
si = \id_A,
\quad
ti = \id_A.
\end{eqn}
The elements of $A_k$ are called \emph{$k$-cells} of $A$.
A \emph{morphism of reflexive globular objects} is an indexed morphism of degree $0$ that commutes with morphisms $s$, $t$ and $i$.
We denote by $\RGlob(\Cr)$ (resp.\ $\Glob(\Cr)$) the category of \emph{reflexive globular objects} (resp.\ \emph{globular objects}) in $\Cr$ and their morphisms.
We denote by $\RGlob_n(\Cr)$ (resp.\ $\Glob_n(\Cr)$) the full subcategory of $\RGlob(\ShuO)$ of reflexive $n$-globular objects (resp.\ $n$-globular objects) in $\Cr$. 
We will denote by
\[
\Vr_n(\Cr) : \Glob_{n+1}(\Cr) \fl \Glob_n(\Cr)
\]
the functor that forgets the $(n+1)$-cells.
For $A$ a globular object and $\ell\geq k \geq 0$, the \emph{$(\ell,k)$-source} and \emph{$(\ell,k)$-target morphisms} 
\[
s_k^\ell: A_\ell \fl A_k,
\qquad
t_k^\ell : A_\ell \fl A_k,
\]
are respectively defined as the following iterated composition of source and target morphisms:
\[
s_k^\ell := s_k\circ \ldots s_{\ell-2}\circ s_{\ell-1},
\qquad
t_k^\ell := t_k \circ \ldots t_{\ell-2} \circ t_{\ell-1}.
\]
We denote by $i_\ell^k : A_k \fl A_\ell$
the iterated identity $i_l^k = i_\ell \circ i_{\ell-1}\ldots \circ i_{k+1}$.
When there is no ambiguity, we will write $s_k$ and $t_k$ for source and target maps respectively, and we will omit $i_\ell^k$ entirely, since $i_\ell^k$ is injective by \eqref{E:GlobularIdentityRelations}. 
For $k\geq 0$, we denote by $A\star_k A$ the following pullback of globular operads
\[
\xymatrix @R=1.75em {
A \star_k A
	\ar[r] ^-{\pi_1}
	\ar[d] _-{\pi_2}
	\pullbackcorner
& A
	\ar[d] ^-{s_k}
\\
A
	\ar[r] _-{t_k}
& A_k
}
\]

Let $A$ be a globular object of some category $\Cr$. For $n\geq 1$, two $n$-cells $a,b$ of $A$ are \emph{parallel} if $s(a)=s(b)$ and $t(a)=t(b)$. An \emph{$n$-sphere of $A$} is a pair $(a,b)$ of parallel $n$-cells.

\subsubsection{Higher categories for operads}
Recall that, for $n\geq 0$, an (internal strict) \emph{$n$-category} in $\Cr$ is a
\begin{enumerate}
\item reflexive $n$-globular object, that is a diagram in $\Cr$ of the form
\[
\xymatrix @C=3.3em @W=2.5em{
A_0 
	\ar [r] |-*+{i_1}
&
A_1 
	\ar@<-1.3ex> [l] _-{s_0}
	\ar@<1.3ex> [l] ^-{t_0}
	\ar [r] |-*+{i_2}
& 
\quad
\cdots
\quad
	\ar@<-1.3ex> [l] _-{s_1}
	\ar@<1.3ex> [l] ^-{t_1}
	\ar [r] |-*+{i_{n-1}}
& 
\;\;A_{n-1}\;\;
	\ar@<-1.3ex> [l] _-{s_{n-2}}
	\ar@<1.3ex> [l] ^-{t_{n-2}}
	\ar [r] |-*+{i_n}
&
A_n
	\ar@<-1.3ex> [l] _-{s_{n-1}}
	\ar@<1.3ex> [l] ^-{t_{n-1}}
}
\]
whose morphisms satisfy globular and identity relations~(\ref{E:GlobularIdentityRelations}),
\item equipped with a structure of category in $\Cr$ on
\[
\xymatrix @C=3em @W=2.5em{
A_k
&
A_\ell
	\ar@<-0.7ex> [l] _-{s_k}
	\ar@<0.7ex> [l] ^-{t_k}
}
\]
for all $k<\ell$, whose \emph{$\star_k$-composition morphism of $\ell$-cells} is denoted by $\star_k^\ell: A_\ell\star_k A_\ell \fl A_\ell$,
\item such that the $2$-globular object
\[
\xymatrix @C=3em @W=2.5em{
A_j
&
A_k 
	\ar@<-0.7ex> [l] _-{s_j}
	\ar@<0.7ex> [l] ^-{t_j}
& 
A_\ell
	\ar@<-0.7ex> [l] _-{s_k}
	\ar@<0.7ex> [l] ^-{t_k}
}
\]
is a $2$-category in $\Cr$ for all $j<k<\ell$. 
\end{enumerate}
We denote by $n\Cat(\Cr)$ the category of $n$-categories in $\Cr$ and their $n$-functors.
The category $\omega\Cat(\Cr)$ of \emph{$\omega$-categories in $\Cr$} is the limit of
\[
0\Cat(\Cr) \leftarrow 1\Cat(\Cr) \leftarrow \cdots \leftarrow n\Cat(\Cr) \leftarrow \cdots
\]
where each arrow forgets the cells of highest dimension.

For $n\in \Nb\cup\{\omega\}$, a \emph{shuffle} (resp.\ \emph{symmetric}) \emph{$n$-operad} is an $n$-category in $\ShuO$ (resp.\ $\SymO$). We denote by $\ShuO_n$ (resp.\ $\SymO_n)$ the corresponding category with internal $n$-functors as morphisms.
We denote by $\Ur_n^\Sha : \ShuO_n \fl \Glob_n(\Ind)$ (resp.\ $\Ur_n^\SS : \SymO_n \fl \Glob_n(\SSColl)$) the forgetful functor that forgets the operadic structure.

Note that, in an $\omega$-shuffle operad $P$, the composition $\star_\ell^k : P_k\star_\ell P_k \fl P_k$ is a morphism of shuffle operads. As a consequence, the composition satisfies the following \emph{exchange relation between $\circ_\Sha$ and $\star_k$}
\begin{eqn}{equation}
\label{E:MultiplicativeExchangeRelation}
\raisebox{0.5cm}{
\xymatrix @C=5em{
(P_k\star_\ell P_k)\circ_\Sha(P_k\star_\ell P_k)
  \ar[rr]^-{(\pi_1\circ_\Sha\pi_1)\star_\ell (\pi_2\circ_\Sha\pi_2)}
  \ar[dr]_-{\mu_k((1\star_\ell 1)\circ_\Sha(1\star_\ell 1))\qquad \quad }
&
&
(P_k\circ_\Sha P_k)\star_\ell (P_k\circ_\Sha P_k)
  \ar[dl]^-{\quad \qquad\mu_k(1\circ_\Sha 1)\star_\ell \mu_k(1\circ_\Sha 1)}
\\
&
P_n
}}
\end{eqn}

As for associative $\omega$-algebras, {\cite[Prop. 1.2.3]{GuiraudHoffbeckMalbos19}}, the interaction between the categorical and linear structures gives useful expressions:
\begin{lemma}\label{star-sum}
Let $P$ be a shuffle (resp.\ symmetric) $\omega$-operad. 
\begin{enumerate}
\item For every $0 \leq k < n$ and $\star_k$-composable pair $(a, b)$ of $P_n$, we have $a \star_k b = a - t_k(a) + b$.
\item For all $n \geq 1$, every $n$-cell $a$ of $P$ is invertible with inverse $a^- := s_{n - 1}(a) - a + t_{n - 1}(a)$.
\end{enumerate}
\end{lemma}

We deduce the following proposition:

\begin{proposition}
The category $\ShuO_\omega$ (resp.\ $\SymO_\omega)$ is isomorphic to the category $\Gpd_\omega(\ShuO)$ (resp.\ $\Gpd_\omega(\SymO)$) of internal $\omega$-groupoids in $\ShuO$ (resp.\ $\SymO$).
\end{proposition}

Finally, we give some categorical properties of the categories of shuffle operads.

\begin{proposition}
\label{P:reflectionLimits}
The forgetful functor $\ShuO\fl\Coll$ (resp.\ $\SymO\fl\SSColl$) reflects all limits, filtered colimits, and reflexive coequalizers.
\end{proposition}
\begin{proof}
The statement for the functor $\SymO\fl\SSColl$ is proven in \cite[Prop. 1.2.4]{Fresse17}, where the limits, filtered colimits, and reflexive coequalizers of symmetric operads are equipped with unique monoidal structures. In particular, it is shown that the monoidal product $\circ_\SS$ in $\SSColl$ preserves all limits, filtered colimits and reflexive coequalizers. This comes from the fact that the tensor product of $\Vect$ preserves these limits and colimits. Given the similarities between the monoidal products $\circ_\SS$ and $\circ_\Sha$, the same arguments apply to the monoidal structure of $\Coll$, and so we conclude that the functor $\ShuO\fl\Coll$ preserves limits, filtered colimits, and reflexive coequalizers.
\end{proof}

\begin{proposition}
\label{P:LocallyFinitelyPresentable}
The category $\ShuO_\omega$ is locally finitely presentable. In particular, it is complete and cocomplete.
\end{proposition}
\begin{proof}
Let us first show that $\ShuO$ is locally finitely presentable by viewing it as the category of algebras over an accessible monad. Recall that there exists an adjunction $\Coll\dashv\ShuO$ where the left adjoint is $\fopShu{-}:\Coll\fl\ShuO$. Therefore, $\ShuO$ is the category of algebras of the monad of free shuffle operads $T:\Coll\fl\Coll$. By the Proposition \ref{P:reflectionLimits}, the forgetful functor $\ShuO\fl\Coll$ preserves filtered colimits, \emph{i.e.} it is finitary, making the monad $T$ finitary. Moreover, the category $\Vect$ of vector spaces is locally finitely presentable, and the category $\Ord$ is a small category, so $\Coll$ is also locally finitely presentable. Thus $T$ is an accessible monad on a locally finitely presentable category. Following~\cite[§ 2.78]{AdamekRosicky94} the category of $T$-algebras $\ShuO$ is locally finitely presentable.

The category $\ShuO_\omega$ of $\omega$-categories internal in $\ShuO$ is the category of models of a finite limit sketch~\cite{Polybook21}, in the locally finitely presentable category $\ShuO$. By \cite[Prop. 1.53]{AdamekRosicky94}, we conclude that $\ShuO_\omega$ is also locally finitely presentable.  
\end{proof}

\subsection{Higher operads as globular bimodules}
\label{S:HigherDiomensionalOperadsAsBimobules}

In Theorem \ref{T:MainTheoremA}, we show that the axioms of the definition of the category $\ShuO_\omega$ of shuffle $\omega$-operads are redundant by proving that it is isomorphic to a category with fewer axioms thanks to the linear exchange relation.
Throughout this section, the compositions $\circ$ and $\part$ designate $\circ_\Sha$ and $\part_\Sha$.

\subsubsection{Partial multiplication}
For $(P,\mu,\eta)$ a shuffle operad, denote by $\iota_P$ the morphism
\[
\iota^{}_P:P\part P\hookrightarrow P\circ(\Unit\oplus P)\xrightarrow{1\circ(\eta+1)}P\circ P.
\]
We equip the operad $P$ with a morphism called \emph{partial multiplication}
\[
\mu^\part:P\part P\xrightarrow{\iota^{}_P}P\circ P\xrightarrow{\mu}P.
\]
As a consequence, we have the equality of morphisms
\[
\mu\varphi^{}_{P}=\sum_p(\mu^\part)^{\part p}:\bigoplus_{p=1}^\infty P\part P^{\part p}\rightarrow P.
\]

\subsubsection{Linear exchange relation}
\label{SSS:LinearExchangeRelation}
Let $(P_n,\mu_n,\eta_n)_{n\geq 0}$ be an $\omega$-operad. By the exchange relation between compositions $\circ$ and $\star_0$, we observe that, for every $n\geq 1$,
\begin{align*}
\mu_n^\part & =\mu_n((1\star_0t_0)\circ(s_0\star_01))\iota^{}_{P_n} \\
& =(\mu_n(1\circ s_0)\star_0\mu_n(t_0\circ 1))\iota^{}_{P_n} \tag{\ref{E:MultiplicativeExchangeRelation}} \\
& =(\mu_n(1\circ s_0)+\mu_n(t_0\circ 1)-\mu_n(t_0\circ s_0))\iota^{}_{P_n} \tag{Lemma \ref{star-sum}} \\
& =\mu_n^\part(1\part s_0)+\mu_n^\part(t_0\part 1)-\mu_n^\part(t_0\part s_0).
\end{align*}
Similarly, we calculate
\[
\mu_n^\part = \mu_n^\part(s_0\part 1)+\mu_n^\part(1\part t_0)-\mu_n^\part(s_0\part t_0).
\]
Regarding $P_n$ as a $P_0$-bimodule, these equations still hold, although we need to introduce new notations for the partial left and right actions of $P_0$ on $P_n$. This motivates the following definitions.

Let $(P,\mu,\eta)$ be a shuffle operad and $(A,\lambda,\rho)$ be a $P$-bimodule such that $(P,A)$ is a reflexive $1$-globular $P$-bimodule. More explicitly, there are morphisms $s,t:A\fl P$ and $i:P\fl A$. We equip $A$ with morphisms called \emph{partial actions}
\begin{align*}
\lambda^\part:P\part A&\xrightarrow{\varphi} P\circ(\Unit\oplus A)\xrightarrow{1\circ(\eta\oplus 1)} P\circ(P\oplus A)\xrightarrow{\lambda} A, \\
\rho^\part:A\part P&\xrightarrow{\varphi}A\circ(\Unit\oplus P)\xrightarrow{1\circ(\eta+1)}A\circ P\xrightarrow{\rho}A.
\end{align*}
We also define the morphisms
\begin{align*}
\mu_A^\uparrow & :=\rho^\part(1\part s)+\lambda^\part(t\part 1)-i\mu^\part(t\part s), \\
\mu_A^\downarrow & :=\lambda^\part(s\part 1)+\rho^\part(1\part t)-i\mu^\part(s\part t)
\end{align*}
and the multiplication
\[
\mu_A:A\circ A\xrightarrow{\sigma_A}\bigoplus_pA\part A^{\part p}\xrightarrow{\sum(\mu_A^\uparrow)^p}A.
\]
We say that a reflexive $1$-globular $P$-bimodule $(P,A)$ satisfies the \emph{linear exchange relation} if the following relation holds:
\begin{eqn}{equation}
\label{E:LinearExchangeRelation}
\mu_A^\uparrow=\mu_A^\downarrow.
\end{eqn}

\subsubsection{Interpretation of morphisms $\mu_A^\uparrow$ and $\mu_A^\downarrow$}
If $(P,A)$ is a reflexive $1$-globular $P$-bimodule, we can interpret the elements of $A$ as rewriting rules that relate elements of $P$: an element $a\in A$ rewrites $s(a)$ as $t(a)$, which we denote by $a:s(a)\fl t(a)$. Via the injection $i$, an element $x$ of $P$ can also seen as a trivial rewriting rule $i(x):x\fl x$.

Let $a,b\in A$. For every compatible elementary composition $\circ_{i,\tau}$, we would like to interpret the composition $a\circ_{i,\tau}b$ as a pair of orthogonal reductions:
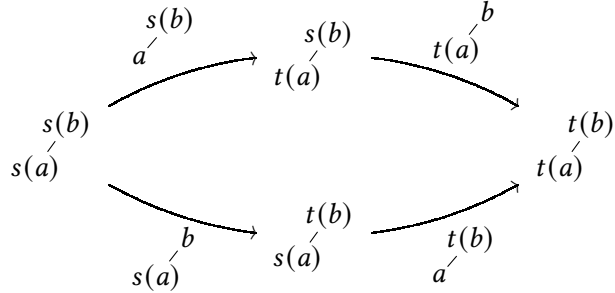
\[
\xymatrix@C=5em @R=0em{
&
\segment{t(a)}{s(b)}
  \ar@/^2ex/ [dr] ^-{\segment{t(a)}{b}}
&
\\
\segment{s(a)}{s(b)}
  \ar@/^2ex/ [ur] ^-{\segment{a}{s(b)}}
  \ar@/_2ex/ [dr] _-{\segment{s(a)}{b}}
&& 
\segment{t(a)}{t(b)}
\\
&
\segment{s(a)}{t(b)}
  \ar@/_2ex/ [ur] _-{\segment{a}{t(b)}}
&
}
\]
where $\segment{s(a)}{s(b)}$ is a graphical representation of $s(a)\circ_{i,\tau}s(b)$, and so on.
While $a\circ_{i,\tau}b$ is not necessarily an element of $A$, we find that
\begin{align*}
\mu_A^\uparrow\left(\segment{a}{b}\right)\:=\:\segment{a}{s(b)}-\segment{t(a)}{s(b)}+\segment{t(a)}{b}, \\
\mu_A^\downarrow\left(\segment{a}{b}\right)\:=\:\segment{s(a)}{b}-\segment{s(a)}{t(b)}+\segment{a}{t(b)}.
\end{align*}
We see that $\mu_A^\uparrow$ applies the rule $a$ first, and $b$ second, while $\mu_A^\downarrow$ does the opposite; this motivates the upwards and downwards arrow notations.

\subsubsection{Example}
Let $a=\tree{[a_0[a_1][a_2]]}$ 
be the element of $A\circ A$ from Example~\ref{SSS:DecompositionExample} with inputs omitted. Then
\[
\mu_A(a)\:=\:\tree{[a_0[s(a_1)][s(a_2)]]}-\tree{[t(a_0)[s(a_1)][s(a_2)]]}+\tree{[t(a_0)[s(a_1)][a_2]]}-\tree{[t(a_0)[s(a_1)][t(a_2)]]}+\tree{[t(a_0)[a_1][t(a_2)]]}.
\]

\begin{lemma}
\label{L:3-implies-2}
Let $(A,\lambda,\rho)$ be a $P$-bimodule such that $(P,A)$ is a reflexive $1$-globular $P$-bimodule satisfying the linear exchange relation. Then $(A,\mu_A,i\eta)$ is an operad.
\end{lemma}
\begin{proof}
Write $\mu_A^\part:=\mu_A^\uparrow=\mu_A^\downarrow$. It suffices to check the associativity and unit axioms of internal monoidal objects. The unit axioms are clearly satisfied, by definition of $\sigma_A$ and by the unit axioms of $P$-bimodules.

To show the associativity axiom, we need to calculate and compare $\mu_A(\mu_A\circ 1)$ and $\mu_A(1\circ\mu_A)$. The key calculation is the following, which generalizes the previous example: for all $p\geq 0$ and $\tau\in\SS_p$, we have the equality of morphisms
\begin{eqn}{equation}
\label{E:KeyCalculation}
\mu_A^\part(1^{\part p})=\sum_{i=1}^pf_{i,1}^\tau\part\cdots\part\undermath{\tau(i)}{1}\part\cdots\part f_{i,p}^\tau-\sum_{i=1}^{p-1}f_{i,1}^\tau\part\cdots\part\undermath{\tau(i)}{t}\part\cdots\part f_{i,p}^\tau,
\end{eqn}
where
\[
f_{i,j}^\tau=\left\{\begin{array}{ll}
	t & \text{if }\tau^{-1}(j)<i, \\
	s & \text{if }\tau^{-1}(j)>i,
\end{array}\right.
\]
and $-\part-$ represents $\lambda^\part(-\part-)$, $\rho^\part(-\part-)$, or $\mu^\part(-\part-)$ depending on the types of arguments, always with bracketing to the left.

We show this equality by induction on $p$. For $p=0,1$, the result is trivial. For $p=2$, for $\tau=\id$,
\[
\mu_A^\part=\mu_A^\uparrow=\rho^\part(1\part s)+\lambda^\part(t\part 1)-i_n\mu^\part(t\part s),
\]
by definition, and for $\tau=(1\,2)$,
\[
\mu_A^\part=\mu_A^\downarrow=\rho^\part(1\part t)+\lambda^\part(s\part 1)-i_n\mu^\part(s\part t)
\]
by hypothesis on $(P,A)$.
	
Let $p\geq 2$ and suppose that we have shown the equality for $p$. Let $\tau\in\SS_{p+1}$, and denote $i_0=\tau^{-1}(p+1)$ and $\tau' = \tau \, (i_0\;p+1)$. If $i_0 < p + 1$, then
\begin{align*}
& (\mu_A^\part)^{p + 1}(1^{\part p + 1}) =
\mu_A^\part((\mu_A^\part)^p(1^{\part p}) \part 1)
\\
={} &
\mu_A^\part \left( \sum_{i = 1}^p f_{i,1}^{\tau'} \part \cdots \part \undermath{\tau'(i)}{1} \part \cdots \part f_{i,p}^{\tau'} \part 1
- \sum_{i = 1}^{p - 1} f_{i, 1}^{\tau'} \part \cdots \part \undermath{\tau'(i)}{t} \part \cdots \part f_{i, p}^{\tau'} \part 1 \right)
\\
={} &
\left( \sum_{i = 1}^{i_0 - 1} f_{i, 1}^{\tau'} \part \cdots \part \undermath{\tau'(i)}{t} \part \cdots \part f_{i, p}^{\tau'} \part 1
+ \sum_{i = i_0}^{p} f_{i, 1}^{\tau'} \part \cdots \part \undermath{\tau'(i)}{s} \part \cdots \part  f_{i, p}^{\tau'} \part 1 \right.
\\
&
+ \sum_{i = 1}^{i_0 - 1} f_{i, 1}^{\tau'} \part \cdots \part \undermath{\tau'(i)}{1} \part \cdots \part f_{i, p}^{\tau'} \part s
+ \sum_{i = i_0}^{p} f_{i, 1}^{\tau'} \part \cdots \part \undermath{\tau'(i)}{1} \part \cdots \part  f_{i, p}^{\tau'} \part t
\\
&
\left. - \sum_{i = 1}^{i_0 - 1} f_{i, 1}^{\tau'} \part \cdots \part \undermath{\tau'(i)}{t} \part \cdots \part f_{i, p}^{\tau'} \part s
- \sum_{i = i_0}^{p} f_{i, 1}^{\tau'} \part \cdots \part \undermath{\tau'(i)}{s} \part \cdots \part  f_{i, p}^{\tau'} \part t \right)
\\
&
- \left( \sum_{i = 1}^{i_0 - 1} f_{i, 1}^{\tau'} \part \cdots \part \undermath{\tau'(i)}{t} \part \cdots \part f_{i, p}^{\tau'} \part 1
+ \sum_{i = i_0 + 1}^{p} f_{i, 1}^{\tau'} \part \cdots \part \undermath{\tau'(i)}{s} \part \cdots \part f_{i, p}^{\tau'} \part 1 \right)
\\ 
={} &
f_{i_0, 1}^{\tau'} \part \cdots \part \undermath{\tau'(i_0)}{s} \part \cdots \part f_{i_0, p}^{\tau'} \part 1 \\ 
&
+ \sum_{\substack{i = 1 \\ i \neq i_0}}^p f_{i, 1}^{\tau'} \part \cdots \part \undermath{\tau'(i)}{1} \part \cdots \part f_{i, p}^{\tau'} \part f_{i, p + 1}^\tau + f_{i_0, 1}^{\tau'} \part \cdots \part \undermath{\tau'(i_0)}{1} \part \cdots \part f_{i_0, p}^{\tau'} \part t
\\ 
&
- \sum_{\substack{i = 1 \\ i \neq i_0}}^{p} f_{i, 1}^{\tau'} \part \cdots \part \undermath{\tau'(i)}{t} \part \cdots \part f_{i, p}^{\tau'} \part f_{i, p + 1}^\tau
- f_{i_0, 1}^{\tau'} \part \cdots \part \undermath{\tau'(i_0)}{s} \part \cdots \part f_{i_0, p}^{\tau'} \part t
\\ 
={} &
\sum_{i = 1}^{p + 1} f_{i, 1}^{\tau} \part \cdots \part \undermath{\tau(i)}{1} \part \cdots \part f_{i, p + 1}^\tau
- \sum_{i = 1}^{p} f_{i, 1}^{\tau} \part \cdots \part \undermath{\tau(i)}{t} \part \cdots \part f_{i, p + 1}^\tau.
\end{align*}
If $i_0 = p + 1$, then
\begin{align*}
& (\mu_A^\part)^{p + 1}(1^{\part p + 1}) =
\mu_A^\part((\mu_A^\part)^p(1^{\part p}) \part 1)
\\
={} & 
\mu_A^\part \left( \sum_{i = 1}^p f_{i,1}^{\tau'} \part \cdots \part \undermath{\tau'(i)}{1} \part \cdots \part f_{i,p}^{\tau'} \part 1
- \sum_{i = 1}^{p - 1} f_{i, 1}^{\tau'} \part \cdots \part \undermath{\tau'(i)}{t} \part \cdots \part f_{i, p}^{\tau'} \part 1 \right)
\\
={} &
\sum_{i = 1}^{p} f_{i, 1}^{\tau'} \part \cdots \part \undermath{\tau'(i)}{t} \part \cdots \part f_{i, p}^{\tau'} \part 1
+ \sum_{i = 1}^{p} f_{i, 1}^{\tau'} \part \cdots \part \undermath{\tau'(i)}{1} \part \cdots \part f_{i, p}^{\tau'} \part s
\\
&
- \sum_{i = 1}^{p} f_{i, 1}^{\tau'} \part \cdots \part \undermath{\tau'(i)}{t} \part \cdots \part f_{i, p}^{\tau'} \part s
- \sum_{i = 1}^{p - 1} f_{i, 1}^{\tau'} \part \cdots \part \undermath{\tau'(i)}{t} \part \cdots \part f_{i, p}^{\tau'} \part 1
\\
={} &
\sum_{i = 1}^{p + 1} f_{i, 1}^{\tau} \part \cdots \part \undermath{\tau(i)}{1} \part \cdots \part f_{i, p + 1}^\tau
- \sum_{i = 1}^{p} f_{i, 1}^{\tau} \part \cdots \part \undermath{\tau(i)}{t} \part \cdots \part f_{i, p + 1}^\tau.
\end{align*}
This proves \eqref{E:KeyCalculation} for $p + 1$. Next, let $a$ be an arbitrary tree monomial of $A\circ A\circ A\simeq (A\circ A)\circ A\simeq A\circ(A\circ A)$. We write
\begin{align*}
a & =((a_0\mid a_1\cdots a_{k_0}) \mid a_{k_0+1} \cdots a_{k_0+k_1} \cdots a_{k_0+\cdots+k_{k_0-1}+1}\cdots a_{k_0+\cdots+k_{k_0}}) \\
& =(a_0\mid (a_1\mid a_{k_0+1}\cdots a_{k_0+k_1}) \cdots (a_{k_0}\mid a_{k_0+\cdots+{k_0+\cdots k_{k_0-1}+1}}\cdots a_{k_0+\cdots+k_{k_0}})),
\end{align*}
where $k_0$ is the arity of $a_0$ and, for $i\in\range{k_0}$, $k_i$ is the arity of $a_i$. We can also understand $a$ via its planar tree $T(a)$:
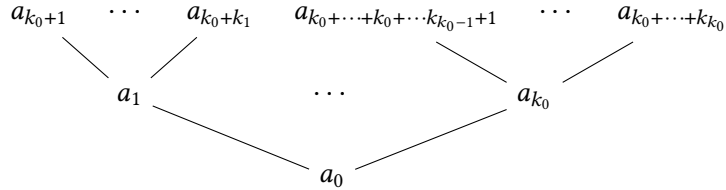
\[
\begin{forest}
[a_0
  [a_1
    [a_{k_0+1}]
    [\cdots, no edge]
    [a_{k_0+k_1}]
  ]
  [\cdots, no edge, before computing xy={s/.average={s}{siblings}}]
  [a_{k_0}
    [a_{k_0+\cdots+{k_0+\cdots k_{k_0-1}+1}}]
    [\cdots, no edge]
    [a_{k_0+\cdots+k_{k_0}}]
  ]
]
\end{forest}
\]
For the rest of this proof, we write $\cdot$ for all elementary compositions $\circ_{i,\tau}$. On the one hand, calculating $\mu_A(\mu_A\circ 1)(a)$ is equal to calculating $(\mu_A^\part)(a_0\cdot a_{k_0}\cdots a_{1}\cdot a_{k_0+\cdots+k_{k_0}} \cdots a_{k_0+1})$ as in the previous calculation, with $p=k_0+\cdots+k_{k_0}$ and the identity permutation. On the other hand, calculating $\mu_A(1\circ\mu_A)(a)$ is equal to calculating $(\mu_A^\part)(a_0\cdot a_{k_0}\cdots a_{1}\cdot a_{k_0+\cdots+k_{k_0}} \cdots a_{k_0+1})$ as in the previous calculation, with $p=k_0+\cdots+k_{k_0}$ and the permutation
\[
[0,k_0,k_0+\cdots+k_{k_0},\ldots,k_0+\cdots+k_{k_0-1}+1,\ldots,1,k_0+k_1,\ldots,k_0+1],
\]
where each integer $i$ represents the position of $a_i$ in the argument of $\mu_A^\part$. Thus, by the previous calculation, $\mu_A(\mu_A\circ1)=\mu_A(1\circ\mu_A)$. We conclude that $(A,\mu_A,i\eta)$ is an operad.

As an aside, calculating $\mu_A(\mu_A\circ 1)$ corresponds to a breadth-first traversal of a tree, while calculating $\mu_A(1\circ\mu_A)$ corresponds to a depth-first traversal.
\end{proof}

\begin{theorem}
\label{T:MainTheoremA}
The following categories are isomorphic:
\begin{enumerate}
\item the category $\ShuO_\omega$,
\item the full subcategory of $\RGlob(\Bimod{\ShuO})$ whose objects are pairs $(P,A)$ where $(P,\mu,\eta)$ is a shuffle operad and $A=(A_n,\lambda_n,\rho_n)_{n\geq 0}$ is a reflexive globular $P$-bimodule such that $A_0=P$ and $(P,A_n)$ satisfies the linear exchange relation \eqref{E:LinearExchangeRelation} for all $n\geq 1$.
\end{enumerate}
\end{theorem}
\begin{proof}
We show that each category is a full subcategory of the other.

($\textbf{i}\subseteq\textbf{ii}$) Let $P=(P_n,\mu_n,\eta_n)_{n\geq 0}$ be an $\omega$-operad. Forgetting the $\star_k$-compositions and operadic multiplications $\mu_n$, $P$ is equipped with a reflexive globular $P_0$-bimodule structure. By the calculations and discussion of \ref{SSS:LinearExchangeRelation}, for all $n\geq 1$, $P_n$ seen as a $P_0$-bimodule satisfies the linear exchange relation \eqref{E:LinearExchangeRelation}. Thus $\ShuO_\omega$ is a full subcategory of the second category.\newline

($\textbf{ii}\subseteq\textbf{i}$)
Let $(P,\mu,\eta)$ be an operad and $(A_n,\lambda_n,\rho_n)_{n\geq 0}$ a globular reflexive $P$-bimodule satisfying the linear exchange relation \eqref{E:LinearExchangeRelation} and such that $A_0=P$. We proceed in two steps: first we equip $A$ with a globular reflexive operad structure, then we equip it with a $\omega$-operad structure.

First, let $n\geq 0$. Equip $A_n$ with the partial multiplication $\mu_n^\part:=\mu_{A_n}^\uparrow=\mu_{A_n}^\downarrow$ and then define the operadic multiplication
\[
\mu_n:A_n\circ A_n\xrightarrow{\sigma_{A_n}}\bigoplus_{p\geq 1} A_n\part A_n^{\part p}\xrightarrow{\sum(\mu_n^\part)^p}A_n.
\]
The multiplication $\mu_n$ satisfies the associativity and unit axioms by Lemma \ref{L:3-implies-2}. Moreover, by construction, $\mu_n^\part$ satisfies the relations
\begin{eqn}{equation}
\label{E:GlobularOperadRelation}
\mu_n^\part=\mu_n^\part(1\part s_0)+\mu_n^\part(t_0\part 1)-\mu_n^\part(t_0\part s_0)=\mu_n^\part(s_0\part 1)+\mu_n^\part(1\part t_0)-\mu_n^\part(s_0\part t_0).
\end{eqn}
This gives $A$ a globular reflexive operad structure. Next, for the $\omega$-operad structure on $A$, we define the $\star_k$-compositions as follows: for all $\star_k$-composable $n$-cells $a,b$, define
\[
a\star_kb:=a-t_k(a)+b.
\]
Let $0\leq k<\ell<n$ be three integers. The target morphism $t_\ell:A_n\rightarrow A_\ell$ is linear, so it commutes with $\star_k$. For all $\star_\ell$-composable pairs $(a,a')$ and $(b,b')$ of $A_n$ such that $(a,b)$ and $(a',b')$ are $\star_k$-composable, we calculate
\begin{align*}
(a\star_\ell a')\star_k(b\star_\ell b') & = (a-t_\ell(a)+a')\star_k(b-t_\ell(b)+b') \\
& = a\star_kb-t_\ell(a)\star_kt_\ell(b)+a'\star_kb' \\
& = a\star_kb-t_\ell(a\star_kb)+a'\star_kb' \\
& = (a\star_kb)\star_\ell(a'\star_kb').
\end{align*}
Thus the $\star_k$-compositions satisfy exchange relations. To show that $A$ is an $\omega$-operad, it suffices to show that the $\star_k$-compositions are morphisms of operads. $A_n\star_kA_n$ is equipped with an operad structure given by the multiplication
\[
(A_n\star_kA_n)\circ(A_n\star_kA_n)\xrightarrow{(\pi_1\circ\pi_1)\star_k(\pi_2\circ\pi_2)}(A_n\circ A_n)\star_k(A_n\circ A_n)\xrightarrow{\mu_k\star_k\mu_k}A_n\star_kA_n,
\]
where $\pi_1,\pi_2$ are the projections of the fiber product $A_n\star_kA_n$. Therefore it suffices to check the exchange relation \eqref{E:MultiplicativeExchangeRelation}:
\[
\xymatrix @C=4em{
(A_n\star_kA_n)\circ(A_n\star_kA_n)
  \ar[rr]^-{(\pi_1\circ\pi_1)\star_k(\pi_2\circ\pi_2)}
  \ar[dr]_-{\mu_n((1\star_k1)\circ(1\star_k1))\qquad \quad }
&
&
(A_n\circ A_n)\star_k(A_n\circ A_n)
  \ar[dl]^-{\quad \qquad\mu_n(1\circ 1)\star_k\mu_n(1\circ 1)}
\\
&
A_n
}
\]
Writing
\[
\mu_n=\mu_n\varphi^{}_{A_n\star_kA_n}\sigma^{}_{A_n\star_kA_n}=\left(\sum_p(\mu_n^\part)^{\part p}\right)\sigma^{}_{A_n\star_kA_n},
\]
we get the diagram
\[
\xymatrix @C=4.8em{
(A_n\star_kA_n)\circ(A_n\star_kA_n)
  \ar[rr]^-{(\pi_1\circ\pi_1)\star_k(\pi_2\circ\pi_2)}
  \ar[d]_-{\sigma^{}_{A_n\star_kA_n}}
&
&
(A_n\circ A_n)\star_k(A_n\circ A_n)
  \ar[d]^-{\sigma^{}_{A_n}\star_k\sigma^{}_{A_n}}
\\
\bigoplus(A_n\star_kA_n)\part(A_n\star_kA_n)^{\part p}
  \ar[rr]^-{(\pi_1\part\pi_1^{\part p})\star_k(\pi_2\part\pi_2^{\part p})}
  \ar[dr]_-{\mu_n^\part((1\star_k1)\part(1\star_k1)^{\part p})\qquad \qquad \qquad}
&
&
\left(\bigoplus A_n\part A_n^{\part p}\right)\star_k\left(\bigoplus A_n\part A_n^{\part p}\right)
  \ar[dl]^-{\qquad \quad \qquad\mu_n^\part(1\part 1^{\part p})\star_k\mu_n^\part(1\part 1^{\part p})}
\\
&
A_n
}
\]
The upper square commutes by naturality of $\sigma$. To show that the lower triangle commutes, it suffices to show that, for all $\star_k$-composable pairs $(a,a')$ and $(b,b')$ of $A_n$ and all elementary compositions $\circ_{i,\tau}$ such that $a\circ_{i,\tau}b$ and $a'\circ_{i,\tau}b'$ are well defined,
\[
(a\star_k a')\circ_{i,\tau}(b\star_k b') = (a\circ_{i,\tau}b)\star_k(a'\circ_{i,\tau}b').
\]
Write $\cdot$ for $\circ_{i,\tau}$. Let us begin with the case $k=0$. By definition of $\star_k$-composition and bidistributivity of $\circ_{i,\tau}$, we have
\begin{align*}
(a+a'-t_0(a))\cdot(b+b'-t_0(b)) ={} & a\cdot b +a\cdot b' + a'\cdot b + a'\cdot b'-t_0(a)\cdot b-t_0(a)\cdot b' \\
& -a\cdot t_0(b)-a'\cdot t_0(b)+t_0(a)\cdot t_0(b)
\end{align*}
By applying \eqref{E:GlobularOperadRelation} to $a\cdot b'$ and $a'\cdot b$, we get
\begin{align*}
a\cdot b' & =t_0(a)\cdot b'+a\cdot s_0(b')-t_0(a)\cdot s_0(b'), \\
a'\cdot b & =s_0(a')\cdot b+a'\cdot t_0(b)-s_0(a')\cdot t_0(b).
\end{align*}
Since $t_0(a)=s_0(a')$ and $t_0(b)=s_0(b')$, we conclude that
\begin{align*}
(a+a'-t_0(a))\cdot(b+b'-t_0(b)) & =a\cdot b+a'\cdot b'-t_0(a)\cdot t_0(b) \\
& =(a\cdot b)\star_0(a'\cdot b).
\end{align*}
	
Now, let $k\geq 1$. In this case, $n\geq 2$, so by globularity, $t_0(a)=t_0(a')$ and $s_0(b)=s_0(b')$. Write $c:=t_0(a)=t_0(a')$ and $d:=s_0(b)=s_0(b')$ . Using the exchange relations between $\part$ and $\star_0$, and between $\star_0$ and $\star_k$, we get
\begin{align*}
(a\cdot b)\star_k(a'\cdot b') & =((a\star_0c)\cdot(d\star_0b))\star_k((a'\star_0c)\cdot(d\star_0b')) \\
& =((a\cdot d)\star_0(c\cdot b))\star_k((a'\cdot d)\star_0(c\cdot b')) \\
& =((a\cdot d)\star_k(a'\cdot d))\star_0((c\cdot b)\star_k(c\cdot b')).
\end{align*}
By definition of $\star_k$ and by bidistributivity of $\circ_{i,\tau}$, we get
\begin{align*}
(a\cdot d)\star_k(a'\cdot d) & =a\cdot d+a'\cdot d-t_k(a)\cdot d=(a\star_ka')\cdot d, \\
(c\cdot b)\star_k(c\cdot b') & =c\cdot b+c\cdot b'-c\cdot t_k(b)=c\cdot(b\star_kb').
\end{align*}
Thus
\begin{align*}
(a\cdot b)\star_k(a'\cdot b') & =((a\star_ka')\cdot d)\star_0(c\cdot(b\star_kb')) \\
& =(a\star_ka')\cdot(b\star_kb').
\end{align*}
Thus the exchange relation is satisfied, and we conclude that $A$ is an $\omega$-operad.
\end{proof}

\section{Shuffle operadic polygraphs}
\label{S:HigherDimensionalShufflePolygraphs}

In this section we introduce the notion of a shuffle polygraph that defines systems of generators and oriented relations for higher shuffle operads. 

\subsection{Shuffle polygraphs}
\label{SS:ShufflePolygraphs}

The structure of polygraph was introduced independently by Street and Burroni as a system of generators for free higher categories~\cite{Street76, Burroni93}. This subsection introduces a version of this structure for higher shuffle operads.

\subsubsection{Cellular extensions}
For $n\geq 0$, a \emph{cellular extension} of a shuffle $n$-operad $P$ is an indexed set $X$ equipped with two morphisms
\[
\xymatrix @C=2.5em {
P_n
&
X
  \ar@<-0.7ex> [l] _-{s_n}
  \ar@<0.7ex> [l] ^-{t_n}
}
\]
in $\Ind$ satisfying the globular relations $s_{n-1}s_n = s_{n-1}t_n$ and $t_{n-1}s_n = t_{n-1}t_n$, for $n\geq 1$, and whose elements are called \emph{$(n+1)$-generators}.
Note that every $n$-operad $P$ has two canonical cellular extensions: the empty one, and the one denoted by $\Sph(P_n)$ that consists of a $(n+1)$-generator $a\fl b$ for every $n$-sphere $(a,b)$ of $P$.

\subsubsection{Extended higher operads}
For $n\geq 0$, the category of \emph{extended $n$-operads} $\ShuO_n^+$ is defined by the following pullback of forgetful functors
\[
\xymatrix @C=2.5em{
\ShuO_n^+
 	\ar[r]
	\ar[d]	
	\pullbackcorner
&
\Glob_{n+1}(\Ind)
  \ar[d] ^-{\Vr_n}
\\
\ShuO_n
  \ar[r] _-{\Ur_n^\Sha}
&
\Glob_n(\Ind)
}
\]
where the functor $\Ur_n$ forgets the shuffle composition, and the functor $\Vr_n$ forgets the $(n+1)$-cells.
Explicitly, an extended shuffle $n$-operad is a pair $(P,X)$ where $P$ is an $n$-operad and $X$ a cellular extension of $P_n$. A morphism of extended $n$-operads $(P,X) \fl (P',X')$ is a morphism of shuffle $n$-operads $f : M \fl M'$ and a morphism $g:X \fl X'$ in $\Ind$ such that the two following square diagrams commute in $\Ind$:
\[
\xymatrix @C=3.3em @R=2.5em{
P_n
  \ar[d] _-{f_n}
&
X
  \ar[d] ^-{g}
  \ar@<-0.8ex>[l] _-{s_n}
  \ar@<+0.8ex>[l] ^-{t_n}
\\
P'_n
&
X'
  \ar@<-0.8ex>[l] _-{s_n}
  \ar@<+0.8ex>[l] ^-{t_n}
}
\]

\begin{proposition}
\label{P:free-n-operad}
Let $P$ be a shuffle $(n-1)$-operad and $X$ be a cellular extension of $P$. 
Let $\fopShu{X}$ denote the coequalizer of the two morphisms
\[
\xymatrix @C=2cm
{
(P_0\langle X\rangle\oplus P_{n-1})\part (P_0\langle X\rangle\oplus P_{n-1})
  \ar@<0.7ex>[r]^-{\mu_{P_0\langle X\rangle\oplus P_{n-1}}^\uparrow}
  \ar@<-0.7ex>[r]_-{\mu_{P_0\langle X\rangle\oplus P_{n-1}}^\downarrow}
&
P_0\langle X\rangle\oplus P_{n-1}
}
\]
in the category $\Bimod{P}$, where the morphisms are defined relative to the pair $(P_0,P_0\langle X\rangle\oplus P_{n-1})$.
Then $(P,\fopShu{X})$ is the free shuffle $n$-operad on $(P,X)$.
\end{proposition}
\begin{proof}
We will progressively enrich the cellular extension $X$ with more and more structure in order to get a reflexive globular bimodule satisfying the linear exchange relation \eqref{E:LinearExchangeRelation}. First, let us define categories of ``enriched cellular extensions'' of $P$. Let $\mathbb D$ be a category such that $P$ is a reflexive $(n-1)$-globular object of $\mathbb D$. Define
\begin{itemize}
\item $\Glob_P(\mathbb D)$ the full subcategory of $\Glob_n(\mathbb D)$ whose objects $X$ satisfy $(X_0,\ldots,X_{n-1})=P$, and
\item $\RGlob_P(\mathbb D)$ the full subcategory of $\RGlob_n(\mathbb D)$ whose objects $X$ satisfy $(X_0,\ldots,X_{n-1})=P$.
\end{itemize}
Finally, denote by $\RGlob_P^\part(\Bimod{P_0})$ the full subcategory of $\RGlob_P(\Bimod{P_0})$ whose objects satisfy the linear exchange relation \eqref{E:LinearExchangeRelation}.
	
Following Theorem \ref{T:MainTheoremA}, given an extended $(n-1)$-operad $(P,X)$, in order to construct the free $n$-operad, it suffices to construct the free object on $(P,X)\in\Glob_P(\Ind)$ in $\RGlob_P^\part(\Bimod{P_0})$. Therefore, it suffices to construct the sequence of free functors
\[
\Glob_P(\Ind) \fl \Glob_P(\Bimod{P_0}) \fl \RGlob_P(\Bimod{P_0}) \fl \RGlob_P^\part(\Bimod{P_0}).
\]
\begin{itemize}
\item Let $(P,X)$ be and extended $(n-1)$-operad. The first free functor is induced by the free functors $\Ind \fl \Coll \fl \Bimod{P_0}$, so it sends $X$ to $P_0\langle X\rangle$.
\item Let $(P,X)$ be an object of $\Glob_P(\Bimod{P_0})$. Since $P$ is already a reflexive $(n-1)$-globular object, the second free functor is induced by the free functor $\Bimod{P_0}\rightarrow P_{n-1}/\Bimod{P_0}$, so it sends $X$ to $X\oplus P_{n-1}$.
\item Let $(P,X)$ be an object of $\RGlob_P(\Bimod{P_0})$. The third free functor sends $X$ to the coequalizer of
\[
\xymatrix
{
X\part X
  \ar@<0.7ex>[r]^-{\mu_X^\uparrow}
  \ar@<-0.7ex>[r]_-{\mu_X^\downarrow}
&
X
}
\]
where the morphisms are defined relative to the pair $(P_0,X)$.
\end{itemize}
By composing these functors, we get
\[
\fopShu{X}:=\text{coeq}\,\big(\mu_{P_0\langle X\rangle\oplus P_{n-1}}^\uparrow,\mu_{P_0\langle X\rangle\oplus P_{n-1}}^\downarrow:(P_0\langle X\rangle\oplus P_{n-1})\part(P_0\langle X\rangle\oplus P_{n-1})\rightrightarrows P_0\langle X\rangle\oplus P_{n-1}\big)
\]
and we conclude that $(P,\fopShu{X})$ is the free $n$-operad on $(P,X)$.
\end{proof}

\subsubsection{Free shuffle $n$-operad}
For $n\geq 1$, the forgetful functor $\Wr_n^\Sha : \ShuO_n \fl \ShuO_{n-1}^+$
that forgets the composition of $n$-cells admits a left adjoint
\begin{eqn}{equation}
\label{E:FreeFunctorLn}
\Lr^\Sha_n :  \ShuO_{n-1}^+ \fl \ShuO_n
\end{eqn}
that associates to an extended $(n-1)$-operad $(P,X)$ the \emph{free $n$-operad over $(P,X)$} given by $\Lr^\Sha_n(P,X)=(P,\fopShu{X})$.
In the sequel, the $n$-operad $(P,\fopShu{X})$ will be denoted by $P[X]$, and its $k$-source and $k$-target maps will be denoted by $s_k$ and $t_k$ respectively.

\subsubsection{Shuffle polygraphs}
\label{SSS:ShufflePolygraphs}
We define the category $\ShuPol_n$ of \emph{$n$-polygraphs} and the free functor 
\[
\Fr_n : \ShuPol_n \fl \ShuO_n,
\]
by induction on $n\geq0$.
For $n=0$, we define $\ShuPol_0$ as the category $\Ind$. 
The free $0$-monoid functor
\[
\Fr_0 : \ShuPol_0 \fl \ShuO_0
\]
is the composite of free functors (\ref{E:FreeFunctorOperadFromSet}).
We suppose that for $n\geq 1$ the category $\ShuPol_{n-1}$ of $(n-1)$-polygraphs is defined and that the free $(n-1)$-operad functor
\[
\Fr_{n-1} : \ShuPol_{n-1} \fl \ShuO_{n-1}
\]
is constructed.
The category $\ShuPol_n$ of $n$-polygraphs is defined as the following pullback in $\Cat$
\begin{eqn}{equation}
\label{E:PullbackPolygraph}
\raisebox{0.8cm}{
\xymatrix @C=2.5em{
\ShuPol_n
  \ar[r] ^-{\widetilde{\Ur}_{n-1}}
  \ar[d] _-{\widetilde{\Vr}_{n-1}} 
  \pullbackcorner
&
\ShuO_{n-1}^+
  \ar[d] ^-{\Wr_n^\Sha}
\\
\ShuPol_{n-1}
  \ar[r] _-{\Fr_{n-1}}
&
\ShuO_{n-1}
}}
\end{eqn}
where the vertical functor on the right forgets the cellular extension of an extended monoid. The free symmetric $n$-operad functor is defined as the composite
\[
\xymatrix{
\ShuPol_n 
  \ar[r] ^-{\widetilde{\Ur}_{n-1}}
&
\ShuO_{n-1}^+ 
  \ar[r] ^-{\Lr^\Sha_n} 
&
\ShuO_n,
}
\]
where $\widetilde{\Ur}_{n-1}$ is the functor defined by the pullback~(\ref{E:PullbackPolygraph}) and $\Lr_n^\Sha$ is the free functor defined in~(\ref{E:FreeFunctorLn}).

The category $\ShuPol_{\omega}$ of $\omega$-polygraphs and the free $\omega$-operad functor $\Fr_{\omega} : \ShuPol_{\omega} \fl \ShuO_{\omega}$ are defined as the limit of the functors:
\[
\cdots 
\fl 
\ShuPol_{n} 
\ofl{\widetilde{\Vr}_{n-1}} 
\ShuPol_{n-1} 
\fl
\cdots
\fl
\ShuPol_1
\ofl{\widetilde{\Vr}_{0}} 
\ShuPol_0,
\]
in the category of categories, where the functors $\widetilde{\Vr}_{n-1}$ are defined by~(\ref{E:PullbackPolygraph}).

In this way, an $n$-polygraph $X$ is defined inductively as a data $(X_0,\ldots,X_n)$, where $X_0$ is an indexed set and for every $0<k <n$, $X_k$ is a cellular extension of the free $(k-1)$-operad generated by $(X_0,\ldots, X_{k-1})$, denoted by
\[
\fopShu{X}_{k} = \fopShu{X}_0[X_1]\cdots [X_k].
\]
For $0\leq p <n$, we will denote by $X_{\leq p}$ the underlying $p$-polygraph $(X_0,\ldots,X_p)$.

\subsubsection{Higher-dimensional monomials}

Let~$X$ be an $\omega$-polygraph. Tree monomials in $\ftree{X_0}$ are called \emph{$0$-monomials of $\fopShu{X}$}, and they form a linear basis of the collection~$\fopShu{X_0}$, which means that every $0$-cell~$a$ of~$\fopShu{X}_0$ can be uniquely written as a (possibly empty) linear combination
\[
a = \sum_{i=1}^p \lambda_i u_i + \lambda\varepsilon
\]
of pairwise distinct $0$-monomials~$u_1$, \dots, $u_p$ of~$\fopShu{X}_0$, with~$\lambda_i\in\K\setminus\{0\}$, $\lambda\in\K$, and $\varepsilon$ denotes the trivial monomial. This expression is called \emph{the canonical decomposition of~$a$}, and we define the \emph{support of~$a$} as the set $\supp(a)=\ens{u_1,\dots,u_p}$.

For~$n\geq 1$, if $\alpha$ is an $n$-cell of $X$, and $\vec v$ is a list of $0$-monomials, we will denote by $(\alpha\mid\vec v)$ the $n$-cell of $\fopShu{X}$ with source $(s(\alpha)\mid\vec v)$ and target $(t(\alpha)\mid\vec v)$. An \emph{$n$-monomial of~$\fopShu{X}$} is an $n$-cell of~$\fopShu{X}$ of the form~$u\circ_{i,\tau}(\alpha\mid\vec v)$, where~$\alpha$ is an $n$-cell of~$X$, and~$u$ and~$\vec v$ are monomials of~$\fopShu{X}$. By construction of the free $n$-operad over~$(\fopShu{X_{n-1}},X_n)$, and by freeness of~$\fopShu{X_{n-1}}$, every $n$-cell~$a$ of~$\fopShu{X}$ can be written as a linear combination
\begin{eqn}{equation}
\label{E:DecompositionNCell}
a = \sum_{i=1}^p \lambda_i a_i + 1_c
\end{eqn}
of pairwise distinct $n$-monomials~$a_1$, \dots, $a_p$ and of an identity $n$-cell~$1_c$ of~$\fopShu{X}$, and this decomposition is unique up to the linear exchange relation (\ref{E:LinearExchangeRelation}).
The \emph{size} of an $n$-cell~$a$ of~$\fopShu{X}$ is the minimal number of $n$-monomials of~$\fopShu{X}$ required to write~$a$ as in~\eqref{E:DecompositionNCell}.

\subsubsection{Graded shuffle polygraps}
In order to define in Section~\ref{S:ShufflePolygraphicResolutionConvergence} minimality and Koszulness properties with respect to a polygraphic resolution, we introduce the notion of a graded shuffle $\omega$-polygraph. Just as we defined shuffle operads as internal monoids in the presheaf category $\Vect^{\op{\Ord}}$ in Section~\ref{SS:ShuffleOperads}, we define \emph{graded shuffle operads} as internal monoids in the presheaf category $\grVect^{\op{\Ord}}$. For $n\in \Nb\cup\{\omega\}$, a \emph{graded shuffle $n$-operad} is an $n$-category in $\grShuO$, and we denote by $\grShuO_n$ the corresponding category with internal $n$-functors as morphisms. In particular, the source, target and composition morphisms of graded shuffle $n$-operad are graded.

The category $\grShuO^+_n$ of \emph{graded extended $n$-operads} is defined similarly to $\ShuO^+_n$: its objects are pairs~$(P,X)$, where~$P$ is a graded $n$-operad, and~$X$ is a graded cellular extension of~$A$, meaning that $X = \amalg_{i\geq 0} X^{(i)}$ and that the source and target of each~$x$ in~$X^{(i)}$ are homogeneous of degree~$i$. In that case, the free $(n+1)$-operad~$P[X]$, defined as in the nongraded case, is also graded.

A \emph{graded $\omega$-polygraph} is an $\omega$-polygraph~$X$ such that each set~$X_n$ is graded, for $n\geq 0$. This notion restricts to $n$-polygraphs, and a $1$-polygraph~$X$ is called \emph{quadratic} if~$X_0$ is concentrated in degree~$1$ and~$X_1$ is concentrated in degree~$2$.
We say that a graded $\omega$-polygraph~$X$ is \emph{concentrated on the superdiagonal} if each graded set~$X_n$, for~$n\geq 0$, is concentrated in degree~$n+1$. In that case, because the source and target maps are graded, for~$n\geq 1$, the source and target of every $n$-cell of~$\fopShu{X}$ are homogeneous $(n-1)$-cells of~$\fopShu{X}$ of degree~$n+1$.

\subsection{Shuffle polygraphic resolutions}

In this subsection we introduce the notion of a polygraphic resolution for shuffle operads.

\subsubsection{Presentation of a shuffle operad}
\label{SSS:PresentationShuffleOperads}
The shuffle operad \emph{presented} by a shuffle $1$-polygraph~$X$ is the coequalizer in the category $\ShuO$ of the following source and target morphisms, denoted by~$\cl{X}$,
\begin{eqn}{equation}
\label{E:presentationShuffleOperad}
\xymatrix@C=1.8em{
\fopShu{X}_1 
  \ar@<+0.7ex> [r]^-{s_0} 
  \ar@<-0.7ex> [r]_-{t_0}
&
\fopShu{X_{0}} 
  \ar[r] ^-{\pi_X}
&
\cl{X}.
}
\end{eqn}
Following Proposition~(\ref{P:reflectionLimits}), the category $\ShuO$ preserves reflexive coequalizers and so the construction is well defined. We say that a shuffle operad~$P$ is \emph{presented by a polygraph~$X$}, or that~$X$ is a \emph{presentation of~$P$}, if~$P$ is isomorphic to~$\cl{X}$ in the category $\ShuO$. 

\subsubsection{Shuffle polygraphic resolutions}
\label{SSS:ShufflePolygraphicResolutions}
For $n\geq 1$, let $P$ be a shuffle $n$-operad. A cellular extension $X$ of $P$ is \emph{acyclic} if for every $n$-sphere $(f,g)$ in $P_n$, there exists a shuffle $(n+1)$-cell $F$ in the $(n+1)$-operad $P[X]$ with source $f$ and target $g$.
A \emph{coherent presentation} of a shuffle polygraph $P$ is a shuffle $2$-polygraph $X$, that presents $P$, and whose cellular extension $X_2$ is acyclic.
A \emph{shuffle} \emph{polygraphic resolution} of a shuffle operad $P$ is a shuffle $\omega$-polygraph~$X$ that presents $P$ and whose cellular extensions $X_k$ are acyclic for~$k\geq 2$. 

\subsubsection{Example: standard polygraphic resolution}
Let $P$ be a shuffle operad, and $B$ a basis of $P$ seen as a collection. We then define the \emph{standard polygraphic resolution} $\Std(P)$ by induction as follows. For $n=0$, we define the indexed set $\Std(P)_0 := B$, and for $u\in B$ we denote by $[u]$ the corresponding element in $\Std(P)_0$. Any element $a=\sum_{u\in B}\lambda_uu$ of $P$ then corresponds to a linear combination of elements $[a]:=\sum_{u\in B}\lambda_u[u]$ in $\fopShu{\Std(P)_0}$. Next, for $n=1$, we set
\[
\Std(P)_1 := \{\ [u]\circ_{i,\tau} [v] \fl [u\circ_{i,\tau} v] \; \mid\; u,v\in B\ \},
\]
so that the pair $(B,\Std(P)_1)$ forms a $1$-polygraph that presents the shuffle operad~$P$.
Now, suppose that $\Std(P)_n$ is defined for $n\geq 1$. Then we set $\Std(P)_{n+1}:=\Sph(\fopShu{\Std(P)_n})$. By construction, the $\omega$-polygraph $\Std(P)$ is a polygraphic resolution of the shuffle operad $P$.

\subsubsection{Tietze equivalence of $\omega$-polygraphs}
\label{SSS:TietzeEquivalence}
We define the notion of a \emph{weak-equivalence} of $\omega$-operads as for $\omega$-categories, defined in \cite{LafontMetayerWorytkiewicz10}.
For $n\geq 0$, two $n$-cells $a,b$ of an $\omega$-operad $P$ are \emph{$\omega$-equivalent} if there exists an $(n+1)$-cell $f:a\fl b$ in $P$. In that case, we write $a\sim_\omega b$.
A morphism of $\omega$-operads $F:P\fl Q$ is a \emph{weak equivalence} if it satisfies the following properties:
\begin{enumerate}[{\bf i)}]
\item For every $0$-cell $a$ of $Q$, there exists a $0$-cell $\widehat a$ in $P$ such that $F(\widehat a)\sim_\omega a$.
\item For every pair of $0$-cells $a,b$ of $P$ and every $1$-cell $f:F(a)\fl F(b)$ of $Q$, there exists a $1$-cell $\widehat f:a\fl b$ of $P$ such that $F(\widehat f)\sim_\omega f$.
\item For $n\geq 1$ and every pair of parallel $n$-cells $a,b$ of $P$ and every $(n+1)$-cell $f:F(a)\fl F(b)$ of~$Q$, there exists an $(n+1)$-cell $\widehat f:a\fl b$ of $P$ such that $F(\widehat f)\sim_\omega f$.
\end{enumerate}

We say that two shuffle $\omega$-polygraphs $X$ and $Y$ are \emph{Tietze equivalent} if the induced free $\omega$-operads $\fopShu{X}$ and $\fopShu{Y}$ are weakly equivalent.
The original notion of Tietze equivalence for $1$-polygraphs is a particular case of this notion for $\omega$-polygraphs. Two $1$-polygraphs $X$ and $Y$ are \emph{Tietze equivalent} if the presented shuffle operads $\cl{X}$ and $\cl{Y}$ are isomorphic. In that case, extending $X$ and $Y$ into $\omega$-polygraphs with identities in higher dimensions gives two Tietze equivalent $\omega$-polygraphs. Tietze equivalence also generalizes the notion of Tietze equivalence between $(3,1)$-polygraphs introduced in \cite{GaussentGuiraudMalbos15}.

\begin{proposition}
\label{P:HigherTietzeEquivalence}
Let $X$ and $Y$ be two acyclic $\omega$-polygraphs. Then the presented shuffle operads $\cl X$ and~$\cl Y$ are isomorphic if, and only if, the polygraphs $X$ and $Y$ are Tietze equivalent.
\end{proposition}
\begin{proof}
$(\Rightarrow)$ Denote by $\varphi:\cl X\fl \cl Y$ the isomorphism. We define a morphism of $\omega$-operads $F:\fopShu{X}\fl\fopShu{Y}$ such that $\pi_YF=\varphi\pi_X$ on $0$-cells and show that it is a weak equivalence simultaneously. Since $\fopShu{X}$ is a free $\omega$-operad, it suffices to define $F$ on the $n$-generators of $X$ for all $n\geq 0$. We proceed by induction on $n\geq 0$.

For $n=0$, define linear maps $i_X:\cl X\fl\fopShu{X}$ and $i_Y:\cl Y\fl\fopShu{Y}$, which are sections of $\pi_X$ and $\pi_Y$, respectively. For $x\in X_0$, we set $F(x_0):=i_Y(\varphi(\pi_X(x_0)))$, and we check that $\pi_YF=\varphi\pi_X$ on $0$-cells. Now, for a $0$-cell $a$ of $\fopShu{Y}$, let $\hat{a}:=\varphi^{-1}i_Y(a)\in\fopShu{X}$. Then $\pi_Y(F(\hat{a}))=\pi_Y(a)$, so $F(\hat{a})\sim_\omega a$.

For $n=1$, for $\alpha:a\fl b$ a $1$-generator of $X$, $\pi_Y(F(a))=\pi_Y(F(b))$, so there exists a $1$-cell $f:F(a)\fl F(b)$ in $\fopShu{Y}$. We set $F(\alpha):= f$. Then, for every pair of $0$-cells $a,b$ of $\fopShu{X}$ and every $1$-cell $f:F(a)\fl F(b)$ of $\fopShu{Y}$, $\pi_Y(F(a))=\pi_Y(F(b))$, which is equivalent to $\pi_X(a)=\pi_X(b)$ via the isomorphism $\varphi$. Therefore there exists $\hat f:a\fl b$ in $\fopShu{X}$, and $F(\hat f):F(a)\fl F(b)$ is parallel to $f$. Since $Y$ is acyclic, $F(\hat f)\sim_\omega f$.

Let $n\geq 1$ and suppose that $F$ is defined on $n$-cells of $\fopShu{X}$. For $\alpha:a\fl b$ an $(n+1)$-generator of $X$, the $n$-cells $a$ and $b$ of $\fopShu{X}$ are parallel, so the $n$-cells $F(a)$ and $F(b)$ of $\fopShu{Y}$ are parallel. By acyclicity of $Y$, there exists an $(n+1)$-cell $f:F(a)\fl F(b)$ in $\fopShu{Y}$, so we set $F(\alpha):=f$. Now, let $a,b$ be two parallel $n$-cells of $\fopShu{X}$ and $f:F(a)\fl F(b)$ an $(n+1)$-cell of $\fopShu{Y}$. By acyclicity of $X$, there exists an $(n+1)$-cell $\hat f:a\fl b$, so $F(\hat f)$ and $f$ are parallel in $\fopShu{Y}$, so by acyclicity of $Y$ we have $F(\hat f)\sim_\omega f$.

We conclude that $F:\fopShu{X}\fl\fopShu{Y}$ is a weak equivalence, so $X$ and $Y$ are Tietze equivalent.

$(\Leftarrow)$ Let $F:\fopShu{X}\fl\fopShu{Y}$ be a weak equivalence. By condition {\bf i)}, $\pi_YF:\fopShu{X}\fl\cl Y$ is surjective. Moreover, if $a\sim_\omega b$ in $\fopShu{X}$, then $\pi_YF(a)=\pi_YF(b)$, so $F$ induces a morphism of $\omega$-operads $\cl F:\cl X\fl\cl Y$. By condition {\bf ii)}, $\cl F$ is injective. Thus $\cl F$ is an isomorphism between $\cl X$ and $\cl Y$.
\end{proof}

\subsubsection{Reduced polygraphs}
Let $X$ be a left-monomial $1$-polygraph. Recall from \cite{Polybook21} that a $1$-generator $\alpha\in X_1$ is \emph{right} (resp.\ \emph{left}) \emph{reduced} if $t_0(\alpha)\in\Red(X_1)$ (resp.\ $s_0(\alpha)\in\Red(X_1\setminus\{\alpha\})$). We say that $X$ is  \emph{reduced} when each of its $1$-generators is left and right reduced. We prove that every (finite) convergent left-monomial $1$-polygraph is Tietze-equivalent to a reduced (finite) convergent left-monomial $1$-polygraph. 

\section{Shuffle operadic rewriting}
\label{S:ShuffleOperadicRewriting}

The first part of this section presents the main rewriting properties of shuffle $1$-polygraphs. We relate the notion of a convergent shuffle polygraph, whose $1$-generators are oriented with respect to a given monomial order, with the notion of Gröbner bases introduced in \cite{DotsenkoKhoroshkin10}, and with the notion of Poincaré-Birkhoff-Witt bases introduced in \cite{Hoffbeck10}. Throughout this section, all operads and polygraphs are shuffle.

\subsection{Rewriting in shuffle operads}

We introduce a concept of rewriting in the context of shuffle operads.

\subsubsection{The terminal indexed set}
Denote by $\Box$ the terminal object of $\Ind$, that is, the indexed set that is a singleton $\Box(k)=\{\Box_k\}$ for each arity $k\geq 1$. 
Denote by $\iota_k:\Set\fl \Ind$ the inclusion functor defined by $\iota_k(X_0)(k) = X_0$ and $\iota_k(X_0)=\varnothing$ for the other arities.  

\subsubsection{One-hole contexts of indexed sets}
A \emph{one-hole context} of an indexed set $X_0$ is an element $\Gamma$ of the free $\fopShu{X_0}$-bimodule $\fopShu{X_0}\langle\Box\rangle$. We say that $\Gamma$ is of \emph{inner arity} $k$ if it is an element of $\fopShu{X_0}\langle\Box(k)\rangle$.

Let $A$ be an $\fopShu{X_0}$-bimodule and $a\in A(k)$. Identifying $A(k)$ with $\Hom_{\Set}(\Box(k),A(k))$, $a$ induces a morphism
\[
\varphi_a:\fopShu{X_0}\langle\Box(k)\rangle \fl \fopShu{X_0}\langle A(k)\rangle.
\]
via the the functor $\fopShu{X_0}\langle\iota_k\rangle:\Set\fl\Bimod{\fopShu{X_0}}$. For $\Gamma$ a one-hole context of $\ftree{X}_0$ of inner arity $k$, we write $\Gamma[a]:=\varphi_a(\Gamma)$. Explicitly, $\Gamma[a]$ is a tree of the form
\[
\begin{forest}
[u
  [v_1]
  [\cdots, no edge]
  [a, l=8ex, edge label={node[midway,fill=white,font=\scriptsize]{$i$}}
    [w_1]
    [\cdots, no edge]
    [w_k]
  ]
  [\cdots, no edge]
  [v_n]
]
\end{forest}
\]
where $k,n\geq 1$, $i\in\range{n}$, $u\in \fopShu{X}_0(n)$, $v_1,\ldots,\check v_i,\ldots,v_n,w_1,\ldots,w_k\in \fopShu{X}_0$ and $\Box_k$ appears in the $i^\text{th}$ position.  The notation $\check v_i$ means that we omit $v_i$.
In this way, every one-hole context $\Gamma$ on $X_0$ can be written $\Gamma:=w\circ_{i,\tau}(\Box_k\mid\vec w)$ with $w,\vec w\in\fopShu{X_0}$.
In this work, we will only consider \emph{monomial} one-hole contexts, that is when $w,\vec w$ are monomials of $\ftree{X_0}$, so we will omit the word monomial.

\subsubsection{Two-hole contexts of indexed sets}
Let $X_0$ be an indexed set. There exists a bifunctor \linebreak $C_2^{X_0}:\Ind\times\Ind \fl \Ind$ which sends a pair of indexed sets $Y,Y'$ to the indexed set of elements of $\ftree{(X_0\sqcup Y\sqcup Y')}$ with one occurrence of both $Y$ and $Y'$.
A \emph{two-hole context of $\ftree{X}_0$} is an element $\Gamma$ of $C_2^{X_0}(\Box,\Box)$. 
We say that $\Gamma$ is of \emph{inner arities} $(k,\ell)$ if it is an element of $C_2^{X_0}(\Box(k),\Box(\ell))$. 

Let $P$ be an operad equipped with a morphism $\pi:\fopShu{X_0}\fl P$, and $a\in P(k),a'\in P(\ell)$. Identifying
\begin{align*}
P(k)\times P(\ell) & \simeq \Hom_{\Set}(\Box(k),P(k))\times\Hom_{\Set}(\Box(\ell),P(\ell)) \\
& \simeq\Hom_{\Set\times\Set}((\Box(k),\Box(\ell)),(P(k),P(\ell))),
\end{align*}
the pair $(a,a')$ induces a morphism
\[
\varphi_{a,a'}:C_2^{X_0}(\Box(k),\Box(\ell)) \fl C_2^{X_0}(P(k),P(\ell))
\]
via the bifunctor $C_2^{X_0}(\iota_k,\iota_\ell):\Set\times\Set \fl \Ind$. Moreover, $f$ induces a morphism $\pi_*:C_2^{X_0}(P(k),P(\ell))\fl P$. For $\Gamma$ a two-hole context of $\ftree{X}_0$ of inner arities $(k,\ell)$, we write $\Gamma(a,a'):=i_*\varphi_{a,a'}(\Gamma)$. Explicitly, $\Gamma(a,a')$ is a tree of one of the following two forms, where the application of $i$ is implicit:
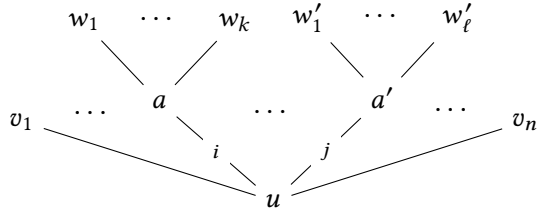
\begin{enumerate}
\item
\[
\begin{forest}
[u
  [v_1]
  [\cdots, no edge, before computing xy={l/.average={l}{siblings}}]
  [a, l=8ex, edge label={node[midway,fill=white,font=\scriptsize]{$i$}}
    [w_1]
    [\cdots, no edge]
    [w_k]
  ]
  [\cdots, no edge, before computing xy={s/.average={s}{siblings}, l/.average={l}{siblings}}]
  [a', l=8ex, edge label={node[midway,fill=white,font=\scriptsize]{$j$}}
    [w'_1]
    [\cdots, no edge]
    [w'_\ell]
  ]
  [\cdots, no edge, before computing xy={l/.average={l}{siblings}}]
  [v_n]
]
\end{forest}
\]
where
\begin{gather*}
n\geq 2,\quad i,j\in\range{n},\quad u\in \fopShu{X}_0(n),\text{ and} \\
v_1,\ldots,\check v_i,\ldots,\check v_j,\ldots,v_n,w_1,\ldots,w_k,w'_1,\ldots,w'_\ell\in \fopShu{X}_0,
\end{gather*}
\item
\[
\begin{forest}
[u
  [v_1]
  [\cdots, no edge]
  [a, l=8ex, edge label={node[midway,fill=white,font=\scriptsize]{$i$}}
    [x_1]
    [\cdots, no edge]
    [x_h
      [x'_1]
      [\cdots, no edge]
      [a', l=8ex, edge label={node[midway,fill=white,font=\scriptsize]{$j$}}
        [w_1]
        [\cdots, no edge]
        [w_\ell]
      ]
      [\cdots, no edge]
      [x'_m]
    ]
    [\cdots, no edge]
    [x_k]
  ]
  [\cdots,no edge]
  [v_n]
]
\end{forest}
\]
where 
\begin{gather*}
m,n\geq 1,\quad i\in\range{n},\quad j\in\range{m},\quad h\in\range{k},\text{ and} \\ 
u,v_1,\ldots,\check v_i,\ldots,v_n,w_1,\ldots,w_\ell,x_1,\ldots,x_k,x'_1,\ldots,\check x'_j,\ldots,x'_m\in \fopShu{X}_0.
\end{gather*}
\end{enumerate}
In this work, we will only consider \emph{monomial} two-hole contexts, that is, when in {\bf i)} and {\bf ii)} the $u$, $\vec v$, $\vec x$, $\vec x'$, $\vec w$, $\vec w'$ are monomials of $\ftree{X_0}$.

\subsubsection{Left-monomiality and homogeneity}
A cellular extension $X_1$ of $\fopShu{X}_0$ is \emph{left-monomial} if, for every $1$-generator $\alpha$ in $X_1$ the source $s_0(\alpha)$ is a non-trivial $0$-monomial, and $s_0(\alpha)\notin\supp(t_0(\alpha))$. A $1$-polygraph is \emph{left-monomial} if $X_1$ is so. We prove that every $1$-polygraph is Tietze equivalent to a left-monomial one.
For $N\geq 1$, a cellular extension $X_1$ of $\fopShu{X}_0$ is \emph{homogeneous} if, for every $1$-generator $\alpha$ in $X_1$ the weight of $s_0(\alpha)$ and  $t_0(\alpha)$ are equal to $N$. A $1$-polygraph is \emph{$N$-homogeneous} if $X_1$ is so. When $N=2$ we say \emph{quadratic} for $N$-homogeneous.

\subsubsection{Rewriting step}
A \emph{rewriting step} of a left-monomial $1$-polygraph $X$ is a $1$-cell $f$ of $\fopShu{X}_1$ of size $1$ of the form
\[
f = \lambda g + 1_c,
\]
where $\lambda\in\K\setminus\{0\}$, $g$ is a $1$-monomial of $\fopShu{X}_1$, and $c$ is a $0$-cell of $\fopShu{X}_0$ such that the $0$-monomial $s_0(u)\notin \supp(c)$. A $1$-cell of $\fopShu{X}_1$ if \emph{positive} if it is the $\star_0$-composition of rewriting steps.

A $0$-cell $a$ of $\fopShu{X}_0$ is \emph{reduced} if there is no rewriting step with source $a$. We denote by $\Red(X)$ the indexed submodule of reduced $0$-cells. The indexed set $\Red_m(X)$ of reduced $0$-monomials of $\fopShu{X}_0$ forms a basis of $\Red(X)$. A \emph{normal form} of $a$ is a reduced $0$-cell $b$ such that there is a positive $1$-cell with source $a$ and target $b$.

\subsubsection{Monomial orders and termination}
\label{termination}
An \emph{indexed poset} $(A,\prec)$ is an indexed set $A$, such that each $A(k)$ is equipped with a partial order $\prec_k$ ; we will omit the index on $\prec$. An indexed poset $(A,\prec)$ is \emph{well-founded} if each $A(k)$ is a well-founded poset. 

Let $X_0$ be an indexed set. An order relation $\prec$ on the free monoid $\ftree{X_0}$ of tree monomials is \emph{stable by product} if, for all $u,u'\in \ftree{X_0}(k)$, $v,v'\in \ftree{X_0}(\ell)$, $i\in\range{k}$, and $\tau\in\Sha(\ell-1,k-i)$, $u\prec u',v\prec v'$ implies $u\circ_{i,\tau}v\prec u'\circ_{i,\tau}v'$.
A total order relation stable by product is called a \emph{monomial order} on $\ftree{X_0}$. Note that this notion also appears in \cite{Hoffbeck10} and \cite{DotsenkoKhoroshkin10}.

For $Y$ a left-monomial cellular extension of $\fopShu{X}_0$, an order relation $\prec$ on $\ftree{X_0}$ is \emph{compatible with $Y$} if, for every $1$-cell $\alpha:u\fl a$ of $Y$ and every monomial $v\in\supp(a)$, $v\prec u$.
The relation $\prec$ can be extended to the free shuffle operad $\fopShu{X}_0$ as follows: for two $0$-cells $a,b$ of $\fopShu{X}_0$, we have $b\prec a$ if the two following conditions are satisfied
\begin{enumerate}
\item $\supp(a)\setminus\supp(b)\neq\varnothing$,
\item for all $v\in\supp(b)\setminus\supp(a)$, there exists $u\in\supp(a)\setminus\supp(b)$ such that $v\prec u$.
\end{enumerate}

For a left-monomial $1$-polygraph $X$, we denote by $\prec_{X_1}$ the smallest partial order relation on $\ftree{X_0}$ stable by product and compatible with $X_1$. 
A $1$-polygraph $X$ is \emph{terminating} if the relation $\prec_{X_1}$ is well-founded.
In that case, for every rewriting step $f$ of $X$, we have $t_0(f) \prec_{X_1} s_0(f)$, and thus there does not exist infinite sequence of rewriting steps of $X$.

\subsubsection{Path-lexicographic order on $1$-monomials}
\label{SSS:PathLexicographicMonomialOrder}
Let $(X_0,\prec)$ be a totally ordered indexed set, and~$\prec_{\SS}$ a total order on permutations. Let us recall from \cite{DotsenkoKhoroshkin10,Hoffbeck10} the \emph{path-lexicographic} monomial order on~$\ftree{X}_0$.
Given a $0$-monomial $u$ of arity $k$, there exists a unique path from its root to each of its inputs. Write such a path as a word $a=x_1\ldots x_n$ in the alphabet $X_0$ made of labels of the vertices of the path. To the $0$-monomial $u$ we associate the pair $(L, f)$ where $L$ is the sequence $(p_1,\ldots,p_k)$ of path words from the root to the inputs of $u$, and $f$ is the permutation of the inputs of $u$.
Then we define the \emph{path-lexicographic order} induced by the orders $\prec$ and $\prec_\SS$, denoted by $\prec_{pl}$, by setting
\[
(L, f) \prec_{pl} (L', f') 
\;\;\text{if}\;\; \big(L \prec_{lex} L'
\;\;\text{or}\;\; \big(L = L'\;\;\text{and}\;\;f\prec_{\SS} f' \big)\big),
\]
where $\prec_{lex}$ denotes the lexicographic order on words on $X_0$ induced by $\prec$.

Now, let $X$ be a left-monomial $1$-polygraph. Consider a total order $\prec$ on $X_0\sqcup X_1$, such that $\alpha\prec x$, for every $x\in X_0$ and $\alpha\in X_1$,  and $\prec_{\SS}$ a total order on permutations.
We will denote by $\prec_{pl}$ the induced path-lexicographic order on $\ftree{(X_0\sqcup X_1)}$, which induces a \emph{path-lexicographic order on $1$-monomials}.

\begin{proposition}
\label{P:WellFoundedMonomialTerminating}
Let $X$ be a left-monomial $1$-polygraph. If $\ftree{X_0}$ admits a well-founded monomial order  $\prec$ compatible with $X_1$, then $X$ is terminating.
\end{proposition}
\begin{proof}
We have $\prec_{X_1}\subseteq\prec$, so $\prec_{X_1}$ is well-founded, so $X$ is terminating.
\end{proof}

However, the converse implication is not true. In general, in order to prove termination when no monomial order is known, it is necessary to use a proof strategy appropriated to the set of rules. The following gives an illustration for one of the simplest strategies.

\begin{proposition}
\label{P:TerminationProofWithoutMonomialOrder}
A left-monomial $1$-polygraph $X$ terminates if, and only if, there exists a well-founded indexed poset $(W,<)$ and a morphism of indexed sets $\Phi:\ftree{X_0}\fl W$ such that $\Phi(\Gamma[v])<\Phi(\Gamma[s(\alpha)])$ holds for every $1$-generator $\alpha\in X_1$, one-hole context $\Gamma$, and $v\in\supp(t(\alpha))$.
\end{proposition}
\begin{proof}
Suppose that the polygraph $X$ terminates. Then $\ftree{X_0}$ is equipped with a well-founded partial order $\prec_{X_1}$, and We set $\Phi$ to be the identity morphism on $\ftree{X_0}$.

Conversely, let $\prec$ be the partial order generated by $v\prec u$ if there exists a rewriting rule $\alpha\in X_1$ and a one-hole context $\Gamma$ such that $u=\Gamma[s(\alpha)]$ and $v\in\supp(\Gamma[t(\alpha)])$. The order $\prec$ is stable by product by considering
\[
u\circ_{i,\tau} v \;\prec\; u'\circ_{i,\tau}v \;\prec\; u'\circ_{i,\tau}v',
\]
for all tree monomials $u,v,u',v'$ in $\ftree{X_0}$, and is compatible with $X_1$ by definition. Thus $\prec_{X_1}\subseteq\prec$, and so the map $\Phi:(\ftree{X_0},\prec_{X_1})\fl W$ is a strictly monotone morphism of indexed posets. Since $W$ is well-founded, $(\ftree{X_0},\prec_{X_1})$ is as well, and so the $1$-polygraph $X$ terminates.
\end{proof}

\subsubsection{Example}
\label{E:xyz}
We consider the polygraph $X$ with three $0$-generators $x,y,z$ and the following $1$-generator:
\[
\xymatrix{
\alpha:\;
\inputtree{[x[y[1][2]][z[3][4]]]}
\ar[r]
&
\inputtree{[x[x[1][2]][x[3][4]]]} + \inputtree{[y[y[1][2]][y[3][4]]]} + \inputtree{[z[z[1][2]][z[3][4]]]}
}
\]
For $u\in \ftree{X_0}$, we set $\Phi(u):=|u|_x + 3|u|_{y-z}$,
where $|u|_x$ denotes the number of occurrences of $x$ in $T(u)$ and $|u|_{y-z}$ the number of inner vertices of $T(u)$ whose two children are, from left to right, $y$ and $z$.

Then $\Phi(\Gamma[s(\alpha)])>\Phi(\Gamma[v])$ for all contexts $\Gamma$ of inner arity $4$ and every $v\in\supp(t(\alpha))$. Indeed, for every $\Gamma = w\circ_{i,\tau}(\Box_4\mid w_1\,w_2\,w_3\,w_4)$, we have
\begin{align*}
\Phi\left(\tree{[w[x[y[w_1][w_2]][z[w_3][w_4]]]]} - \tree{[w[x[x[w_1][w_2]][x[w_3][w_4]]]]}\right)
&
= \left|\tree{[x[y][z]]} - \tree{[x[y][z]]}\right|_x + 3\left|\tree{[x[y][z]]} - \tree{[x[x][x]]}\right|_{y-z}
\\
&
= (1-1) + (3-0) = 1,
\\
\Phi\left(\tree{[w[x[y[w_1][w_2]][z[w_3][w_4]]]]} - \tree{[w[y[y[w_1][w_2]][y[w_3][w_4]]]]}\right)
&
= \left|\tree{[x[y][z]]}-\tree{[y[y][y]]}\right|_x + 3\left|\tree{[w[x[y][z]]]} - \tree{[w[y[y][y]]]}\right|_{y-z}
\\
&
= 1 + 3\left|\tree{[x[y][z]]}\right|_{y-z} + 3|w|_{y-z} - 3|w\circ_{i,\tau}y|
\\
&
=\begin{cases}
4 & \text{if }|w\circ_{i,\tau}y|_{y-z}=|w|_{y-z}, \\
1 & \text{if }|w\circ_{i,\tau}y|_{y-z}=|w|_{y-z}+1,
\end{cases}
\\
\Phi\left(\tree{[w[x[y[w_1][w_2]][z[w_3][w_4]]]]} - \tree{[w[z[z[w_1][w_2]][z[w_3][w_4]]]]}\right)
&
=\begin{cases}
4 & \text{if }|w\circ_{i,\tau}z|_{y-z}=|w|_{y-z}, \\
1 & \text{if }|w\circ_{i,\tau}z|_{y-z}=|w|_{y-z}+1.
\end{cases}
\end{align*}
Following Proposition~\ref{P:TerminationProofWithoutMonomialOrder} the polygraph $X$ terminates. Note that, there is no monomial order that orients this rule in this way. Indeed, every orientation compatible with a monomial order reduces first one of the term of right hand side.

\subsection{Confluence of shuffle polygraphs}

In this subsection we define and algebraically characterize the property of confluence of a shuffle polygraph. We prove the coherent critical branching theorem for shuffle polygraphs involving a restricted notion of critical branchings. The definitions and results of this section do not differ much from the case of associative algebras in \cite{GuiraudHoffbeckMalbos19}, and indeed associative algebras can be seen as operads concentrated in arity $1$, but since the notion of contexts is not made explicit in previous works, we restate all of the definitions and results.

\subsubsection{Branchings}
A \emph{branching} of a left-monomial $1$-polygraph $X$ is a pair $(f,g)$ of positive $1$-cells of $\fopShu{X}_1$ where $f$ and $g$ have the same source $s_0(f)=s_0(g)$, which we denote by $s_0(f,g)$. The branching~$(f,g)$ is said to be \emph{local} if $f$ and $g$ are both rewriting steps.

Let~$X$ be a $2$-polygraph such that $X_{\leq 1}$ is left-monomial. A branching~$(f,g)$ of the $1$-polygraph~$X_{\leq 1}$ is \emph{$X_2$-coherently confluent}, or $(f,g)$ is coherently confluent for short, if there exist positive $1$-cells~$h$ and~$k$ of~$\fopShu{X}_1$ and a $2$-cell~$F$ of~$\fopShu{X}_2$ as in the following diagram
\[
\vcenter{\xymatrix @R=0.6em{
& t_0(f)
	\ar@/^/ [dr] ^-{h}
	\ar@2 []!<0pt,-12.5pt>;[dd]!<0pt,12.5pt> ^-*+{F}
\\
s_0(f)=s_0(g)
	\ar@/^/ [ur] ^-{f}
	\ar@/_/ [dr] _-{g}
&&
t_0(h)=t_0(k)
\\
&
t_0(g)
	\ar@/_/ [ur] _-{k}
}}
\]
If~$u$ is a $0$-cell of~$\fopShu{X}_0$, the $2$-polygraph~$X$ is \emph{coherently confluent} (resp.\ \emph{locally coherently confluent}) \emph{at~$u$} if every branching (resp.\ local branching) of~$X$ of source~$u$ is coherently confluent. 
The $2$-polygraph~$X$ is \emph{coherently confluent} (resp.\ \emph{locally coherently confluent}) if it is so at every $0$-cell of~$\fopShu{X}_0$, and that~$X$ is \emph{coherently convergent} if it is terminating and coherently confluent. 
A left-monomial $1$-polygraph~$X$ is \emph{confluent} (resp.\ \emph{locally confluent}) if the $2$-polygraph $(X_0,X_1,\Sph(\fopShu{X}_1))$ has the corresponding coherent property, and it is \emph{convergent} if it is both terminating and confluent.

\subsubsection{Classification of local branchings}
\label{SSS:ClassificationLocalBranchings}
We distinguish the following four types of local branchings of a left-monomial $1$-polygraph~$X$:
\begin{enumerate}
\item \emph{aspherical branchings}: $(f,f)$, where $f$ is a rewriting step of $X$,
\item \emph{additive branchings}: $(\lambda f + \mu 1_v+1_c,\lambda 1_u+\mu g+1_c)$, where $f:u\fl a$ and $g:v\fl b$ are $1$-monomials of $\fopShu{X}_1$, $\lambda$ and $\mu$ are nonzero scalars, $c$ is a $0$-cell of $\fopShu{X}_0$, $u\neq v$, and $u,v\notin\supp(c)$.
\item \emph{multiplicative branchings}: $(\lambda\Gamma[f,1_v]+1_c,\lambda\Gamma[1_u,g]+1_c)$, where $\Gamma$ is a two-hole context of $\ftree{X}_0$, $f:u\fl a$ and $g:v\fl b$ are $1$-monomials of $\fopShu{X}_1$, $\lambda$ is a nonzero scalar, $c$ is a $0$-cell of $\fopShu{X}_0$, and $\Gamma[u,v]\notin\supp(c)$.
\item \emph{intersecting branchings}: local branchings that are neither aspherical, additive, nor multiplicative.
\end{enumerate}

We define a well-founded partial order $\sqsubseteq$ on branchings of $X$ as follows: for every one-hole context~$\Gamma$ of~$\ftree{X}_0$ and every $0$-cell $c$ of $\fopShu{X}_0$, we set
\[
(f,g)\sqsubseteq(\Gamma[f]+1_c,\Gamma[g]+1_c).
\]
The \emph{critical branchings} are the minimal intersecting branchings for this order. 
We denote the intersecting branchings by $(\Gamma[\alpha]+1_c,\Delta[\beta]+1_c)$, where $\alpha,\beta$ are $1$-generators of $X$, $\Gamma,\Delta$ are one-hole contexts of~$\ftree{X}_0$, and $s_0(\Gamma[\alpha])=s_0(\Delta[\beta])$. 
An \emph{essential branching} is a critical branching $(\Gamma[\alpha],\Delta[\beta])$ where $\Gamma[\alpha]$ and $\Delta[\beta]$ are consecutive $1$-monomials for the path-lexicographic monomial order $\prec_{pl}$ defined in \textsection~\ref{SSS:PathLexicographicMonomialOrder}.

Let~$X$ be a $2$-polygraph such that $X_{\leq 1}$ is left-monomial, and $u$ be a $0$-cell of~$\fopShu{X}_1$. We say that~$X$ is \emph{essentially coherently confluent at~$u$} if every essential branching of~$X$ of source~$u$ is coherently confluent, and that~$X$ is \emph{essentially coherently confluent} if it is so at every $0$-cell of~$\fopShu{X}_0$. 
A left-monomial $1$-polygraph $X$ is \emph{essentially confluent} if the $2$-polygraph $(X_0,X_1,\Sph(\fopShu{X}_1))$ is essentially coherently confluent.

As for polygraphs of associative algebras, we have:
\begin{lemma}[{{\cite[Lemmata 3.1.3 and 4.1.2.]{GuiraudHoffbeckMalbos19}}}]
\label{L:Papillon}
Let $X$ be a $2$-polygraph such that $X_{\leq 1}$ is left-monomial, and a $0$-cell $a$ in $\fopShu{X}_0$ such that $X$ is coherently confluent at $b$ for any $b\prec_{X_1}a$. If $f$ is a $1$-cell of $\fopShu{X}_1$ that decomposes
\[
a_0\xrightarrow{f_1}a_1\xrightarrow{f_2}\cdots\xrightarrow{f_p}a_p
\]
into $1$-cells of size $1$, with $a_i\prec_{X_1}a$ for all $i\in\range{p-1}$, then there exists a $0$-cell $a'$, $1$-cells $g,h$, and a $2$-cell $F$ in $\fopShu{X}_2$ such that 
\[
\xymatrix{
&
a_p
  \ar@/^/[dr] ^-{h}
  \ar@2 []!<0pt,-20pt>;[d]!<0pt,5pt> ^-*+{F}
\\
a_0
  \ar@/^/[ur] ^(0.5){f}
  \ar@/_/[rr] _-{g}
&
{}
&
a'
}
\]
When $p=1$, then $F$ is an identity $2$-cell.
\end{lemma}

\begin{theorem}[Coherent essential branchings theorem]
\label{T:CriticalBranchings}
Let~$X$ be a $2$-polygraph such that $X_{\leq 1}$ is terminating and left-monomial.
If~$X$ is essentially coherently confluent, then it is coherently confluent.
\end{theorem}
\begin{proof}
The structure of the proof is the same as for the similar result for associative algebras given in~\cite{GuiraudHoffbeckMalbos19}. The primary difference is that we prove that we can restrict the hypotheses to the critical branchings that are essential. 
Suppose that $X$ is an essentially coherently confluent $2$-polygraph. We proceed by well-founded induction on the sources of the branchings of $X_{\leq 1}$, with respect to the order $\prec_{X_1}$, to prove that~$X$ is coherently confluent at every $0$-cell of~$\fopShu{X}_0$. For each source, we first prove local coherent confluence and then deduce coherent confluence by Newman's lemma, exactly as for associative algebras in~\cite[Thm.\ 4.2.1]{GuiraudHoffbeckMalbos19}.

A reduced $0$-cell cannot be the source of a local branching, so~$X$ is coherently confluent at reduced $0$-cells.
Now, fix a nonreduced $0$-cell~$a_0$ of~$\fopShu{X}_0$, and assume that~$X$ is coherently confluent at every~$b\prec_{X_1} a_0$. Then we proceed by case analysis on the type of the local branchings.
If we show that critical branchings are coherently confluent, then the cases of aspherical, additive, multiplicative, and non-critical intersecting branchings are handled exactly as in~\cite[Thm.\ 4.2.1]{GuiraudHoffbeckMalbos19}. There remains to show the coherent confluence of critical branchings.

Let $(\Gamma[\alpha],\Delta[\beta])$ be a critical branching of source~$a_0$ as in \textsection~\ref{SSS:ClassificationLocalBranchings}. 
We proceed by induction on the size of the source of the branching and the number of $1$-monomials $\Lambda[\gamma]$ of $\fopShu{X}_1$ such that $\Gamma[\alpha]\prec_{pl} \Lambda[\gamma] \prec_{pl} \Delta[\beta]$.
If the branching is essential, then it is coherently confluent by hypothesis. Otherwise, there exists a $1$-monomial $\Lambda[\gamma]$ of $\fopShu{X}_1$ between $\Gamma[\alpha]$ and $\Delta[\beta]$. We get two branchings $(\Gamma[\alpha],\Lambda[\gamma])$ and $(\Lambda[\gamma],\Delta[\beta])$. 

The branching $(\Gamma[\alpha],\Lambda[\gamma])$ is either multiplicative or intersecting. If it is multiplicative, then it is coherently confluent by the multiplicative case. Otherwise, it is either non-minimal with respect to the order $\sqsubseteq$ or a critical branching. In the non-minimal case, there exists a factorisation 
\[
\xymatrix @R=0.8em @C=2.5em{
& \Gamma_1\Gamma_0[t(\alpha)]
	\ar@/^2ex/ [dr] ^-{\Gamma_1[f_0]}
		\ar@2 []!<0pt,-12.5pt>;[dd]!<0pt,12.5pt> ^-*+{F}
\\
\Gamma_1\Gamma_0[s(\alpha)]=\Gamma_1\Lambda_0[s(\gamma)]
	\ar@/^3ex/ [ur] ^-{\Gamma_1\Gamma_0[\alpha]}
	\ar@/_3ex/ [dr] _-{\Gamma_1\Lambda_0[\gamma]}
&& \Gamma_1[b_0]
\\
& \Gamma_1\Lambda_0[t(\gamma)]
	\ar@/_2ex/ [ur] _-{\Gamma_1[g_0]}
}
\]
where $(\Gamma_0[\alpha],\Lambda_0[\gamma])$ is a critical branching, and by induction on the size of its source, $\Gamma_1$ is a one-hole context of $\ftree{X}_0$, and $f_0,g_0$ are positive $1$-cells of $\fopShu{X}_1$ making this critical branching coherently confluent.
In the critical case, the branching is either essential, or not. If it is essential, then it is coherently confluent by hypothesis. Otherwise, it is coherently confluent by induction hypothesis.

We proceed similarly for the branching $(\Lambda[\gamma],\Delta[\beta])$: in every case, we write
\begin{itemize}
\item the branching $(\Gamma[\alpha],\Lambda[\gamma])$ as $(\Gamma_1\Gamma_0[\alpha],\Gamma_1\Lambda_0[\gamma])$, with $\Gamma_1$ a one-hole context of $\ftree{X}_0$ and $(\Gamma_0[\alpha],\Lambda_0[\gamma])$ a confluent branching,
\item the branching $(\Lambda[\gamma],\Delta[\beta])$ as $(\Delta_1\Lambda_0'[\gamma],\Delta_1\Delta_0[\beta])$, with $\Delta_1$ a one-hole context of $\ftree{X}_0$ and $(\Lambda_0'[\gamma],\Delta_0[\beta])$ a confluent branching. 
\end{itemize}
We then construct the following coherently confluent diagram
\[
\xymatrix @R=2em @C=2em {
& \Gamma_1\Gamma_0[t(\alpha)]
	\ar@/^/ [rr] ^-{\Gamma_1[f_0]}
	\ar@2 []!<+30pt,+12.5pt>;[d]!<+30pt,+8pt> ^-*+{F}
&
& \Gamma_1[b_0]
	\ar@/^/ [drr] ^-{h_1}
	\ar@2 []!<+5pt,-20pt>;[dd]!<+5pt,+30pt> ^-*+{H}
\\
\Gamma[s(\alpha)]
	\ar@/^2ex/ [ur] ^-*+{\Gamma[\alpha]}
	\ar@/_2ex/ [dr] _-*+{\Delta[\beta]}
	\ar[rr] |-*+{\Lambda[\gamma]}
&
         \ar@2 []!<+30pt,-12.5pt>;[d]!<+30pt,+8pt> ^-*+{G}
& \Lambda[t(\gamma)]
	\ar@/^1ex/ [ur] ^-{\Gamma_1[g_0]}
	\ar@/_1ex/ [dr] _-{\Delta_1[f_0']}
&&& c
\\
& \Delta_1\Delta_0[t(\beta)]
	\ar@/_/ [rr] _-{\Delta_1[g_0']}
&& \Delta_1[b_0']
	\ar@/_/ [urr] _-{h_1'}
}
\]
where the $2$-cells $F$ and $G$ are defined by the aforementioned coherently confluent branchings and the $2$-cell $H$ is given by induction hypothesis.
\end{proof}

\subsubsection{Operad presented by an ideal}
\label{SSS:OperadPresentedIdeal}

Let $X_0$ be an indexed set and $I$ an ideal of the free shuffle operad $\fopShu{X_0}$. We equip the collection $\Unit\oplus I$ with a shuffle operad structure, with unit $\eta:\Unit\hookrightarrow\Unit\oplus I$ and multiplication given by the following composition
\begin{align*}
\mu_I:(\Unit\oplus I)\circ_\Sha(\Unit\oplus I) & \simeq (\Unit\oplus I)\oplus (I\circ_\Sha(\Unit\oplus I)) \\
& \rightarrow(\Unit\oplus I)\oplus(I\circ_\Sha\fopShu{X_0}) \\
& \xrightarrow{1\oplus\rho} (\Unit\oplus I)\oplus I \\
& \rightarrow \Unit\oplus I,
\end{align*}
where $\rho$ is the right action of $I$ as an $\fopShu{X_0}$-bimodule. Denote by $\fopShu{X_0}/I$ the coequalizer of the morphisms of shuffle operads
\[
\xymatrix@C=3.5em{
\Unit\oplus I
  \ar@<0.7ex> [r] ^-{\eta \oplus 1}
  \ar@<-0.7ex> [r] _-{\eta \oplus 0}
&
\fopShu{X_0}
}
\]
in $\ShuO$. Note that the underlying collection of $\fopShu{X_0}/I$ is the cokernel of the inclusion $I\hookrightarrow\fopShu{X_0}$ in $\Coll$.

Let~$X$ be a $1$-polygraph. The \emph{boundary} of a $1$-generator $\alpha$ in $X$ is the $1$-cell $\partial(\alpha):=s_0(\alpha)-t_0(\alpha)$, and we set $\partial(X_1):=\{\,\partial(\alpha)\;|\;\alpha\in X_1\,\}$. 
We denote by~$I(X)$ the ideal of the free operad~$\fopShu{X}_0$ generated by the set of boundaries of the $1$-generators of~$X$, that is the free $\fopShu{X}_0$-bimodule generated by $\partial(X_1)$. Explicitly, the ideal~$I(X)$ is made of all the linear combinations
\[
\sum_{i=1}^p\lambda_i\Gamma_i[\dr(\alpha_i)]
\]
where $\lambda_i$ is a scalar and $\Gamma_i$ is a one-hole context. Note that the operad~$\cl{X}$ presented by~$X$ is isomorphic to~$\fopShu{X_0}/I(X)$.

We also have the result corresponding to~\cite[Prop.\ 3.3.4]{GuiraudHoffbeckMalbos19}:
\begin{proposition}
\label{P:ConfluenceAlgebraicCharacterisation}
For a terminating left-monomial $1$-polygraph $X$, the following assertions are equivalent :
\begin{enumerate}
\item $X$ is confluent.
\item $\Red(X)\cap I(X)=0$.
\item $\fopShu{X}_0 = \Red(X)\oplus I(X)$.
\end{enumerate} 
\end{proposition}

\subsection{Monomial-ordered shuffle polygraphs}
\label{SS:MonomialOrder}

In this subsection, we consider $1$-polygraphs whose orientations of $1$-generators are compatible with a fixed monomial order. We relate these polygraphs to the notion of Gröbner bases for operads introduced in~\cite{DotsenkoKhoroshkin10}. From Theorem~\ref{T:CriticalBranchings} we deduce a completion procedure for these polygraphs as in~\cite{DotsenkoKhoroshkin10}, but by resolving only essential branchings instead of all critical branchings.

\subsubsection{Gröbner bases~\cite{DotsenkoKhoroshkin10}}
Let~$X_0$ be an indexed set and~$\preccurlyeq$ be a monomial order on the free operad~$\fopShu{X}_0$. If~$a$ is a nonzero $0$-cell of~$\fopShu{X}_0$, the \emph{leading monomial of~$a$} is the maximum element~$\lm(a)$ of~$\supp(a)$ with respect to~$\preccurlyeq$, and $0$ when $\supp(a)$ is empty.
The \emph{leading coefficient of~$a$} is the coefficient~$\lc(a)$ of~$\lm(a)$ in~$a$, and the \emph{leading term of~$a$} is the element $\lt(a):=\lc(a)\lm(a)$ of~$\fopShu{X}_0$. 
Observe that, for~$a,b$ in~$\fopShu{X}_0$, we have~$a \prec b$ if, and only if, either $\lm(a) \prec \lm(b)$ or $(\lt(a)=\lt(b)\text{ and }a-\lt(a) \prec b-\lt(b))$. For $Y$ an indexed subset of $\fopShu{X}_0$, we denote by $\lm(Y)$ the indexed set of leading monomials of elements of~$Y$.

Let~$I$ be an ideal of the free operad~$\fopShu{X}_0$. A \emph{Gröbner basis for~$I$ with respect to~$\preccurlyeq$} is an indexed subset~$\Gr$ of~$I$ such that the ideals of~$\fopShu{X}_0$ generated by~$\lm(I)$ and by~$\lm(\Gr)$ coincide.

\begin{proposition}
\label{P:ConvergentGröbner}
If~$X$ is a convergent left-monomial $1$-polygraph, and~$\preccurlyeq$ is a monomial order on~$\fopShu{X}_0$ that is compatible with~$X_1$, then the indexed set~$\dr(X_1)$ forms a Gröbner basis of~$I(X)$. 

Conversely, let~$X_0$ be an indexed set, let~$\preccurlyeq$ be a monomial order on~$\fopShu{X}_0$, let~$I$ be an ideal of~$\fopShu{X}_0$ and~$\Gr$ be a subset of~$I$. Define $X(\Gr)$ as the $1$-polygraph with $0$-generators~$X_0$ and a $1$-generator
\[
\alpha_a\: : \: \lm(a) \fl \lm(a) - \frac{1}{\lc(a)} a
\]
for each~$a$ in~$\Gr$. If~$\Gr$ is a Gröbner basis for~$I$, then~$X(\Gr)$ is a convergent left-monomial presentation of~$\fopShu{X_0}/I$, such that $I(X(\Gr))=I$, and~$\preccurlyeq$ is compatible with~$X(\Gr)_1$.
\end{proposition}
\begin{proof}
Suppose that~$X$ is convergent. For every $1$-generator $\alpha$ of $X$, $\dr(\alpha)$ is in~$I(X)$. Since~$\preccurlyeq$ is compatible with~$X_1$, we have $\lm(\dr(\alpha))=s(\alpha)$ for every $1$-cell~$\alpha$ of~$X$. If $a$ is a nonzero $0$-cell of $I(X)$, then by Proposition~\ref{P:ConfluenceAlgebraicCharacterisation}, there exists a positive $1$-cell $a\fl 0$. By compatibility between $X_1$ and $\preccurlyeq$, any rewriting rule that does not reduce the leading monomial $\lm(a)$ will reduce $a$ into a $0$-cell with the same leading monomial. Thus, in order to rewrite to $0$, we must apply a rewriting rule to $\lm(a)$ at some point, and so $\lm(a)$ belongs to the ideal generated by the leading monomials of $\dr X_1$. 
Thus~$\dr(X_1)$ is a Gröbner basis for~$(I(X),\preccurlyeq)$.

Conversely, assume that~$\Gr$ is a Gröbner basis for~$(I,\preccurlyeq)$. 
The monomial order~$\preccurlyeq$ is compatible with~$X(\Gr)_1$, hence by Proposition~\ref{P:WellFoundedMonomialTerminating}, the polygraph $X(\Gr)$ terminates. Moreover, we have $I(X(\Gr))=I$, so the algebra presented by~$X(\Gr)$ is indeed isomorphic to $\fopShu{X_0}/I$. Moreover, the reduced monomials of~$\fopShu{X(\Gr)}$ are the monomials of~$\fopShu{X_0}$ that cannot be decomposed as~$\Gamma[\lm(a)]$ with~$a$ in~$\Gr$ and $\Gamma$ a one-hole context of $\ftree{X}_0$. Thus, if a reduced $0$-cell~$a$ of~$\fopShu{X}_0$ is in~$I$, its leading monomial must be~$0$, because~$\Gr$ is a Gröbner basis of~$(I,\preccurlyeq)$. By proposition~\ref{P:ConfluenceAlgebraicCharacterisation}, we get that the polygraph~$X(\Gr)$ is confluent. 
\end{proof}

\subsubsection{Poincaré-Birkhoff-Witt bases \cite{Hoffbeck10}}
Let $P$ be a operad, let $X_0$ be a generating indexed set of $P$, and let $\preccurlyeq$ be a monomial order of $\fopShu{X_0}$. A \emph{Poincaré-Birkhoff-Witt} (\emph{PBW}) \emph{basis for $(P,X_0,\preccurlyeq)$} is an indexed subset $\Br$ of $\ftree{X_0}$ such that:
\begin{enumerate}
\item $\Br$ is a linear basis of $P$, for $u\in \ftree{X_0}$, we write $[u]_\Br:=\sum_i\lambda_iw_i$ its decomposition in $P$ on the basis $\Br$,
\item for all $u,v$ in $\Br$ and all compatible elementary compositions $\circ_{i,\tau}$, either $u\circ_{i,\tau}v$ belongs to $\Br$ or $u\circ_{i,\tau}v \succ [u\circ_{i,\tau}v]_\Br$,
\item a tree monomial $u$ of $\ftree{X_0}$ is in $\Br$ if, and only if, for every decomposition $u=\Gamma(x\circ_{i,\tau}x')$ of $u$ where $x,x'\in X_0$ and $\Gamma$ is a one-hole context of $\fopShu{X}_0$, $x\circ_{i,\tau}x'\in\Br$.
\end{enumerate}

\begin{proposition}
\label{P:ConvergentPBW}
If $X$ is a convergent left-monomial quadratic presentation of an operad $P$, and $\preccurlyeq$ is a monomial order on $\fopShu{X_0}$ compatible with $X_1$, then the indexed set $\red(X)$ is a PBW basis for $(P,X_0,\preccurlyeq)$.

Conversely, let $P$ be a quadratic operad, $X$ a generating indexed set of $P$, $\preccurlyeq$ a monomial order on $\fopShu{X_0}$, and $\Br$ a PBW basis of $(A,X_0,\preccurlyeq)$. Define $X(\Br)$ as the $1$-polygraph with $0$-genrators $X_0$ and with a $1$-generator
\[
x\circ_{i,\tau}x' \xrightarrow{\alpha_{x\circ_{i,\tau}x'}} [x\circ_{i,\tau}x']_\Br
\]
for all $x,x'$ in $X_0\cap\Br$ such that $x\circ_{i,\tau}x'\neq[x\circ_{i,\tau}x']_\Br$ in $\fopShu{X_0}$. Then $X(\Br)$ is a quadratic convergent left-monomial presentation of $P$ such that $\red(X(\Br))=\Br$ and $\preccurlyeq$ is compatible with $X(\Br)_1$.
\end{proposition}
\begin{proof}
Suppose that $X$ is a quadratic convergent left-monomial presentation of  an operad $P$. By proposition \ref{P:ConfluenceAlgebraicCharacterisation}, we have the following exact sequence of collections:
\[
0 \fl I(X) \fl \fopShu{X_0} \fl \Red(X) \fl 0.
\]
Since $P$ is isomorphic to $\fopShu{X_0}/I(X)$ as an operad, it is also isomorphic to $\Red(X)$ as a collection, and therefore $\red(X)$ is a basis of $P$. The fact that $\preccurlyeq$ is compatible with $X_1$ implies axiom \Item{ii} of PBW bases. Axiom \Item{iii} comes from the definition of a reduced monomial for a quadratic left-monomial $1$-polygraph.

Conversely, assume that $\Br$ is a PBW basis for $(P,X,\preccurlyeq)$. By definition, $X(\Br)$ is quadratic and left-monomial, and axiom \Item{iii} of PBW bases implies $\red(X(\Br))\cap I(X(\Br))=0$. Termination of $X(\Br)$ is given by axiom \Item{ii} of PBW bases because $\preccurlyeq$ is well-founded. By proposition \ref{P:ConfluenceAlgebraicCharacterisation}, it is sufficient to prove that $\Red(X(\Br))\cap I(X(\Br))=0$ to get confluence: on the one hand, a reduced $0$-cell~$u$ of~$\Red(X(\Br))$ is a linear combination of $0$-cells of~$\Br$, so that~$u$ is its only normal form; and, on the other hand, if~$u$ belongs to~$I(X(\Br))$, then~$u$ admits~$0$ as a normal form. Finally, the operad presented by~$X(\Br)$ is isomorphic to~$\Red(X(\Br))$, that is to~$\K\Br$, hence to~$P$ by the previous exact sequence and because~$\Br$ is a linear basis of~$P$.
\end{proof}

\subsubsection{Completion procedure}
\label{SSS:CompletionProcedure}
For a $2$-polygraph $X$ where $X_2=\Sph(\fopShu{X}_1)$, Theorem~\ref{T:CriticalBranchings} leads to a completion procedure for $1$-polygraphs that reaches a convergent polygraph by resolving essential branchings. Given a terminating $1$-polygraph $X$, and a monomial order $\prec$ on~$\ftree{X}_0$ compatible with $X_1$, the procedure works as follows 
\begin{enumerate}[{\bf (i)}]
\item for every essential branching $(f,g)$ of $X$, the $0$-cells $t_0(f)$ and $t_0(g)$ are reduced to some normal forms $\rep{t_0(f)}$ and $\rep{t_0(g)}$. If $\rep{t_0(f)}\neq \rep{t_0(g)}$:
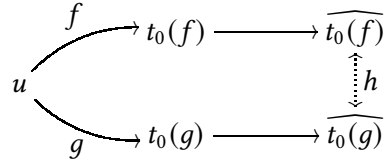
\[
\xymatrix @C=3.5em @R=0.4em{
& t_0(f)
\ar[r]
& 
\rep{t_0(f)}
\ar@{<..>}[dd] ^-{h}
\\
u
  \ar@/^2ex/ [ur] ^-{f}
  \ar@/_2ex/ [dr] _-{g}
  &&
\\
& t_0(g)
\ar[r] 
&
\rep{t_0(g)}
}
\]
a $1$-generator $h : \lm(a) \fl a-\lm(a)$, where $a=\rep{t_0(f)} - \rep{t_0(g)}$, is added to reach confluence of the branching ;
\item the addition of $1$-generators in the step {\bf (i)} can create new essential branchings, whose confluence must also be completed as in {\bf (i)} ;
\item Repeat the previous steps until there are no non-confluent essential branchings.
\end{enumerate}

As a consequence of Theorem~\ref{T:CriticalBranchings}, we have

\begin{proposition}
The procedure~\textsection~\ref{SSS:CompletionProcedure} on a $1$-polygraph $X$ produces a (possibly infinite) convergent polygraph that presents the operad $\cl{X}$.
\end{proposition}

An analogue completion procedure for non-symmetric operads has been described in detail with an explicit handling of the critical branchings in~\cite[Algorithm~2]{MalbosRen21}.

\section{Shuffle polygraphic resolutions from convergence}
\label{S:ShufflePolygraphicResolutionConvergence}

In this section, unless otherwise specified, all operads and polygraphs are shuffle. We recall from \cite{GuiraudHoffbeckMalbos19} the characterization of the property of acyclicity for an $\omega$-polygraph through the existence of a homotopical contraction. Subsection~\ref{SS:OverlappingPolygraphicResolution}, presents the main result of this article, Theorem~\ref{T:MainTheoremD}, that extends a reduced convergent left-monomial $1$-polygraph into a polygraphic resolution of the presented operad. 
In Subsection~\ref{SS:BimodulesResolution}, given a polygraphic resolution of an operad, we construct a bimodule resolution for the operad.
Finally, in Subsection~\ref{SS:ConfluenceKoszul} we prove a criterion of Koszulness in terms of quadratic convergence. 

\subsection{Polygraphic resolutions and contractions}
\label{SS:PolygraphicResolutionsContracations}

In this first subsection, we extend to $\omega$-operads the notion of homotopy developed in \cite{GuiraudHoffbeckMalbos19} for $\omega$-algebras, see also \cite{AraMaltsiniotis15} and~\cite{Guiraud19}. Then we introduce the notion of a \emph{contraction} of a polygraph, which allow us to characterize acyclic $\omega$-polygraphs.

\subsubsection{Homotopies}
Let $P$ and $Q$ be $\omega$-operads and $F,G:P\fl Q$ be morphisms of $\omega$-operads. A \emph{homotopy} from $F$ to $G$ is a graded linear map
\[
\eta : P \to Q
\]
of degree $1$, \emph{i.e.}, $\eta$ sends $n$-cells to $n+1$-cells), such that, writing $\eta_a$ for $\eta(a)$,
\begin{enumerate}
\item for every $n\geq 0$, for every $n$-cell $a$ of $P$,
\begin{eqn}{align}
\label{E:SourceHomotopy}
s_n(\eta_a) & = F(a)\star_0\eta_{t_0(a)}\star_1\cdots\star_{n-1}\eta_{t_{n-1}(a)} \\
\label{E:TargetHomotopy}
t_n(\eta_a) & = \eta_{s_{n-1}(a)}\star_{n-1}\cdots\star_1\eta_{s_0(a)}\star_0G(a),
\end{eqn}
\item for all $0\leq k<n$ and every $\star_k$-composable pair $(a,b)$ of $n$-cells of $P$,
\begin{align*}
\eta_{a\star_kb} ={} & F(s_{k+1}(a))\star_0\eta_{t_0(b)}\star_1\cdots\star_{k-1}\eta_{t_{k-1}(b)}\star_k\eta_b \\
& \star_{k+1}\eta_a\star_k\eta_{s_{k-1}(a)}\star_{k-1}\cdots\star_1\eta_{s_0(a)}\star_0G(t_{k-1}(b)),
\end{align*}
\item for all $n\geq 0$ and every $n$-cell $a$ of $P$,
\[
\eta_{1_a}=1_{\eta_a}.
\]
\end{enumerate}
In order for this definition to be licit, we need to check that the $\star_k$-compositions of \Item{i} are well defined. See {\cite[Appendix B.8]{AraMaltsiniotis15}} or {\cite[\textsection~5.1.1]{GuiraudHoffbeckMalbos19}} for the verification.
Note that the mappings~$a\mapsto s(\eta_a)$ and~$a\mapsto t(\eta_a)$ are operads morphisms because both are composites of operad morphisms. The globularity of $\eta_a$ follows from 
\begin{align*}
& ss(\eta_a)
	= s(F(a)) \star_0 \eta_{t_0(a)} \star_1 \cdots \star_{n-2} \eta_{t_{n-2}(a)} 
	= s(\eta_{s(a)}) 
	= st(\eta_a)
\\[0.5ex]
\text{and}\quad 
&
ts(\eta_a) 
	= t(\eta_{t(a)})
    = \eta_{s_{n-2}(a)} \star_{n-2} \cdots \star_1 \eta_{s_0(a)} \star_0 t(G(a))
	= tt(\eta_a).
\end{align*}

\subsubsection{}
Let us expand the homotopy~$\eta$ in low dimension. It maps a $1$-cell~$f:a\fl a'$ of~$P$ to a $2$-cell
\[
\xymatrix @R=0.25em @C=1.5em {
& F(a')
	\ar@/^/ [dr] ^-{\eta_{a'}}
	\ar@2 []!<0pt,-15pt>;[dd]!<0pt,15pt> ^-*+{\eta_f}
\\
F(a)
	\ar@/^/ [ur] ^-{F(f)}
	\ar@/_/ [dr] _-{\eta_a}
&& G(a')
\\
& G(a)
	\ar@/_/ [ur] _-{G(f)}
}
\]
of~$Q$, and a $2$-cell $A:f\dfl f':a\fl a'$ of~$P$ to the following $3$-cell of $Q$
\[
\vcenter{\xymatrix @R=0.25em @C=4em {
& F(a')
	\ar@/^/ [dr] ^-{\eta_{a'}}
		\ar@2 []!<5pt,-15pt>;[dd]!<5pt,15pt> ^-*+{\eta_{f'}}
\\
F(a)
	\ar@/^5ex/ [ur] ^(0.2){F(f)} _{}="s"
	\ar@/_/ [ur] |-*+{F(f')} ^{}="t"
		\ar@2 "s"!<-6pt,-10pt>;"t"!<-9pt,10pt> ^-{F(A)}
	\ar@/_/ [dr] _{\eta_a}
&& {G(a')} 
\\
& G(a)
	\ar@/_/ [ur] _-{G(f')}
}}
\quad\overset{\eta_A}{\Rrightarrow}\quad
\vcenter{\xymatrix @R=0.25em @C=4em  {
& F(a')
	\ar@/^/ [dr] ^-{\eta_{a'}}
		\ar@2 []!<-10pt,-15pt>;[dd]!<-10pt,15pt> ^-*+{\eta_{f}}
\\
F(a)
	\ar@/^/ [ur] ^-{F(f)}
	\ar@/_/ [dr] _{\eta_a}
&& {G(a')} 
\\
& G(a)
	\ar@/^/ [ur] |-*+{G(f)} _-{}="s"
	\ar@/_5ex/ [ur] _(0.8){G(f')} ^{}="t"
		\ar@2 "s"!<-7pt,-12pt>;"t"!<-10pt,7pt> ^-{G(A)}
}}
\]

\subsubsection{Unital sections and contractions}
Let $X$ be a $\omega$-polygraph. A \emph{unital section} of $X$ is a morphism of $\omega$-operads $\iota:\cl X\fl \fopShu{X}$, which is a section of the canonical projection $\pi:\fopShu{X}\twoheadrightarrow\cl X$, and such that $\iota_{1}=1$, where $1\in\K\subseteq \fopShu{X}(1)$.
The morphism $\iota$ assigns to every $0$-cell $a$ of $\overline X$ a representative $0$-cell $\iota_a$ in $\fopShu{X}$, in such a way that is the identity on the unit $\K$. Note that a unital section is not necessarily compatible with shuffle composition. For $a$ an $n$-cell of $\fopShu{X}$, we will write $\widehat a$ for $\iota\pi(a)$. Note that $\rep a=1_{\rep{s_0(a)}}$ for $n\geq 1$.

Fix $\iota$ a unital section of $X$. An \emph{$\iota$-contraction} of $X$ is a homotopy $\sigma:\id_{\fopShu{X}}\fl\iota\pi$ such that $\sigma_a=1_a$ for every $n$-cell $a$ of $\fopShu{X}$ that belongs to the image of $\iota$ or~$\sigma$.
We say that $\sigma$ is a \emph{right} $\iota$-contraction if, for all $n\geq 0$, $n$-cells $f,g$ of $\fopShu{X}$, and compatible elementary composition $\circ_{i,\tau}$,
\begin{eqn}{equation}
\label{E:RightContraction}
\sigma_{f\circ_{i,\tau}g}=(s_0(f)\circ_{i,\tau}\sigma_g)\star_0\sigma_{f\circ_{i,\tau}\widehat g}.
\end{eqn}

\begin{lemma}
\label{L:AbelianContraction}
Let $\sigma$ be a $\iota$-contraction. For $n\geq 1$ and every $n$-cell $a$ of $\fopShu{X}$,
\begin{eqn}{equation}
\label{E:ContractionForm}
s_n(\sigma_a)=a-t_{n-1}(a)+\sigma_{t_{n-1}(a)}\quad\text{and}\quad t_n(\sigma_a)=\sigma_{s_{n-1}(a)}.
\end{eqn}
\end{lemma}
Note that for $a$ a $0$-cell of $\fopShu{X}$, $s_0(\sigma_a)=a$ and $t_0(\sigma_a)=\rep a$.
\begin{proof}
Let us first prove
\[
a\star_0\sigma_{t_0(a)}\star_1\cdots\star_k\sigma_{t_k(a)}=a-t_k(a)+\sigma_{t_k(a)}
\]
by induction on $k\in\{0,\ldots,n-1\}$. The result is clear for $k=0$. For $k\geq 1$, we calculate
\begin{align*}
a\star_0\eta_{t_0(a)}\star_1\cdots\star_k\sigma_{t_k(a)} ={} & a\star_0\sigma_{t_0(a)}\star_1\cdots\star_{k-1}\sigma_{t_{k-1}(a)} \\
& - t_k(a\star_0\sigma_{t_0(a)}\star_1\cdots\star_{k-1}\sigma_{t_{k-1}(a)}) \\
& + \sigma_{t_k(a)} \\
={} & (a-t_{k-1}(a)+\sigma_{t_{k-1}(a)})-t_k(a-t_{k-1}(a)+\sigma_{t_{k-1}(a)})+ \sigma_{t_k(a)} \\
={} & a-t_k(a)+\sigma_{t_k(a)},
\end{align*}
the last equality coming from the fact that $t_kt_{k-1}(a)=t_{k-1}(a)$ and $t_k(\sigma_{t_{k-1}(a)})=\sigma_{t_{k-1}(a)}$. Applying $k=n-1$ and \eqref{E:SourceHomotopy} with $F=\id_{\fopShu{X}}$, we conclude that
\[
s_n(\sigma_a)=a-t_{n-1}(a)+\sigma_{t_{n-1}(a)}.
\]
For the second equation, we proceed similarly to show that, for all $k\in\{0,\ldots,n-1\}$,
\[
\sigma_{s_k(a)}\star_k\cdots\star_1\sigma_{s_0(a)}\star_0\rep a =\rep a-s_k(\rep a)+\sigma_{s_k(a)}=\sigma_{s_k(a)}
\]
because $\rep a=1_{\rep{s_0(a)}}$. Applying $k=n-1$ and \eqref{E:TargetHomotopy} with $G=\iota\pi$, we conclude that
\[
t_n(\sigma_a)=\sigma_{s_0(a)}.
\]
\end{proof}

\subsubsection{Reduced and essential monomials}
Let $\iota$ be an unital section of $X$, and $\sigma$ an $\iota$-contraction of an $\omega$-polygraph $X$. A $0$-monomial $u$ of $\fopShu{X}$ is \emph{$\iota$-reduced} if $\widehat u=u$. A non-$\iota$-reduced $0$-monomial $u$ of the free $\omega$-operad $\fopShu{X}$ is \emph{$\iota$-essential} if $u=(x\mid\vec v)$ where $x$ is a $0$-generator of $X$ and $v_1,\ldots,v_k$ are $\iota$-reduced $0$-monomials of $\fopShu{X}$.
When the underlying $1$-polygraph $X_{\leq 1}$ is convergent, and the section $\iota$ sends a $0$-monomial on its unique normal form with respect to $X_1$, the $\iota$-reduced $0$-monomials coincide with reduced ones.

For $n\geq 0$, an $n$-monomial $a$ of $\fopShu{X}$ is \emph{$\sigma$-reduced} if it is an identity or in the image of $\sigma$. If $\sigma$ is a right $\iota$-contraction of $X$ and $n\geq 0$, then a non-$\sigma$-reduced $n$-monomial $a$ of $\fopShu{X}$ is \emph{$\sigma$-essential} if $a=(\alpha\mid\vec v)$, where $\alpha$ is a $n$-generator of $X$ and $v_1,\ldots,v_k$ are $\iota$-reduced $0$-monomials of the $\omega$-operad $\fopShu{X}$.

\begin{lemma}
\label{L:EssentialValues}
Let $X$ be an $\omega$-polygraph and $\iota$ a unital section of $X$. A right $\iota$-contraction $\sigma$ of $X$ is uniquely and entirely determined by its values on the $\iota$-essential $0$-monomials and, for $n\geq 1$, on the $\sigma$-essential $n$-monomials of $\fopShu{X}$.
\end{lemma}
\begin{proof}
The proof follows the same arguments as in the case of associative algebras given in {\cite[Section 5.2]{GuiraudHoffbeckMalbos19}}, and it is divided in two steps:
\begin{enumerate}
\item First, we prove that a homotopy $\eta:F\fl G$ between morphisms of $\omega$-operads $F,G:\fopShu{X}\fl\fopShu{X}$ is uniquely and entirely determined by its values on $n$-monomials for all $n\geq 0$, provided it satisfies the following relation:
\begin{eqn}{equation}
\label{E:HomotopyExchange}
\eta_{\mu_{\fopShu{X_n}}^\uparrow} = \eta_{\mu_{\fopShu{X_n}}^\downarrow},
\end{eqn}
where $\mu_{\fopShu{X_n}}^\uparrow$ and $\mu_{\fopShu{X_n}}^\downarrow$ are defined considering $\fopShu{X_n}$ as a $\fopShu{X_0}$-bimodule.
\item Next, we prove that the values of a right $\iota$-contraction on $n$-monomials are uniquely and entirely determined by the values on $\iota$-essential and $\sigma$-essential monomials, and that the resulting values satisfy \eqref{E:HomotopyExchange}.
\end{enumerate}
\Item{i} Proceed by induction on $n\geq 0$. For~$n=0$, assume that $\eta_u:F(u)\fl G(u)$ is a fixed $1$-cell of~$\fopShu{X}$ for every $0$-monomial~$u$ of~$\fopShu{X}$. Extend~$\eta$ uniquely to every $0$-cell~$a$ of~$\fopShu{X}$ by linearity.

Now fix $n\geq 1$ and assume that an $(n+1)$-cell $\eta_u$ of $\fopShu{X}$ has been chosen for every $n$-monomial $u$ of $\fopShu{X}$, with source and target given by the definition of homotopies, such that \eqref{E:HomotopyExchange} holds for $n$-monomials. By construction, the $n$-cells of $\fopShu{X}$ are linear combinations of $n$-monomials of $\fopShu{X}$ and of identities of $(n-1)$-cells of $\fopShu{X}$ up to the relation
\[
\mu_{\fopShu{X_0}\langle X_n\rangle\oplus\fopShu{X_{n-1}}}^\uparrow = \mu_{\fopShu{X_0}\langle X_n\rangle\oplus\fopShu{X_{n-1}}}^\downarrow.
\]
Thus we can extend $\eta$ to all $n$-cells $a$ of $\fopShu{X}$ by choosing a decomposition of $a$ into a linear combination of $n$-monomials and an identity, and using \eqref{E:HomotopyExchange} to ensure that the resulting cell does not depend on the choice of decomposition. We check that the source and target of the resulting $(n+1)$-cell $\eta_a$ match the definition of homotopies by linearity of $F,G$ and the $\star_k$-compositions.

\Item{ii} First, we construct $\sigma$ as a graded linear map by induction on $n$. For $n=0$, if $u$ is a non-$\iota$-essential monomial, then either $u=\rep{u}$, or $u=(x\mid\vec v)$ where $x$ is a $0$-cell of~$X$ and some $v_i$ is a non-$\iota$-reduced monomial. In the former case, $\sigma_u=1_u$ is forced because~$u$ is $\iota$-reduced. In the latter case, take $i$ maximal. Writing $(x\mid\vec v)=(x\mid v_1\cdots 1\cdots v_k)\circ_{i,\tau}v_i$ for some shuffle permutation $\tau$, \eqref{E:RightContraction} imposes
\[
\sigma_{(x\mid\vec v)} = ((x\mid v_1\cdots 1\cdots v_k)\circ_{i,\tau}\sigma_{v_i}) \star_0 \sigma_{(x\mid v_1\cdots\rep{v_i}\cdots v_k)}.
\]
Then proceed by induction on the weight of the~$v_i$ to define~$\sigma_{v_i}$ from the values of~$\sigma$ on $\iota$-reduced monomials.

Now let $n\geq 1$. For every $n$-monomial $\Gamma[\alpha]$, with $\alpha$ a $n$-generator of $X$ and $\Gamma$ a one-hole context of~$\ftree{X}_0$, writing
\[
\Gamma[\alpha]=u\circ_{i,\tau}(\alpha\mid\vec v)
\]
and
\[
(\alpha\mid\vec v)=\alpha\circ_{k,\tau_k}v_k\circ_{k-1,\tau_{k-1}}\cdots\circ_{1,\tau_1}v_1,
\]
the equation \eqref{E:RightContraction} imposes that we set
\begin{align*}
\sigma_{\Gamma[\alpha]}:= {}& (u\circ_{i,\tau}(s_0(\alpha)\mid\sigma_{v_1}\,v_2\cdots v_k))\star_0\cdots\star_0(u\circ_{i,\tau}(s_0(\alpha)\mid\rep{v_1}\cdots\rep{v_{k-1}}\,\sigma_{v_k})) \\
& \star_0(u\circ_{i,\tau}\sigma_{(\alpha\mid\vec{\rep{v}})})\star_0\sigma_{u\circ_{i,\tau}\rep{(\alpha\mid\vec v)}},
\end{align*}
where $(\alpha\mid\vec{\rep{v}})$ is a shortcut for $(\alpha\mid\rep{v_1}\cdots\rep{v_k})$.  Let us check that this definition is well-founded. The $\sigma_{v_i}$ are defined by induction on the weight of the $v_i$, and $\sigma_{u\circ_{i,\tau}\rep{(\alpha\mid\vec v)}}=\sigma_{u\circ_{i,\tau}\rep{(s_0(\alpha)\mid\vec v)}}$ is defined by induction on $n$. It remains to check that $\sigma_{(\alpha\mid\vec{\rep{v}})}$ is defined. If $(\alpha\mid\vec{\rep{v}})$ is $\sigma$-essential, then it is defined by hypothesis. Otherwise, $(\alpha\mid\vec{\rep{v}})$ is $\sigma$-reduced, in which case $(\alpha\mid\vec{\rep{v}})=\sigma_b$ for some $(n-1)$-cell $b$ of $\fopShu{X}$, which imposes $\sigma_{(\alpha\mid\vec{\rep{v}})}:=1_{\sigma_b}$.

Now it remains only to show \eqref{E:HomotopyExchange} and then apply the first point. More explicitly, we need to show
\[
\sigma_{u\circ_{i,\tau}s_0(v)}+\sigma_{t_0(u)\circ_{i,\tau}v}-\sigma_{t_0(u)\circ_{i,\tau}s_0(v)} = \sigma_{u\circ_{i,\tau}t_0(v)}+\sigma_{s_0(u)\circ_{i,\tau}v}-\sigma_{s_0(u)\circ_{i,\tau}t_0(v)}
\]
for all $u,v$ two $n$-monomials of $\fopShu{X}$ and compatible elementary composition $\circ_{i,\tau}$. Write $a=s_0(u)$, $a'=t_0(u)$, $b=s_0(v)$, $b'=t_0(v)$, and $\cdot=\circ_{i,\tau}$. On the one hand,
\begin{align*}
\sigma_{u\cdot b}+\sigma_{a'\cdot v}-\sigma_{a'\cdot b} & = (a\cdot\sigma_b)\star_0\sigma_{u\cdot\rep{b}}+(a'\cdot\sigma_v)\star_0\sigma_{a'\cdot\rep{v}}-(a'\cdot\sigma_b)\star_0\sigma_{a'\cdot\rep{b}} \\
& = a\cdot\sigma_b+\sigma_{u\cdot\rep{b}}-a\cdot\rep{b}+a'\cdot\sigma_v-a'\cdot\sigma_b,
\end{align*}
and on the other hand,
\begin{align*}
\sigma_{u\cdot b'}+\sigma_{a\cdot v}-\sigma{a\cdot b'} & = (a\cdot\sigma_{b'})\star_0\sigma_{u\cdot\rep{b'}} + (a\cdot\sigma_{v})\star_0\sigma_{a\cdot\rep{v}} - (a\cdot\sigma_{b'})\star_0\sigma_{a\cdot\rep{b'}} \\
& = \sigma_{u\cdot\rep{b'}}+a\cdot\sigma_v-a\cdot\rep{v} \\
& = \sigma_{u\cdot\rep{b}} + a\cdot\sigma_v - a\cdot\rep{b}
\end{align*}
Therefore it remains to prove
\begin{eqn}{equation}
\label{E:ExchangeC}
a\cdot\sigma_b + a'\cdot\sigma_v = a\cdot\sigma_v + a'\cdot\sigma_b.
\end{eqn}
Since $\star_0$-composition in $\fopShu{X}$ is a morphism of $\omega$-operads, we have
\[
u\cdot\sigma_b = u\cdot b\star_0a'\cdot\sigma_b = a\cdot\sigma_b\star_0u\cdot\rep{b}.
\]
Using the linear expression of $\star_0$-composition, we get
\[
u\cdot b+a'\cdot\sigma_b - a'\cdot b = a\cdot\sigma_b + u\cdot\rep{b} - a\cdot\rep{b}
\]
Similarly, considering $u\cdot\sigma_v$, we get
\[
u\cdot b+a'\cdot\sigma_v-a'\cdot b = a\cdot\sigma_v+u\cdot\rep{b}-a\cdot\rep{b}.
\]
Taking the difference of the two previous equations gives us \eqref{E:ExchangeC}.
\end{proof}

\begin{proposition}
\label{P:AcyclicPolygraph}
Let $X$ be an $\omega$-polygraph with a fixed unital section $\iota$. Then $X$ is a polygraphic resolution of the $\omega$-operad $\cl{X}$ if, and only if, $X$ admits a right $\iota$-contraction.
\end{proposition}
\begin{proof}
Suppose that $X$ is a polygraphic resolution of the operad $\cl X$, and define a right $\iota$-contraction $\sigma$ of $X$. Using Lemma~\ref{L:EssentialValues}, we shall define $\sigma$ on $\iota$- and $\sigma$-essential $n$-monomials of $\fopShu{X}$ by induction on $n\geq 0$. If $(x\mid\vec v)$ is an $\iota$-essential $0$-monomial, then $\pi_X(x\mid\vec v)=\pi_X(\rep{(x\mid\vec v)})$ in $\cl X$, hence there exists a $1$-cell $\sigma_{(x\mid\vec v)}:(x\mid\vec v)\fl\rep{(x\mid\vec v)}$ in $\fopShu{X}$. Now assume that $\sigma$ is defined on the $n$-cells of $\fopShu{X}$ for $n\geq 0$ and let $(\alpha\mid\vec v)$ be a $\sigma$-essential $(n+1)$-monomial of $\fopShu{X}$. The $n$-cells defining $s(\sigma_{(\alpha\mid\vec v)})$ and $t(\sigma_{(\alpha\mid\vec v)})$ as in~\eqref{E:ContractionForm} are parallel, so, by acyclicity of $X$, there exists an $(n+2)$-cell $\sigma_{(\alpha\mid\vec v)}$ with this source and target in~$\fopShu{X}$.

Conversely, let $\sigma$ be a right $\iota$-contraction of the polygraph $X$, and let $a$, $b$ be parallel $n$-cells of $\fopShu{X}$ for $n\geq 1$. We have $t(\sigma_a)=\sigma_{s(a)}=\sigma_{s(b)}=t(\sigma_b)$ by \eqref{E:ContractionForm}, so the $(n+1)$-cell $\sigma_a\star_n\sigma_b^-$ is well defined, with source $s(\sigma_a)$ and target $s(\sigma_b)$. Since $t_k(a)=t_k(b)$ for $k\in\{0,\ldots,n-1\}$, we find that
\[
(\sigma_a\star_n\sigma_b^-)\star_{n-1}\sigma_{t_{n-1}(a)}^-\star_{n-2}\cdots\star_1\sigma_{t_0(a)}^-
\]
is a well defined $(n+1)$-cell of $\fopShu{X}$ of source $a$ and target $b$, thus proving that $X_{n+1}$ is an acyclic extension of $\fopShu{X}_n$. Thus $X$ is a polygraphic resolution of $\cl X$.
\end{proof}

\subsection{Polygraphic resolution from a convergent presentation}
\label{SS:OverlappingPolygraphicResolution}

This subsection contains the main result of this article. We show how to extend a reduced left-monomial convergent shuffle $1$-polygraph into a shuffle polygraphic resolution of its presented operad. The $n$-generators of the resolution correspond to certain overlappings of the $1$-generators of the polygraph. 

\subsubsection{Higher-dimensional overlappings}
\label{SSS:HigherDimensionalOverlappings}
Let $X$ be a left-monomial $1$-polygraph, and consider the path-lexicographic order $\prec_{pl}$ on $1$-monomials of $X$ defined in \textsection~\ref{SSS:PathLexicographicMonomialOrder}. We define the family of indexed sets $\Ov(X) = (\Ov(X)_n)_{n\geq 0}$ by induction on $n \geq 0$. The elements of $\Ov(X)_n$ are called \emph{$n$-overlappings of $X$}, and for an $n$-overlapping $u_n$ we will also define its \emph{source} $s_0(u_n)$ and its set of \emph{branches} $B(u_n)$.

For $n=0$, define a $0$-overlapping $u_0$ as a $0$-generator in $X_0$. Define its source as $s_0(u_0):=u_0$ and its set of \emph{branches} as $B(u_0):=\varnothing$.

Now suppose that $n$-overlappings are defined for $n\geq 0$. Let $u_n$ be an $n$-overlapping and $B(u_n) = \{\Gamma_1[\alpha_1] \prec_{pl} \cdots \prec_{pl} \Gamma_n[\alpha_n]\}$ its set of branches, where each $\Gamma_k$ is a one-hole context and $\alpha_k$ is a $1$-generator in $X_1$. Given $0$-monomials $\vec v_{n + 1}$, we define
\[
E(u_n, \vec v_{n + 1}) := \left\{\;
\Gamma[\alpha]
\; \left| \;
\begin{array}{l}
\Gamma \text{ one-hole context}, \alpha \in X_1, \\
\Gamma[s_0(\alpha)] = (s_0(u_n) \mid \vec v_{n + 1}), \\
\Gamma[\alpha] \succ_{pl} (\Gamma_n[\alpha_n] \mid \vec v_{n + 1})
\end{array}
\right.\;\right\}.
\]
An \emph{$(n + 1)$-overlapping} is a tuple $(u_n, \vec v_{n + 1})$, denoted by $u_n \tri \vec v_{n + 1}$, where $u_n$ is an $n$-overlapping, and~$\vec v_{n + 1}$ is a list of reduced $0$-monomials such that, for any list of rooted submonomials $\vec w_{n + 1} \subsetneq \vec v_{n + 1}$, $\#E(u_n, \vec w_{n + 1}) < \#E(u_n, \vec v_{n + 1})$. We then define its source as $s_0(u_n \tri \vec v_{n + 1}) := (s_0(u_n) \mid \vec v_{n + 1})$ and its set of branches as $B(u_n \tri \vec v_{n + 1}) := \{(\Gamma_k[\alpha_k] \mid \vec v_{n + 1}) \mid 1\leq k\leq n\} \cup \{\max E(u_n,\vec v_{n+1})\}$.

\subsubsection{Crowns}
\label{SSS:Crowns}
An $(n + 1)$-overlapping $u_{n} \tri\vec v_{n+1}$ can be represented graphically as
\[
\begin{tikzpicture}[scale=.6, baseline={([yshift=-.5ex]current bounding box.center)}]
\draw
  (0,0) -- (-1.5,1) -- (1.5,1) -- cycle
  (-1.5,1) -- (-3,2) -- (3,2) -- (1.5,1);
\node at (0,.6) {$u_n$};
\node at (-1.3,1.6) {$v_{n+1,1}$};
\node at (0,1.6) {$\cdots$};
\node at (1.3,1.6) {$v_{n+1,k}$};
\end{tikzpicture}
=
\begin{tikzpicture}[scale=.6, baseline={([yshift=-.5ex]current bounding box.center)}]
\draw
  (0, 0) -- (-3, 2) -- (3, 2) -- cycle
  (-1.2, 2) -- (.9, .6);
\node at (-.8, 1.1) {$\Gamma$};
\node at (1, 1.6) {$s_0(\alpha)$};
\end{tikzpicture}
\]
where $\Gamma[\alpha]$ is the maximal element of $E(u_n, \vec v_{n + 1})$. We call the list $\vec v_{n+1}$ of reduced $0$-monomial a \emph{crown on $u_{n}$}.
Given an $n$-overlapping $u_n$ and a list of $0$-monomials $\vec v_{n+1}$, we define
\[
C(u_n,\vec v_{n+1}):=\{\vec w_{n+1}\subseteq\vec v_{n+1}\mid u_n\tri\vec w_{n+1}\in\Ov(X)_{n+1}\},
\]
that is, the set of crowns $\vec w_{n+1}$ on $u_n$ included in $\vec v_{n+1}$. This set is equipped with the total order defined by $\vec w_{n+1} \prec \vec w'_{n+1}$ if $\max_{\prec_{pl}} B(u_n\tri\vec w_{n+1}) \prec_{pl} \max_{\prec_{pl}} B(u_n\tri\vec w'_{n+1})$. Note that $C(u_n,\vec v_{n+1})$ is empty if, and only if, $E(u_n,\vec v_{n+1})$ is empty.

\subsubsection{Description in low dimensions}
\label{SSS:ReducedLowDimension}
Let us look at the definitions of $n$-overlappings in low dimensions.
A $1$-overlapping $u_0\tri\vec v_1$ is associated to a single branch $\Gamma[\alpha]$. Since the $0$-monomials $\vec v_1$ are reduced and minimal, the context $C$ must be trivial. Thus $\Ov(X)_1$ is in bijection with $X_1$, and this bijection is given by taking the unique branch of the $1$-overlapping.
Next, a $2$-overlapping $u_0\tri\vec v_1\tri\vec v_2$ corresponds to a pair of branches $(\Gamma_1[\alpha_1], \Gamma_2[\alpha_2])$ which form an critical branching in context. Since the crown $\vec v_2$ must be minimal, this context is trivial. Thus $\Ov(X)_2$ is in bijection with the set of critical branchings.

\subsubsection{Overlappings as paths of crowns}
\label{SSS:OverlappingsAsPaths}
Given a left-monomial $1$-polygraph $X$, the $n$-overlappings can be defined inductively as certain paths of length $n$ in the directed graph $\Gr(X)$ defined as follows. Its vertices are the $0$-monomials of $\ftree{X}_0$, and its edges are 
\[
\xymatrix @C=2.5em{
u
  \ar[r] ^-{\tri \vec v}
&
(u | \vec v)
},
\]
such that $u$ is the source of an overlapping and $\vec v$ is a crown. 
Then the indexed set of $n$-overlappings of~$X$ corresponds to a subset of paths of $\Gr(X)$, starting in $X_0$ and of length $n$, where each step of the path corresponds to the addition of a crown. That is $u_0\tri\vec v_1\tri\cdots\tri\vec v_n$ corresponds to a path 
\[
u_0 \fl (u_0\mid\vec v_1) \fl \cdots \fl (u_0\mid\vec v_1\mid\cdots\mid\vec v_n).
\]

\subsubsection{Examples}
\Item{i} Consider the following binary quadratic $1$-polygraph
\[
X:=\left\langle
x\in X_0(2) \;\left|\; \inputtree{[x[x[1][2]][3]]}\fl0,\; \inputtree{[x[x[1][3]][2]]}\fl0,\;\inputtree{[x[1][x[2][3]]]} \fl 0\;
\right.\right\rangle.
\]
It has $15$ critical branchings, which correspond to all possible critical branchings in the quadratic binary case. Let us draw the part of the directed graph $\Gr(X)$ corresponding to~$\Ov(X)_2$:
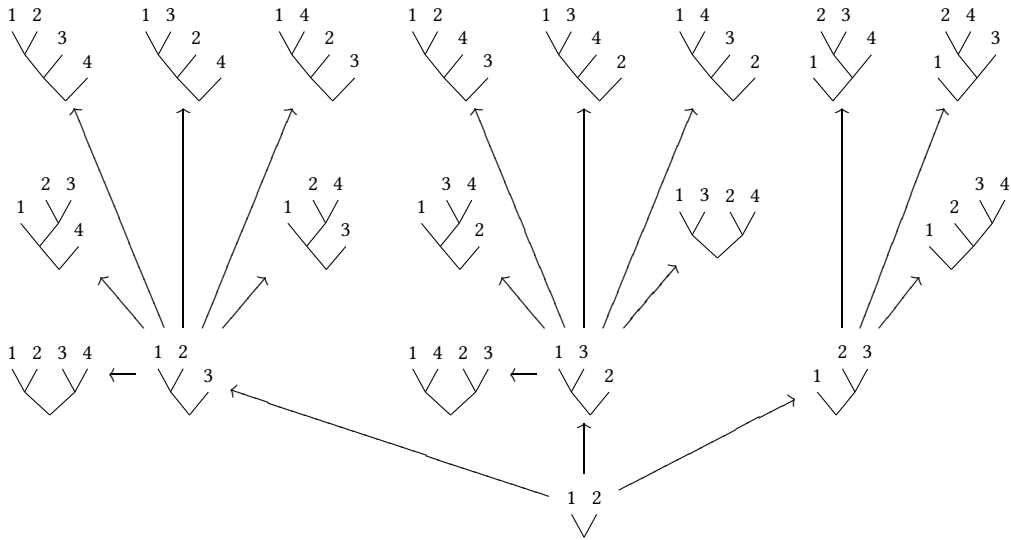
\[
\xymatrix @C = 1ex @R = 4ex{
\skele{[[[[1][2]][3]][4]]}
&
\skele{[[[[1][3]][2]][4]]}
&
\skele{[[[[1][4]][2]][3]]}
&
\skele{[[[[1][2]][4]][3]]}
&
\skele{[[[[1][3]][4]][2]]}
&
\skele{[[[[1][4]][3]][2]]}
&
\skele{[[1][[[2][3]][4]]]}
&
\skele{[[1][[[2][4]][3]]]}
\\
\skele{[[[1][[2][3]]][4]]}
&&
\skele{[[[1][[2][4]]][3]]}
&
\skele{[[[1][[3][4]]][2]]}
&&
\skele{[[[1][3]][[2][4]]]}
&&
\skele{[[1][[2][[3][4]]]]}
\\
\skele{[[[1][2]][[3][4]]]}
&
\skele{[[[1][2]][3]]}
  \ar[uul]
  \ar[uu]
  \ar[uur]
  \ar[ul]
  \ar[ur]
  \ar[l]
&&
\skele{[[[1][4]][[2][3]]]}
&
\skele{[[[1][3]][2]]}
  \ar[uul]
  \ar[uu]
  \ar[uur]
  \ar[ul]
  \ar[ur]
  \ar[l]
&&
\skele{[[1][[2][3]]]}
  \ar[uu]
  \ar[uur]
  \ar[ur]
\\
&&&&
\skele{[[1][2]]}
  \ar[ulll]
  \ar[u]
  \ar[urr]
}
\]
Every internal vertex of every tree monomial is $x$, so we omit them.

\Item{ii} Next, consider the following binary cubic $1$-polygraph
\[
X := \left\langle
x \in X_0(2)
\;\left|\;
\inputtree{[x [x [x [1] [2]] [3]] [4]]} \xrightarrow{\alpha} 0,
\inputtree{[x [x [1] [x [2] [3]]] [4]]} \xrightarrow{\beta} 0,
\inputtree{[x [1] [x [x [2] [3]] [4]]]} \xrightarrow{\gamma} 0,
\inputtree{[x [1] [x [2] [x [3] [4]]]]} \xrightarrow{\delta} 0
\right.\right\rangle,
\]
and consider the path-lexicographic order $\prec_{pl}$ on $1$-monomials where $\alpha \prec \beta \prec \gamma \prec \delta$ (see \textsection~\ref{SSS:PathLexicographicMonomialOrder}). Let us draw a part of the directed graph $\Gr(X)$ around $s_0(\alpha)$:
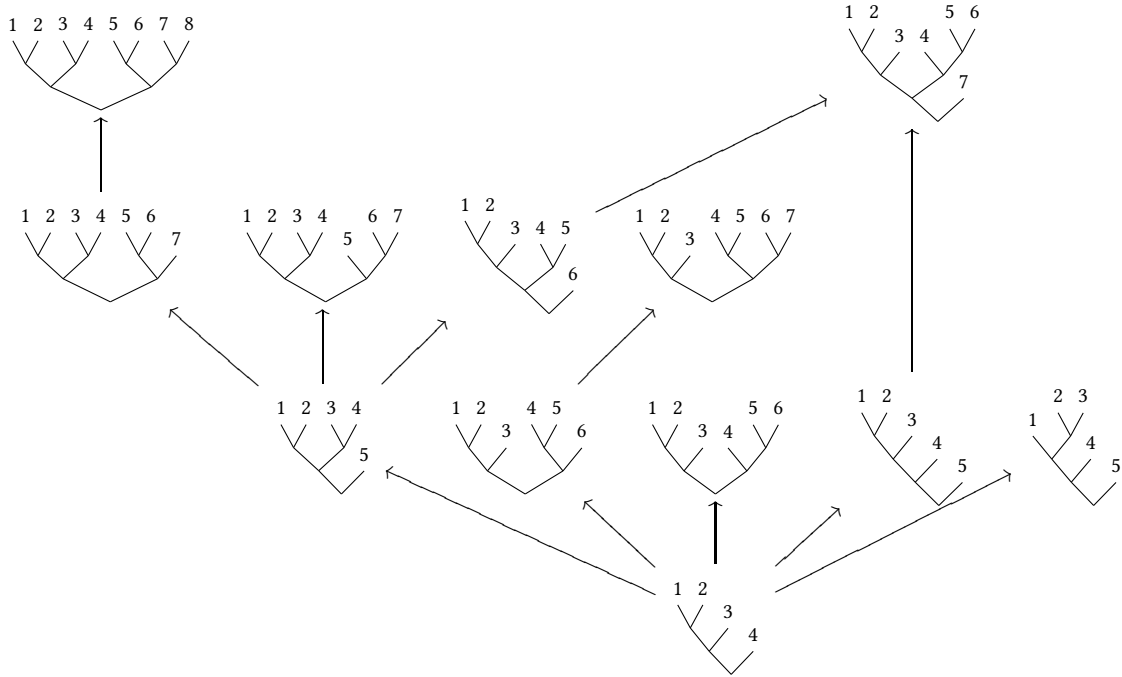
\[
\xymatrix @C = 1ex @R = 4ex{
\skele{[[[[1] [2]] [[3] [4]]] [[[5] [6]] [[7] [8]]]]}
&&&&
\skele{[[[[[1] [2]] [3]] [[4] [[5] [6]]]] [7]]}
\\
\skele{[[[[1] [2]] [[3] [4]]] [[[5] [6]] [7]]]}
  \ar[u]
&
\skele{[[[[1] [2]] [[3] [4]]] [[5] [[6] [7]]]]}
&
\skele{[[[[[1] [2]] [3]] [[4] [5]]] [6]]}
  \ar[urr]
&
\skele{[[[[1] [2]] [3]] [[[4] [5]] [[6] [7]]]]}
\\
&
\skele{[[[[1] [2]] [[3] [4]]] [5]]}
  \ar[ul]
  \ar[u]
  \ar[ur]
&
\skele{[[[[1] [2]] [3]] [[[4] [5]] [6]]]}
  \ar[ur]
&
\skele{[[[[1] [2]] [3]] [[4] [[5] [6]]]]}
&
\skele{[[[[[1] [2]] [3]] [4]] [5]]}
  \ar[uu]
&
\skele{[[[[1] [[2] [3]]] [4]] [5]]}
\\
&&&
\skele{[[[[1] [2]] [3]] [4]]}
  \ar[ull]
  \ar[ul]
  \ar[u]
  \ar[ur]
  \ar[urr]
}
\]
This drawing is not exhaustive, but presents some interesting phenomena. For example, due to our choice of order $\prec_{pl}$, the top left $4$-overlapping can only be obtained in one way, by adding $\alpha$, $\beta$, $\gamma$, and $\delta$ in order. In addition, the top right monomial can be reached from $s_0(\alpha)$ by a path of length $2$ or $3$, and so corresponds to a $3$-overlapping and a $4$-overlapping, depending on if $\beta$ is present or not.

\begin{theorem}[Overlapping polygraphic resolution]
\label{T:MainTheoremD}
Let $X$ be a reduced, convergent, left-monomial $1$-polygraph and $\iota$ the unital section sending every monomial to its reduced form. Then there exist a unique $\omega$-polygraph structure on $\Ov(X)$ and a unique right $\iota$-contraction $\sigma$ of $\Ov(X)$ such that, for all $n$-overlappings $u_n$ of $\Ov(X)$ and reduced $0$-monomials $\vec v_{n+1}$ of $\ftree{X_0}$,
\begin{eqn}{equation}
\label{E:ContractionCondition}
\sigma(u_n\mid\vec v_{n+1}) =
\begin{cases}
u_n\tri\vec v_{n+1} & \text{if }u_n\tri\vec v_{n+1}\in\Ov(X)_{n+1}, \\
\text{an identity} & \text{if }C(u_n,\vec v_{n+1})=\varnothing, \\
\sigma(u_n\mid\vec v_{n+1}) & \text{otherwise (tautological condition).}
\end{cases}
\end{eqn}
As a consequence, $\Ov(X)$ is a polygraphic resolution of the operad $\cl X$.
\end{theorem}

\begin{proof}
By induction on $n\geq 0$, we simultaneously construct the source and target maps of the $\omega$-polygraph structure on the $(n+1)$-generators of $\Ov(X)$ and the right $\iota$-contraction $\sigma:\fopShu{\Ov(X)}_n\fl\fopShu{\Ov(X)}_{n+1}$. By Lemma \ref{L:EssentialValues}, it suffices to define $\sigma$ on the $\iota$- and $\sigma$-essential $n$-monomials of the $\omega$-operad $\fopShu{\Ov(X)}$.

Let $n=0$. The $\iota$-essential $0$-monomials of $\fopShu{\Ov(X)}$ are the $(u_0\mid\vec v_1)$ where $u_0$ is a $0$-generator of~$X$ and the $v_{1,i}$ are reduced $0$-monomials of $\fopShu{X}_0$ such that $(u_0\mid\vec v_1)$ is not reduced. By \eqref{E:ContractionForm}, it suffices to define $\sigma(u_0\mid\vec v_1)$ such that $s_0\sigma(u_0\mid\vec v_1) = (u_0\mid\vec v_1)$ and $t_0\sigma(u_0\mid\vec v_1) = \rep{(u_0\mid\vec v_1)}$.
If $u_0\tri v_1$ is a $1$-overlapping, then we set
\begin{align*}
s_0(u_0\tri\vec v_1) & := (u_0\mid\vec v_1), & t_0(u_0\tri\vec v_1) & := \rep{(u_0\mid\vec v_1)},
\end{align*}
and the first case of \eqref{E:ContractionCondition} imposes $\sigma(u_0\mid\vec v_1):=u_0\tri\vec v_1$. Otherwise, since we have supposed $(u_0\mid\vec v_1)$ reducible, $C(u_0,\vec v_1)$ is nonempty. Let $\vec w_1=\max C(u_0,\vec v_1)$ and write $(u_0\mid\vec w_1\mid\vec w_2)=(u_0\mid\vec v_1)$. Then $u_0\tri\vec v_1\in\Ov(X)_1$ and $E(u_0\tri\vec w_1,\vec w_2)=\varnothing$, so by the second case of \eqref{E:ContractionCondition}, the (not yet defined) $2$-cell $\sigma(u_0\tri\vec w_1\mid\vec w_2)$ is an identity. By \eqref{E:ContractionForm}, we know that the target of $\sigma(u_0\tri\vec w_1\mid\vec w_2)$ is $\sigma(s_0(u_0\tri\vec w_1)\mid\vec w_2)=\sigma(u_0\mid\vec v_1)$. Thus we set
\[
\sigma(u_0\mid\vec v_1):=t_0\sigma(u_0\tri\vec w_1\mid\vec w_2) = s_0\sigma(u_0\tri\vec w_1\mid\vec w_2) = (u_0\tri\vec w_1\mid\vec w_2)\star_0\sigma(\rep{(u_0\mid\vec w_1)}\mid\vec w_2).
\]
Since $X$ is terminating, we define $\sigma(\rep{(u_0\mid\vec w_1)}\mid\vec w_2):(\rep{(u_0\mid\vec w_1)}\mid\vec w_2) \fl \rep{(u_0\mid\vec w_1\mid\vec w_2)}$ by well-founded induction on $\prec_{X_1}$, so this definition is licit.

Now let $n\geq 1$. The essential $n$-cells of $\fopShu{\Ov(X)}$ are the $(u_n\mid\vec v_{n+1})$ where $u_n$ is an $n$-overlapping and the $v_{n+1,i}$ are reduced $0$-monomials of $\fopShu{X_0}$ such that $(u_n\mid\vec v_{n+1})$ is not $\sigma$-reduced. Denote the branches of $u_n$ by $(\Gamma_1[\alpha_1],\ldots,\Gamma_n[\alpha_n])$. We distinguish the three cases of \eqref{E:ContractionCondition}. The induction step for the $\omega$-polygraph structure on $\Ov(X)$ is entirely contained within the first case.

\subsubsection*{First case}
First, suppose that $u_n\tri\vec v_{n+1}$ is an $(n+1)$-overlapping. Since condition \eqref{E:ContractionCondition} imposes $u_n\tri\vec v_{n+1}=\sigma(u_n\mid\vec v_{n+1})$, and \eqref{E:ContractionForm} gives us the source and target of the (not yet defined) $(n+1)$-cell $\sigma(u_n\mid\vec v_{n+1})$, we set
\begin{align*}
s_n(u_n\tri\vec v_{n+1}) & :=(u_n\mid\vec v_{n+1})-(t_{n-1}(u_n)\mid\vec v_{n+1})+\sigma(t_{n-1}(u_n)\mid\vec v_{n+1}),\\
t_n(u_n\tri\vec v_{n+1}) & := \sigma(s_{n-1}(u_n)\mid\vec v_{n+1}),
\end{align*}
which are indeed globular, and define $\sigma(u_n\mid\vec v_{n+1}):=u_n\tri\vec v_{n+1}$. This gives us the polygraphic structure on the $(n+1)$-overlappings.

\subsubsection*{Second case}
Next, suppose that $C(u_n,\vec v_{n+1})=\varnothing$. Writing $u_n=u_{n-1}\tri\vec v_n$, we make the following observations:
\begin{itemize}
\item the pair $(u_{n-1},(\vec v_n\mid\vec v_{n+1}))$ is not an $n$-overlapping.
\item $C(u_{n-1},(\vec v_n\mid\vec v_{n+1}))$ is nonempty, since it includes $\vec v_n$.
\item The $0$-monomials of $(\vec v_n\mid\vec v_{n+1})$ are reduced. Indeed, if not, then there would exist $\Gamma_{n+1}[\alpha_{n+1}] \succ_{pl} (\Gamma_n[\alpha_n]\mid\vec v_{n+1})$ in $E(u_n,\vec v_{n+1})$, which contradicts the fact that $C(u_n,\vec v_{n+1})$ is empty.
\end{itemize}
In particular, the third observation says that $(u_{n-1}\mid\vec v_n\mid\vec v_{n+1})$ is an essential $(n-1)$-monomial. Thus we are in the third case of the induction hypothesis. Following the calculations of the induction hypothesis in the third case below, let $\vec w_n$ be the maximal element of $C(u_{n-1},(\vec v_n\mid\vec v_{n+1}))$ and let $\vec w_{n+1}$ be $0$-monomials such that $(u_{n-1}\mid\vec w_n\mid\vec w_{n+1})=(u_{n-1}\mid\vec v_n\mid\vec v_{n+1})$. Then, by induction, the source and target of $\sigma(u_{n-1}\tri\vec w_n\mid\vec w_{n+1})$ are equal.

Suppose by contradiction that $\vec v_n\neq\vec w_n$. Let $\Gamma_n'[\alpha_n']$ be the last branch associated to $\vec w_n$. Then there exists a $(n+1)$-overlapping $u_{n-1}\tri\vec v_n\tri\vec w_{n+1}'\in\Ov(X)_{n+1}$, whose branches are 
\[
\{(\Gamma_1[\alpha_1] \mid \vec w_{n + 1}'), \ldots, (\Gamma_n[\alpha_n] \mid \vec w_{n+1}'), (\Gamma_n'[\alpha_n'] \mid \vec v_n')\}
\] 
where $\vec v_n'$ is the appropriate list of $0$-monomials. Thus $\vec w_{n+1}'\in C(u_n,\vec v_{n+1})$, which contradicts the hypothesis that $C(u_n,\vec v_{n+1})$ is empty. Therefore $\vec v_n=\vec w_n$, and we conclude that the source and target of $\sigma(u_{n-1}\tri\vec v_n\mid\vec v_{n+1})=\sigma(u_n\mid\vec v_{n+1})$ are equal, allowing us to define $\sigma(u_n\mid\vec v_{n+1})$ as an identity.

\subsubsection*{Third case}
Otherwise, $C(u_n,\vec v_{n+1})$ is nonempty. Let $\vec w_{n+1}$ be its maximal element, $\Gamma_{n+1}[\alpha_{n+1}]$ the associated $1$-monomial, and write $(u_n\mid\vec v_{n+1})=(u_n\mid\vec w_{n+1}\mid\vec w_{n+2})$. Then 
\[
E(u_n\tri\vec w_{n+1},\vec w_{n+2}) \subseteq \{\Gamma'[\alpha']\in E(u_n,\vec v_{n+1}) \mid \Gamma'[\alpha']\succ_{pl} \Gamma_{n+1}[\alpha_{n+1}]\}=\varnothing. 
\]
In addition, the monomials $\vec w_{n+2}$ are reduced, so this is exactly the condition of the second case, so we have the constraint that the source and target of the (not yet defined) $(n+2)$-cell $\sigma(u_n\tri\vec w_{n+1}\mid\vec w_{n+2})$ are equal.

The rest of this case is rather technical, so we summarize our strategy here. We prove and use Lemma~\ref{L:ExplicitBoundaryOfContraction} in order to get an explicit expression of $\partial \sigma(u_n \tri \vec w_{n + 1} \mid \vec w_{n + 2})$, which must be equal to $0$, as we have observed. This expression consists of many terms, including $\sigma(u_n \mid \vec v_{n + 1})$, which is the term that we have to define. We then proceed by well-founded induction on the terms of $\partial \sigma(u_n \tri \vec w_{n + 1} \mid \vec w_{n + 2})$ to define $\sigma(u_n \mid \vec v_{n + 1})$ using the other terms.

Let $k\geq 1$ and $\vec v_0,\vec v_1\ldots,\vec v_k$ $0$-cells of $\fopShu{X}$ such that $\{(v_{0,i} \mid \vec v_1^i \mid \cdots \mid \vec v_n^i)\}_i$ is a well-defined list of $0$-monomials of $\fopShu{X}$ where $\vec v_0$ is the list of roots and, for $\ell\in\range{k}$, $\vec v_\ell^i$ is the sublist of $\vec v_\ell$ of ancestor~$v_{0,i}$. We denote this list by $(\vec v_0 \mid \cdots \mid \vec v_n)$. Similarly, we denote by $\rep{(\vec v_0\mid\vec v_1)}$ the list of reduced $0$-cells $\rep{(v_{0,i}\mid\vec v_1^i)}$. Finally, we denote by $(\vec v_0\|\cdots\|\vec v_k)$ the list of $k$-cells
\[
\underbrace{\sigma(\sigma(\cdots\sigma(\sigma}_k(v_{0,i}\mid\vec v_1^i)\mid\vec v_2^i)\mid\cdots\vec v_{k-1}^i)\mid\vec v_k^i).
\]
Note that, if $u_0\tri\vec v_1\tri\cdots\tri\vec v_k$ is an $k$-overlapping, then $u_0\tri\vec v_1\tri\cdots\tri\vec v_k=(u_0\|\vec v_1\|\cdots\|\vec v_k)$.

\subsubsection{Lemma}
\label{L:ExplicitBoundaryOfContraction}
\emph{
For $n\geq 2$ and $\vec v_0,\ldots,\vec v_n$ $0$-monomials of $\fopShu{X}_0$, we have the equality of $(n-1)$-cells
\begin{align*}
\partial(\vec v_0\|\cdots\|\vec v_n) ={} & ((\vec v_0\|\cdots\|\vec v_{n-1})\mid\vec v_n) \\
& + \sum_{k=1}^n(-1)^k(\vec v_0\|\cdots\|\rep{(\vec v_{n-k}\mid\vec v_{n-k+1})}\|\cdots\|\vec v_n) \\
& +  (-1)^{n+1} (\vec v_0\mid(\vec v_1\|\cdots\|\vec v_n)) + 1_c,
\end{align*}
where $\partial=s-t$ and $c$ is some $(n-2)$-cell.
}
\begin{proof}
Proceed by induction on $n\geq 2$. According to \eqref{E:ContractionForm}, for every $n$-cell $a$,
\[
\partial\sigma(a) = a - \sigma(\partial a) + 1_c,
\]
where $c=-t_{n-1}(a)$ is an $(n-1)$-cell. For $n=2$, applying this equality to $(\vec v_0\|\vec v_1\|\vec v_2)$ gives
\begin{align*}
\partial(\vec v_0\|\vec v_1\|\vec v_2) ={} & \partial\sigma(\sigma(\vec v_0\mid\vec v_1)\mid\vec v_2) \\
={} & (\sigma(\vec v_0\mid\vec v_1)\mid\vec v_2)-\sigma(\partial\sigma(\vec v_0\mid\vec v_1)\mid\vec v_2)+1_c \\
={} & (\sigma(\vec v_0\mid\vec v_1)\mid\vec v_2)+\sigma(\rep{(\vec v_0\mid\vec v_1)}\mid\vec v_2)-\sigma(\vec v_0\mid\vec v_1\mid\vec v_2)+1_c \\
={} & (\sigma(\vec v_0\mid\vec v_1)\mid\vec v_2)+\sigma(\rep{(\vec v_0\mid\vec v_1)}\mid\vec v_2)-(\vec v_0\mid\sigma(\vec v_1\mid\vec v_2))-\sigma(\vec v_0\mid\rep{(\vec v_1\mid\vec v_2)})+1_{c'}+1_c \\
={} & ((\vec v_0\|\vec v_1)\mid\vec v_2) + (-1)^1(\vec v_0\|\rep{(\vec v_1\mid\vec v_2)})+ (-1)^2(\rep{(\vec v_0\mid\vec v_1)}\|\vec v_2) +(-1)^3(\vec v_0\mid(\vec v_1\|\vec v_2))+1_{c+c'}.
\end{align*}
Let $n\geq 2$. Recall that, for all $(n-1)$-cells $u$ and $0$-cells $\vec v$,
\[
\sigma(u\mid\vec v)=(s_0(u)\mid\sigma(\vec v))\star_0(u\mid\rep{\vec v})=\sigma(u\mid\rep{\vec v})+1_c
\]
with $c$ an $(n-1)$-cell. We calculate
\begin{align*}
\partial(\vec v_0\|\cdots\|\vec v_{n+1}) ={} & \partial\sigma((\vec v_0\|\cdots\|\vec v_n)\mid\vec v_{n+1}) \\
={} & ((\vec v_0\|\cdots\|\vec v_n)\mid\vec v_{n+1})-\sigma(\partial(\vec v_0\|\cdots\|\vec v_n)\mid\vec v_{n+1})+1_c \\
={} & ((\vec v_0\|\cdots\|\vec v_n)\mid\vec v_{n+1})-\sigma((\vec v_0\|\cdots\|\vec v_{n-1})\mid\vec v_n\mid\vec v_{n+1}) \\
& -\sum_{k=1}^n(-1)^k\sigma((\vec v_0\|\cdots\|\rep{(\vec v_{n-k}\mid\vec v_{n-k+1})}\|\cdots\|\vec v_n)\mid\vec v_{n+1}) \\
& -(-1)^{n+1}\sigma(\vec v_0\mid(\vec v_1\|\cdots\|\vec v_n)\mid\vec v_{n+1}) + \sigma(1_{c'}) + 1_c \\
={} & ((\vec v_0\|\cdots\|\vec v_n)\mid\vec v_{n+1})- \sigma((\vec v_0\|\cdots\|\vec v_{n-1})\mid\rep{(\vec v_n\mid\vec v_{n+1})})-1_{c''} \\
& -\sum_{k=1}^n(-1)^k(\vec v_0\|\cdots\|\rep{(\vec v_{n-k}\mid\vec v_{n-k+1})}\|\cdots\|\vec v_{n+1}) \\
&  -(-1)^{n+1}(\vec v_0\mid\sigma((\vec v_1\|\cdots\|\vec v_n)\mid\vec v_{n+1})) + 1_{\sigma(c')} + 1_c \\
={} & ((\vec v_0\|\cdots\|\vec v_n)\mid\vec v_{n+1}) \\
& +\sum_{k=1}^{n+1}(-1)^k(\vec v_0\|\cdots\|\rep{(\vec v_{n-k+1}\mid\vec v_{n-k+2})}\|\cdots\|\vec v_{n+1}) \\
& + (-1)^{n+2}(\vec v_0\mid(\vec v_1\|\cdots\|\vec v_{n+1})) + 1_{c+\sigma(c')-c''},
\end{align*}
which concludes the induction step, and the proof of the lemma.
\end{proof}

Writing $u_n=u_0\tri\vec v_1\tri\cdots\tri\vec v_n$, we apply the lemma to $\sigma(u_n\tri\vec w_{n+1}\mid\vec w_{n+2})=(u_0\|\vec v_1\|\cdots\|\vec v_n\|\vec w_{n+1}\|\vec w_{n+2})$ to get the equation of $(n+1)$-cells
\begin{align*}
\partial\sigma(u_n\tri\vec w_{n+1}\mid\vec w_{n+2}) = 0 ={} & ((u_0\|\vec v_1\|\cdots\|\vec v_n\|\vec w_{n+1})\mid\vec w_{n+2}) \\
& + \sum_{k=1}^{n+2}(-1)^k(u_0\|\vec v_1\|\cdots\|\rep{(\vec v_{n-k+2}\mid\vec v_{n-k+3})}\|\cdots\|\vec v_n\|\vec w_{n+1}\|\vec w_{n+2}) \\
& +  (-1)^{n+3} (u_0\mid(\vec v_1\|\cdots\|\vec v_n\|\vec w_{n+1}\|\vec w_{n+2})) + 1_c,
\end{align*}
where $c$ is an $n$-cell. On the righthand side, the $(n+1)$-cell
\[
(u_0\|\vec v_1\|\cdots\|\vec v_n\|\rep{(\vec w_{n+1}\mid\vec w_{n+2})})=\sigma(u_n\mid\vec v_{n+1})
\]
appears. We want to define this $(n+1)$-cell using the other $(n+1)$-cells that appear, that is,
\begin{eqn}{equation}
\label{E:OtherCells}
\begin{gathered}
(u_0\|\vec v_1\|\cdots\|\vec v_n\| w_{n+1}),\quad(\vec v_1\|\cdots\|\vec v_n\|\vec w_{n+1}\|\vec w_{n+2}), \\
(u_0\|\vec v_1\|\cdots\|\rep{(\vec v_{n-k+2}\mid\vec v_{n-k+3})}\|\cdots\|\vec v_n\|\vec w_{n+1}\|\vec w_{n+2}),\quad k\in\{2,\ldots,n+2\}.
\end{gathered}
\end{eqn}
We define a well-founded order $\prec$ on $(n+1)$-cells of the form $(u_0\|\vec v_1\|\cdots\|\vec v_n)$ by setting $(u_0\|\vec v_1\|\cdots\|\vec v_n)\prec (u_0'\|\vec v_1'\|\cdots\|\vec v_n')$ if
\begin{enumerate}
\item $T(u_0\mid\vec v_1\mid\cdots\mid\vec v_n)$ is a proper submonomial of $T(u_0'\mid\vec v_1'\mid\cdots\mid\vec v_n')$, or
\item $T(u_0\mid\vec v_1\mid\cdots\mid\vec v_n)=T(u_0'\mid\vec v_1'\mid\cdots\mid\vec v_n')$ and there exist $i,j$ such that $u_0=u_0',\vec v_1=\vec v_1',\ldots,\vec v_{i-1}=\vec v_{i-1}'$, $v_{i,1}=v_{i,1}',\ldots,v_{i,j-1}=v_{i,j-1}'$, and the weight of $v_{i,j}$ is less than that of $v_{i,j}'$, or
\item there exists a positive $1$-cell $f:(u_0'\|\vec v_1'\|\cdots\|\vec v_n')\fl b$ of $\fopShu{X}$ such that $(u_0\|\vec v_1\|\cdots\|\vec v_n)\in\supp(b)$.
\end{enumerate}
The relation $\prec$ is an order because the $1$-polygraph $X$ is supposed reduced (so we cannot rewrite a $0$-monomial into a larger $0$-monomial). The relation $\prec$ is well-founded because every sequence $((u_0^i\|\vec v_1^i\|\cdots\|\vec v_n^i))_{i\geq 0}$ that decreases for $\prec$ can be rearranged into the concatenation of a decreasing sequence for {\bf iii)} followed by a decreasing sequence for the lexicographic order induced by {\bf i)} and {\bf ii)} (if we can rewrite a submonomial of a $0$-monomial, then we can rewrite the $0$-monomial following the same rule).

We initialize our well-founded induction on the $(n+1)$-overlappings, since $u_0\tri\vec v_1\tri\cdots\tri\vec v_{n+1}=(u_0\|\vec v_1\|\cdots\|\vec v_{n+1})$ is already defined. We then check that all of the $(n+1)$-cells of \eqref{E:OtherCells} are smaller than $\sigma(u_n\mid\vec v_{n+1})$ for the order $\prec$: $(u_0\|\vec v_1\|\cdots\|\vec v_n\|\vec w_{n+1})$ and $(\vec v_1\|\cdots\|\vec v_n\|\vec w_{n+1}\|\vec w_{n+2})$ satisfy {\bf i)}, and $(u_0\|\vec v_1\|\cdots\|\rep{(\vec v_{n-k+2}\mid\vec v_{n-k+3})}\|\cdots\|\vec v_n\|\vec w_{n+1}\|\vec w_{n+2})$ satisfies {\bf ii)} if $(\vec v_{n-k+2}\mid\vec v_{n-k+3})$ is reduced, and {\bf iii)} otherwise, by confluence of $X$. Thus we can define $\sigma(u_n\mid\vec v_{n+1})$ by well-founded induction.

Finally, by Proposition~\ref{P:AcyclicPolygraph}, the $\omega$-polygraph $\Ov(X)$ is acyclic. Since $X$ is reduced, by the discussion of \textsection~\ref{SSS:ReducedLowDimension}, $\Ov(X)_{\leq 1}$ coincides with $X$. Therefore $\Ov(X)$ is a polygraphic resolution of the operad~$\cl X$.
\end{proof}

\begin{corollary}
Let $X$ be an essentially confluent reduced, terminating, left-monomial $1$-polygraph. Then there exists a $\omega$-polygraph structure on $\Ov(X)$ making it a polygraphic resolution of $\cl X$.
\end{corollary}
\begin{proof}
By Theorem \ref{T:CriticalBranchings}, the polygraph $X$ is convergent. Thus, following Theorem \ref{T:MainTheoremD}, $\Ov(X)$ is equipped with a $\omega$-polygraph structure and is a polygraphic resolution of $\cl X$.
\end{proof}

\subsubsection{Coherent presentations from convergence}
In \cite{Squier94}, Squier showed how to compute a coherent presentation of a monoid from a convergent one. This construction is described in the case of associative algebras in {\cite[Thm 4.3.2]{GuiraudHoffbeckMalbos19}}, and in the case of shuffle operads by using the following result. For a convergent left-monomial $1$-polygraph $X$, and a cellular extension~$Y$ of~$\fopShu{X}_1$ that contains a $2$-generator $A_{f,g}$ of shape
\[
\vcenter{\xymatrix @R=0.5em {
& b
	\ar@/^/ [dr] ^-{h}
	\ar@2 []!<0pt,-12.5pt>;[dd]!<0pt,12.5pt> ^-*+{A_{f,g}}
\\
a
	\ar@/^/ [ur] ^-{f}
	\ar@/_/ [dr] _-{g}
&& d
\\
& c
	\ar@/_/ [ur] _-{k}
}}
\]
with~$h$ and~$k$ positive $1$-cells of~$\fopShu{X}_1$, for every critical branching~$(f,g)$ of~$X$, then the $2$-polygraph $(X,Y)$ is acyclic. The $2$-generator $A_{f,g}$ is called a \emph{generating confluence} associated to the critical branching~$(f,g)$. Note that such a generating confluence depends on the choice of the positives cells $h$ and $k$ and the orientation of the $2$-cell $A_{f,g}$. 
The proof of this result is done in two steps. First, we show that the $2$-polygraph $(X,Y)$ is coherently confluent, then we prove acyclicity of the cellular extension~$Y$, see~{\cite[Thm. 4.3.2]{GuiraudMalbos18}}.

The proof of Theorem~\ref{T:MainTheoremD} makes explicit the generating confluences, which can be represented as follows:
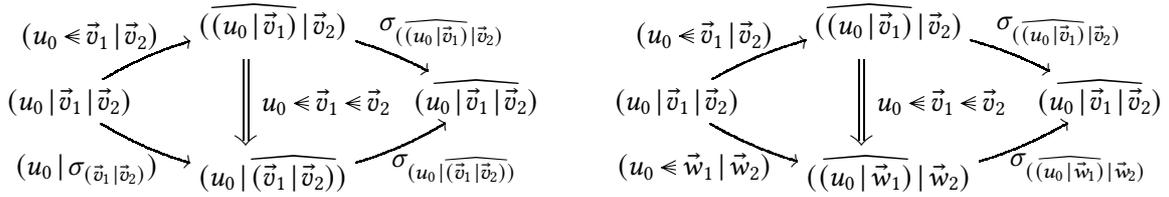
\[
\xymatrix @C=1.7em @R=0.8em {
& (\rep{(u_0\!\mid\!\vec v_1)}\!\mid\!\vec v_2)
	\ar@/^/ [dr] ^-{\quad\sigma_{(\rep{(u_0|\vec v_1)}|\vec v_2)}}
	\ar@2 []!<-10pt,-12.5pt>;[dd]!<-10pt,12.5pt> ^-*+{u_0\tri\vec v_1\tri\vec v_2}
\\
(u_0\!\mid\!\vec v_1\!\mid\!\vec v_2)
	\ar@/^/ [ur] ^-{(u_0\tri\vec v_1\!\mid\!\vec v_2)\quad}
	\ar@/_/ [dr] _-{(u_0\!\mid\!\sigma_{(\vec v_1|\vec v_2)})\quad}
&& \rep{(u_0\!\mid\!\vec v_1\!\mid\!\vec v_2)}
\\
& (u_0\!\mid\!\rep{(\vec v_1\!\mid\!\vec v_2)})
	\ar@/_/ [ur] _-{\quad\sigma_{(u_0|\rep{(\vec v_1|\vec v_2)})}}
}
\qquad
\xymatrix @C=1.7em @R=0.8em {
& (\rep{(u_0\!\mid\!\vec v_1)}\!\mid\!\vec v_2)
	\ar@/^/ [dr] ^-{\quad\sigma_{(\rep{(u_0|\vec v_1)}|\vec v_2)}}
	\ar@2 []!<-10pt,-12.5pt>;[dd]!<-10pt,12.5pt> ^-*+{u_0\tri\vec v_1\tri\vec v_2}
\\
(u_0\!\mid\!\vec v_1\!\mid\!\vec v_2)
	\ar@/^/ [ur] ^-{(u_0\tri\vec v_1\!\mid\!\vec v_2)\quad}
	\ar@/_/ [dr] _-{(u_0\tri\vec w_1\!\mid\!\vec w_2)\quad}
&& \rep{(u_0\!\mid\!\vec v_1\!\mid\!\vec v_2)}
\\
& (\rep{(u_0\!\mid\!\vec w_1)}\!\mid\!\vec w_2)
	\ar@/_/ [ur] _-{\quad\sigma_{(\rep{(u_0|\vec w_1)}|\vec w_2)}}
}
\]
where $u_0\tri\vec v_1\tri\vec v_2$ is a $2$-overlapping of $X$.
If $(\vec v_1\mid\vec v_2)$ is reducible, we have the left diagram. If $(\vec v_1\mid\vec v_2)$ is reduced, we take $\vec w_1=\max C(u_0,(\vec v_1\mid\vec v_2))$ and get the right diagram.

\subsubsection{Examples}
\begin{enumerate}[{\bf i)}]
\item The terminating reduced $1$-polygraph $X$ defined in Example~\ref{E:xyz} does not have critical branchings. As a consequence, it is convergent and can be extended into a polygraphic resolution $\Ov(X)$, with $\Ov(X)_n$ empty for $n\geq 2$.
\item Following Proposition~\ref{P:HigherTietzeEquivalence}, if $X$ is an acyclic $\omega$-polygraph whose underlying $1$-polygraph $X_{\leq 1}$ is left-monomial and convergent, then $X$ is Tietze equivalent to the $\omega$-polygraph $\Ov(X_{\leq 1})$.
In particular, for every operad $P$, the $\omega$-polygraphs $\Std(P)$ and $\Ov(\Std(P)_{\leq 1})$ are Tietze equivalent.
\end{enumerate}

\subsection{Bimodule resolutions from polygraphic resolutions}
\label{SS:BimodulesResolution}

In this subsection, we show how to deduce the homology of a shuffle operad with coefficients in bimodules from a shuffle polygraphic resolution of the operad.

\subsubsection{Construction of a chain complex}
Let $X$ be a shuffle $\omega$-polygraph, and denote $P$ the shuffle operad presented by $X$. Denote by $\cl{u}$ the image of $u$ by the canonical projection $\pi_X$ defined in~\eqref{E:presentationShuffleOperad}. Consider the chain complex $(P\langle X\rangle_n)_{n\geq -1}$ in the category $\Bimod{P}$ generated by $X$, that is, for all $n\geq -1$, $P\langle X\rangle_n:=P\langle X_n\rangle$ is the free $P$-bimodule on $X_n$, where $X_{-1}$ is the unit indexed set $\mathbf{1}$ defined in\textsection~\ref{SSS:ShuffleCompositionOmegaSets}. The generators of $P\langle X \rangle_n$, are denoted $[x]$ for $x$ in $X_n$. The boundary map $\delta_n : P\langle X\rangle_{n+1} \fl P\langle X\rangle_n$ is defined as follows
\begin{enumerate}
\item For $x\in X_0(k)$, we set
\begin{eqn}{equation}
\delta_{-1}([x])=(\varepsilon\mid\overline x)-\sum_{i=1}^k(\overline x\mid\overline 1\cdots \underset{i}{\varepsilon}\cdots\overline 1).
\end{eqn}
\item Consider the derivation  $[\;]:\fopShu{X_0}\rightarrow P\langle X_0\rangle$ defined by induction on the weight of monomials $u$ in $\fopShu{X_0}$, by setting $[1]:=0$, $[x]:=x$ and
\[
[u\mid\vec v]:=([u]\mid\overline v_1\cdots\overline v_k)+\sum_{i=1}^k(\overline u\mid\overline v_1\cdots [v_i]\cdots\overline v_k).
\]
We set $\delta_0$ the \emph{Fox differential} defined for every $1$-generator $\alpha$ in $X_1$ by
\[
\delta_0([\alpha]):=[s_0(\alpha)]-[t_0(\alpha)].
\]
\item For $n\geq 1$, we define the map $[\;]:\fopShu{X_n}\fl P\langle X_n\rangle$ by setting, for
\[
f=\sum_{i=1}^p\lambda_i\Gamma_i[\alpha_i]+1_c
\] 
an $n$-cell of $\fopShu{X}$, where $\alpha_i\in X_n$ and $\Gamma_i$ is a one-hole context of $\ftree{X}_0$, 
\[
[f]=\sum_{i=1}^p\lambda_i\overline\Gamma_i\big[[\alpha_i]\big],
\] 
where $\overline\Gamma_i$ is the one-hole context of $P$ induced by the context $\Gamma_i$. Note that $[f]$ does not depend on the choice of decomposition, so $[\;]$ is well defined.
We set for every $(n+1)$-generator $A$ in $X_{n+1}$
\[
\delta_n([A]):=[s_n(A)]-[t_n(A)].
\]
\end{enumerate}
As a consequence of the globularity of the polygraph $X$, for all $n\geq -1$, we have $\delta_{n+1}\delta_n = 0$ and thus~$P\langle X \rangle$ forms a chain complex.

\begin{lemma}
For every $0$-monomial $u\in \ftree{X_0}(k)$, we have
\[
\delta_{-1}([u])=(\varepsilon\mid \overline u)-\sum_{i=1}^k(\overline u\mid\overline 1\cdots\undermath i\varepsilon\cdots\overline 1).
\]
\end{lemma}
\begin{proof}
Proceed by induction on the depth of the $0$-monomial $u$. The equality is true by definition for $x\in X_0$. For the induction step, consider $(u\mid\vec v)$ with $u\in X_0^\ell(k)$, $v_i\in X_0^\ell(\ell_i)$ for all $1\leq i\leq k$ :
\begin{align*}
\delta_{-1}([u\mid\vec v]) ={} & (\delta_{-1}([u])\mid\overline v_1\cdots\overline v_k)+\sum_{i=1}^k(\overline u\mid\overline v_1\cdots\delta_{-1}([v_i])\cdots\overline v_k) \\
={} & (\varepsilon\mid\overline u\mid\overline v_1\cdots\overline v_k)-\sum_{i=1}^k(\overline u\mid(\overline 1\mid\overline v_1)\cdots (\varepsilon\mid v_i)\cdots(\overline 1\mid\overline v_k)) \\
& +\sum_{i=1}^k(\overline u\mid\overline v_1\cdots(\varepsilon\mid\overline v_i)\cdots\overline v_k)-\sum_{i=1}^k\sum_{j=1}^{\ell_i}(\overline u\mid\overline v_1\cdots(\overline v_i\mid\overline 1\cdots\undermath j\varepsilon\cdots\overline 1)\cdots\overline v_k) \\
={} & (\varepsilon\mid\overline{(u\mid\vec v)})-\sum_{i=1}^{\ell_1+\cdots+\ell_k}(\overline{(u\mid\vec v)}\mid\overline 1\cdots\undermath i\varepsilon\cdots\overline 1).
\end{align*}
\end{proof}

\subsubsection{Trivial $P$-bimodule}
Define the \emph{trivial $P$-bimodule}, denoted by $\Omega_P$, as the free $P$-bimodule generated by the unit indexed set $\mathbf{1}$ quotiented by the relations 
\begin{eqn}{equation}
(\varepsilon\mid\overline u)=\sum_{i=1}^k\overline u\circ_i\varepsilon
\end{eqn}
for every $k\geq 1$ and $\overline u\in P(k)$. Every element of the $P$-bimodule $\Omega_P$ can be written as a linear combination of monomials of the form $\overline u\circ_i\varepsilon$ where $k\geq 1$, $\overline u\in P(k)$, and $1\leq i\leq k$.

\begin{proposition}
\label{P:Resolution}
Let $X$ be an acyclic shuffle $\omega$-polygraph and $P$ the shuffle operad presented by~$X$. Then the chain complex $P\langle X\rangle$ is a resolution of $\Omega_P$ in the category $\Bimod{P}$.
\end{proposition}
\begin{proof}
Note that $\Omega_P$ is exactly the cokernel of $\delta_{-1}$. Thus it suffices to show that the chain complex $P\langle X\rangle$ is exact.

Let us fix $\iota$ a unital section of $X$. Following Proposition~\ref{P:AcyclicPolygraph}, the acyclicity of the polygraph $X$ implies that it admits a right $\iota$-contraction. Let $\sigma$ be such a right $\iota$-contraction.
We define the linear map $i_0:P\langle \mathbf 1\rangle\fl P\langle X_0\rangle$ by
\[
i_0(u\mid v_1\cdots\underset{i}{(\varepsilon\mid w)}\cdots v_k) := (u\mid v_1\cdots\underset{i}{[\rep w]}\cdots v_k),
\]
for $u,v_1,\ldots,\check v_i,\ldots,v_n,w\in P$, and, for $n\geq 1$, the linear map $i_n:P\langle X_{n-1}\rangle\fl P\langle X_n\rangle$ by
\[
i_n(u\mid v_1\cdots\underset{i}{([x]\mid w_1\cdots w_\ell)}\cdots v_k) := (u\mid v_1\cdots\underset{i}{[\sigma_{(x\mid\rep w_1\cdots\rep w_\ell)}]}\cdots v_k),
\]
 and $u,v_1,\ldots,\check v_i,\ldots,v_k,w_1,\ldots,w_\ell\in P$.
 Note that the linear maps $i_n$ are compatible with the left action of $P$. Hence, we prove that the maps $i_n$ define a contracting homotopy of the complex $P\langle X\rangle$, by showing that the identity $i_n\delta_{n-1}+\delta_ni_{n+1}=\id_{P\langle X_n\rangle}$ holds on generators of the $P$-bimodule $P\langle X_n\rangle$ as follows.
 
For $n=0$, on the one hand, we have
\begin{align*}
i_0\delta_{-1}([x]\mid w_1\cdots w_n) & = i_0(\varepsilon\mid\overline x\mid \vec w)-\sum_{i=1}^ni_0(\overline x\mid(\overline 1\mid w_1)\cdots(\varepsilon\mid w_i)\cdots (\overline 1\mid w_n)) \\
& = [\rep{(\overline x\mid\vec w)}]-\sum_{i=1}^n(\overline x\mid w_1\cdots [\rep w_i]\cdots w_n).
\end{align*}
On the other, we have
\begin{align*}
\delta_0i_1([x]\mid w_1\cdots w_n) & = \delta_0[\sigma_{(x\mid\rep w_1\cdots\rep w_n)}] \\
& =\delta_0[(x\mid\rep w_1\cdots\rep w_n)]-\delta_0[\rep{(x\mid\vec w)}] \\
& =([x]\mid w_1\cdots w_n)+\sum_{i=1}^n(\overline x\mid w_1\cdots [\rep w_i]\cdots w_n) -[\rep{(x\mid\vec w)}],
\end{align*}
proving the equality $\delta_0i_1+i_0\delta_{-1}=\id_{P\langle X_0\rangle}$.

For $n\geq 1$, by definition of the right $\iota$-contraction $\sigma$, we show that, for every $(n-1)$-cells $u,w_1,\ldots,w_n$ of $\fopShu{X}$,
\[
i_n[(u\mid w_1\cdots w_n)]=[\sigma_{(u\mid\rep w_1\cdots\rep w_n)}].
\]
Therefore, for every $n$-generator $A:a\rightarrow b$ in $X_n$, we have
\begin{align*}
i_n\delta_{n-1}([A]\mid w_1\cdots w_n) & =i_n[(a\mid\vec w)]-i_n[(b\mid\vec w)]=[\sigma_{(a\mid\rep w_1\cdots\rep w_n)}]-[\sigma_{(b\mid\rep w_1\cdots\rep w_n)}], \\
\delta_ni_{n+1}([A]\mid w_1\cdots w_n) & =[(A\mid\rep w_1\cdots\rep w_n)\star_0\sigma_{(b\mid\rep w_1\cdots\rep w_n)}]-[\sigma_{(a\mid\rep w_1\dots\rep w_n)}] \\
& =([A]\mid w_1\cdots w_n)+[\sigma_{(b\mid\rep w_1\cdots\rep w_n)}]-[\sigma_{(a\mid\rep w_1\cdots\rep w_n)}],
\end{align*}
proving that $i_n\delta_{n-1}+\delta_ni_{n+1}=\id_{P\langle X_n\rangle}$.
\end{proof}

\subsubsection{Homology of shuffle operads}
Recall that the Cartan-Eilenberg homology of a shuffle operad~$P$ with coefficients in a $P$-bimodule $A$ is defined by
\[
H^{CE}_\bullet(P,A) := \Tor^{\Bimod{P}}_\bullet(\Omega_P,A).
\]
In addition, the Quillen homology of $P$ is defined with coefficients in $\Abel(\ShuO / P)$, the category of abelian groups internal to $\ShuO/P$ \cite{Quillen68}. The category $\Abel(\ShuO / P)$ is equivalent to the category $\Bimod{P}$ of $P$-bimodules \cite{BauesJibladzeTonks97}, and we define the Quillen homology of $P$ with coefficients in a $P$-bimodule $A$ by setting
\[
H^Q_\bullet(P,A) := H_\bullet(\Abel(\mathcal{X}) \otimes_P A),
\]
where $\mathcal{X}$ is a simplicial cofibrant resolution of the operad $P$ in the category $ \ShuO/P$, and $\Abel(-) : \ShuO/P \fl \Abel(\ShuO/P) \approx \Bimod{P}$ is the abelianization functor. Following {\cite[Thm. 4.1]{Barr96}}, see also \cite[Thm. 6.2.1]{Barr02}, these two homologies are isomorphic up to shift in degree:
\[
H^Q_\bullet(P,A) \simeq H^{CE}_{\bullet+1}(P,A).
\]

\subsubsection{Finite homological type}
From Theorem~\ref{T:MainTheoremD} we deduce a generalization of Squier's homological finiteness condition~\cite{Squier87}, for finite convergence in the case of operads.
We say that a shuffle operad $P$ has \emph{finite homological type}, $FP_\infty$ for short, if the $P$-bimodule $\Omega_P$ has a resolution in $\Bimod{P}$ by finitely generated projective bimodules.
If $P$ admits a finite convergent presentation $X$, then by Theorem~\ref{T:MainTheoremD}, the overlapping polygraphic resolution $\Ov(X)$ is finite and the complex $P\langle \Ov(X)\rangle$ is a finitely generated free resolution of $\Omega_P$. Thus, $P$ has homological type $FP_\infty$.

\subsubsection{Minimal resolutions}
A \emph{minimal bimodule resolution} of an operad $P$ is a minimal free $P$-bimodule resolution $(A_\bullet,\delta)$ of its trivial $P$-bimodule $\Omega_P$. The minimal condition means that the sequence $(A_\bullet\otimes_P\K,\delta\otimes_P\id)$ has a null differential, where $\K$ denotes the $P$-bimodule concentrated in degree~$0$, whose left and right actions vanish.

\begin{proposition}
\label{P:concentratedDiagonal}
Let $X$ be an acyclic shuffle $\omega$-polygraph and $w:\Nb\fl\Nb\setminus\{0\}$ an increasing function such that $X_n$ is concentrated in weight~$w(n)$. Then $P\langle X\rangle$ is a minimal $P$-bimodule resolution of the operad $P$ presented by $X$.
\end{proposition}
\begin{proof}
The $1$-generators of $X$ are of homogeneous weight, so $P$ is equipped with a weight grading. Given an $n$-generator $u_n$ in $X_n$, we have $\delta_{n-1}[u_n]=\sum_i\lambda_i\cl{\Gamma_i}[u_{n-1,i}]+1_c$, where the $\lambda_i$ are scalars, the $\Gamma_i$ are one-hole contexts of $\ftree{X_0}$, the $u_{n-1,i}$ are $(n-1)$-overlappings, and $c$ is an $(n-2)$-cell in $\fopShu{X}_{n-2}$.
The map~$\delta_{n-1}$ preserves weight, and since the $(n-1)$-overlappings $u_{n-1,i}$ are of strictly smaller weight than $u_n$, it follows that the $\cl{\Gamma}_i$ are nontrivial.
As a consequence, tensoring over $P$ by $\K$ sends $\delta_{n-1}(u_n)$ to $0$, so $P\langle X\rangle$ is minimal.
\end{proof}

\subsubsection{Examples}
Consider the convergent $1$-polygraph $X$ with one $1$-generator $x\in X(2)$ of weight $1$ and one $2$-generator 
\[
\tree{[x[x][x]]} \fl 0
\]
Then $\Ov(X)_n$ is a polygraphic resolution concentrated in degree $2n+1$, so by Proposition~\ref{P:concentratedDiagonal} $P\langle X\rangle$ is a minimal bimodule resolution of the operad~$P$ presented by $X$.
Note however that, for the convergent $1$-polygraph with one $1$-generator $x\in X(1)$ of weight $1$ and one $2$-generator $x^3 \fl 0$, the resolution induced by $\Ov(X)$ is not minimal.

\subsection{Confluence and Koszulness}
\label{SS:ConfluenceKoszul}

In this subsection, we show that shuffle operads presented by a quadratic convergent $1$-polygraph are Koszul. Our result does not suppose that the rewriting rules are oriented with respect to a monomial order. In this way, it generalizes the result obtained by Dotsenko and Khoroshkin in \cite{DotsenkoKhoroshkin10} for shuffle operads with quadratic Gröbner bases. We first begin by recalling the Koszul property for operads.

\subsubsection{Koszul operads}
\label{SSS:KoszulOperads}
Let $P$ be a (connected and graded) symmetric operad. We denote by $\overline B(P)$ for the reduced bar complex on $P$. Recall from \cite[Def. 5.2.3]{Fresse04} that the \emph{Koszul complex} on $P$ is defined by
\[
K(P)_{(s)} := H_s(\overline B(P)_{(s)},\delta) = \ker(\delta : \overline B_s(P)_{(s)} \fl \overline B_{s-1}(P)_{(s)}),
\]
the second equality coming from the fact that $\overline B_n(P)_{(s)}=0$ when $n>s$, and where $(s)$ denotes the degree of $P$. By definition, the complex $K(P)_{(s)}$ is concentrated in degree $s$.
The operad $P$ is \emph{Koszul} if the inclusion morphism $K(P)\hookrightarrow\overline  B(P)$ is a quasi-isomorphism~\cite[Def. 5.2.8]{Fresse04}, or equivalently the homology of the reduced bar complex of $P$ is concentrated on the diagonal \cite[Thm. 5.3.3]{Fresse04}, that is,
\[
\ho_n(\overline B(P)_{(s)}) = 0,
\qquad
\text{for $n\neq s$.}
\]
Recall that the bar-cobar construction on $P$ is a resolution, whose abelianization is the reduced bar complex \cite[\textsection~1.1]{DotsenkoKhoroshkin13}, so that the operad $P$ is Koszul if, and only if, its Quillen homology is concentrated on the diagonal.
Finally, recall from {\cite[Cor.~1.5]{DotsenkoKhoroshkin13}} that for a symmetric operad $P$, there is an isomorphism
\[
\ho_\bullet(\overline B(P))^u \simeq \ho_\bullet(\overline B(P^u)).
\]
As a consequence, the Koszulness of a symmetric operad can be proved via its shuffle version as follows.

\begin{theorem}
\label{T:MainTheoremE}
Let $P$ be a quadratic symmetric operad.
If the associated shuffle operad $P^u$ has a quadratic convergent presentation, then $P$ is Koszul.
\end{theorem}
\begin{proof}
Let $X$ be a quadratic convergent $1$-polygraph presenting $P^u$. By definition, the $0$-generators in $X_0$ are concentrated in degree $1$ and the $1$-generators in $X_1$ in degree $2$. By construction, the polygraphic resolution $\Ov(X)$, constructed in Theorem \ref{T:MainTheoremD}, is concentrated on the superdiagonal, that is for $n\geq 2$, the $n$-generators in $\Ov(X)_n$ are of degree $n+1$; these cells are the generators of the $P$-bimodules of the resolution $P\langle\Ov(X)\rangle$ of $\Omega_P$ of Theorem~\ref{P:Resolution}. Thus, the Cartan-Eilenberg homology of $P$ is concentrated in degree $n+1$, and so the Quillen homology is concentrated on the diagonal. Following~\refeq{SSS:KoszulOperads}, we conclude that $P$ is Koszul.
\end{proof}

\subsubsection{Remark}
If we consider a quadratic symmetric operad whose generators are all of arity one, using Theorem~\ref{T:MainTheoremE} we recover the similar result for quadratic associative algebras: every algebra having a quadratic convergent presentation is Koszul, as proved in {\cite[Prop. 7.2.2]{GuiraudHoffbeckMalbos19}} by a polygraphic construction, see also {\cite[Sec. 4.3]{LodayVallette12}}, and \cite{Berger98} for a such a criterion with the rewriting rules ordered with respect to a monomial order.

\subsubsection{Koszul associative algebra without monomial order}
Let $A$ be the associative algebra presented by
\[
\left\langle \;w,x,y,z\; \left|\; w^2=wx, \; x^2=yx,\; y^2=yz, \;z^2=wz\;\right.\right\rangle.
\]
If we orient the relations according to a monomial order, say the order generated by $w<x<y<z$, this gives a $1$-polygraph with two critical branching that are non-confluent 
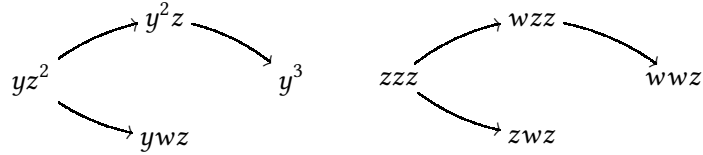
\[
\xymatrix @C=2.5em @R=0.5em{
& y^2z
\ar@/^1ex/[dr]
\\
yz^2
  \ar@/^1ex/ [ur]
  \ar@/_1ex/ [dr]
&&
y^3
\\
& ywz
}
\qquad
\xymatrix @C=2.5em @R=1em{
& wzz
\ar@/^1ex/[dr]
\\
zzz
  \ar@/^1ex/ [ur]
  \ar@/_1ex/ [dr]
&&
wwz
\\
& zwz
}
\]
Moreover, we show that any alphabetic order conduces to a similar situation of non-confluent critical branching.
Instead, consider the following $1$-polygraph:
\[
X:=\left\langle\;
w,x,y,z \;\left| \;
wx \fl w^2,\;
yx \fl x^2,\;
yz \fl y^2,\;
wz \fl z^2\;
\right.\right\rangle.
\]
The termination of $X$ is equivalent to the termination of the following $1$-polygraph
\[
\left\langle\;
w,x,y,z \;\left| \;
wx \fl w,\;
yx \fl x,\;
yz \fl y,\;
wz \fl z\;
\right.\right\rangle,
\]
and this second $1$-polygraph clearly terminates by considering the lengths of words, so $X$ terminates. Moreover, $X$ has no critical branchings, so it is confluent. Thus $X$ is a convergent quadratic $1$-polygraph, so by \cite[Prop. 7.2.2]{GuiraudHoffbeckMalbos19}, thus the algebra $A$ is Koszul.

\subsubsection{Koszul operad without monomial order}
Following the previous example, let $P$ be the symmetric operad presented by
\[
\left\langle\;
w,x,y,z\in P(2)\;
\left|\;\begin{array}{l}
w(1\,2)=w,\; x(1\,2)=x,\; y(1\,2)=y,\; z(1\,2)=z, \\
w\circ_{1,\id}w=w\circ_{1,\id}x,\;
x\circ_{1,\id}x=y\circ_{1,\id}x, \\
y\circ_{1,\id}y=y\circ_{1,\id}z,\;
z\circ_{1,\id}z=w\circ_{1,\id}z,
\end{array}
\;\right.\right\rangle.
\]
Consider the associated shuffle operad $P^u$, which is presented by
\[
\left\langle\;
w,x,y,z\in P^u(2)\;
\left|\;\begin{array}{l}
\inputtree{[w[w[1][2]][3]]} = \inputtree{[w[x[1][2]][3]]},\;
\inputtree{[w[w[1][3]][2]]} = \inputtree{[w[x[1][3]][2]]},\;
\inputtree{[w[1][w[2][3]]]} = \inputtree{[w[1][x[2][3]]]}, \\
\inputtree{[x[x[1][2]][3]]} = \inputtree{[y[x[1][2]][3]]},\;
\inputtree{[x[x[1][3]][2]]} = \inputtree{[y[x[1][3]][2]]},\;
\inputtree{[x[1][x[2][3]]]} = \inputtree{[y[1][x[2][3]]]}, \\
\inputtree{[y[y[1][2]][3]]} = \inputtree{[y[z[1][2]][3]]},\;
\inputtree{[y[y[1][3]][2]]} = \inputtree{[y[z[1][3]][2]]},\;
\inputtree{[y[1][y[2][3]]]} = \inputtree{[y[1][z[2][3]]]}, \\
\inputtree{[z[z[1][2]][3]]} = \inputtree{[w[z[1][2]][3]]},\;
\inputtree{[z[z[1][3]][2]]} = \inputtree{[w[z[1][3]][2]]},\;
\inputtree{[z[1][z[2][3]]]} = \inputtree{[w[1][z[2][3]]]}
\end{array}
\;\right.\right\rangle.
\]
If we orient the induced relations according to a monomial order, say an order where $w<x<y<z$, then in particular we get the rewriting rule
\[
\inputtree{[z[z[1][2]][3]]} \fl \inputtree{[w[z[1][2]][3]]}.
\]
Comparing the this rewriting rule with those of the previous example, we find that this rule creates a non-confluent critical pair. By Proposition \ref{P:ConvergentGröbner}, this also means that this presentation of $P$ does not admit a quadratic Gröbner basis.

Instead, if we orient every relation from right to left, we get a shuffle $1$-polygraph $X$ with $0$-generators $w,x,y,z\in X_0(2)$ and with twelve $1$-generators. With arguments similar to previous example, we show that $X$ is terminating. Moreover, there are no critical branchings, so by Theorem~\ref{T:CriticalBranchings} the $1$-polygraph $X$ is confluent. In this way, $X$ is a convergent quadratic $1$-polygraph, so by Theorem~\ref{T:MainTheoremE}, the operad $P$ is Koszul.

\begin{small}
\renewcommand{\refname}{\Large\textsc{References}}
\bibliographystyle{plain}
\bibliography{biblioCURRENT}
\end{small}

\clearpage

\quad

\vfill

\begin{footnotesize}

\bigskip
\auteur{Philippe Malbos}{malbos@math.univ-lyon1.fr}
{Universit\'e Claude Bernard Lyon 1\\
CNRS UMR 5208, Institut Camille Jordan\\
43 blvd. du 11 novembre 1918\\
F-69622 Villeurbanne cedex, France}

\bigskip
\auteur{Isaac Ren}{isaacren@kth.se}
{Department of Mathematics\\
KTH Royal Institute of Technology\\
S-10044 Stockholm, Sweden
}
\end{footnotesize}

\vspace{1.5cm}

\begin{small}---\;\;\today\;\;-\;\;\hhmm\;\;---\end{small}
\end{document}